\documentclass[12pt]{extarticle}
\usepackage[letterpaper, total={7in, 9in}]{geometry}

\usepackage{amsmath}
\usepackage{amsthm}
\usepackage{amssymb}
\usepackage{mathtools}
\usepackage[dvipsnames]{xcolor}
\usepackage{hyperref}
\usepackage{enumerate}

\def\I{\mathcal{I}}
\def\J{\mathcal{J}}
\def\L{\mathcal{L}}

\newtheorem{thm}{Theorem}[section]
\newtheorem{conj}[thm]{Conjecture}
\newtheorem{prop}[thm]{Proposition}
\newtheorem{claim}{Claim}
\newtheorem{lemma}[thm]{Lemma}
\newtheorem{cor}[thm]{Corollary}

\newtheorem{defi}[thm]{Definition}
\newtheorem{rem}[thm]{Remark}

\title{Maximal independent sets in the middle two layers of the Boolean lattice}
\author{J\'ozsef Balogh\thanks{Department of Mathematics, University of Illinois Urbana-Champaign, Urbana, IL, USA. Research supported in part by NSF grants RTG DMS-1937241, FRG DMS-2152488, the Arnold O. Beckman Research Award (UIUC Campus Research Board RB 24012), and the Simons Fellowship. Chen and Garcia are partially supported by NSF grant RTG DMS-1937241 and UIUC Campus Research Board RB 24012.\\
\text{ }{\hskip1.1em\relax} Email: \texttt{\{jobal, cechen4, rig2\}@illinois.edu}.} , Ce Chen$^{*}$, Ramon I. Garcia$^{*}$}
\date{\today}

\begin{document}

\maketitle

\begin{abstract}
    Let $B(2d-1, d)$ be the subgraph of the hypercube $\mathcal{Q}_{2d-1}$ induced by its two largest layers. Duffus, Frankl and R\"odl proposed the problem of finding the asymptotics for the logarithm of the  number of maximal independent sets in $B(2d-1, d)$. Ilinca and Kahn determined the logarithmic asymptotics and reiterated the question of what their order of magnitude is. We show that the number of maximal independent sets in $B(2d-1,d)$ is
    \[
\left(1+o(1)\right)(2d-1)\exp\left(\frac{(d-1)^2}{2^{2d-1}}\binom{2d-2}{d-1}\right)\cdot 2^{\binom{2d-2}{d-1}},
    \]
    and describe their typical structure. The proof uses a new variation of Sapozhenko's Graph Container Lemma, a new isoperimetric lemma,  a theorem of  Hujter and Tuza on the number of maximal independent sets in triangle-free graphs and a stability version of their result by Kahn and Park, among other tools. 
\end{abstract}

\section{Introduction}

In a graph $G$, an \emph{independent set} $I$ is a vertex set $I\subseteq V(G)$ such that there is no edge $xy\in E(G)$ with $x,y\in I$. Enumerating independent sets and determining their typical structure has become a well-researched topic in extremal combinatorics. One of the fundamental results in this direction is due to Kleitman~\cite{kleitmanantichain}, who determined the log-asymptotics of the number of antichains in $\{0,1\}^n$.

Many classical theorems in extremal combinatorics have been extended to enumerative and structural results in the past decades. In particular, the \emph{hypergraph container method} of Balogh--Morris--Samotij~\cite{balogh2015independent} and Saxton--Thomason~\cite{saxton2015hypergraph} has seen particular success with such problems.

Denote by $\mathcal{Q}_n$ the \textit{discrete hypercube} of dimension $n$, which is the graph defined on the family of subsets of $[n]=\{1,\ldots,n\}$, where two sets are adjacent if and only if they differ in exactly one element. Independent sets in the discrete hypercube have received a lot of attention in recent decades, see~\cite{2025kang, galvin2011threshold,jenssen2024refined, jenssen2019independent}.

A classical result in the discrete hypercube is due to Korshunov and Sapozhenko~\cite{ korshunov1983number},  who showed that the number of independent sets in the Boolean lattice $\mathcal{Q}_{n}$ is  $(1+o(1))\cdot2\sqrt{e}\cdot2^{2^{n-1}}$. Sapozhenko~\cite{sapozhenko1987} developed an influential method to obtain this asymptotics, and this tool is now known as Sapozhenko's Graph Container Method.  

Sapozhenko's Graph Container Method has been applied successfully in several settings. Some remarkable applications of it include enumeration of Lipschitz functions on weak expanders \cite{krueger2024lipschitz}, of independent sets in Abelian Cayley graphs \cite{potukuchi2021abelian}, and of graph homomorphisms from the discrete torus \cite{Jenssen2023Homorphisms}. 


An independent set $I\subseteq V(G)$ is \textit{maximal} if for every $v\in V(G)\setminus I$, the vertex set $I\cup\{v\}$ is not an independent set. Let $\mathcal{I}(G)\coloneqq\textup{MIS}(G)$ be the set system consisting of  the maximal independent sets in $G$ and $\textup{mis}(G)\coloneqq|\textup{MIS}(G)|$. Recently, Kahn and Park~\cite{Kahn2022} obtained precise asymptotics for $\textup{mis}(\mathcal{Q}_n)$.

\begin{thm}[Kahn-Park~\cite{Kahn2022}]\label{thm:Kahn-Park mis}
    \[
    \textup{mis}(\mathcal{Q}_n)=(1+o(1))\cdot 2n\cdot 2^{2^{n-2}}.
    \]
\end{thm}

Our goal is to prove an analogue result for the graph spanned by the middle two layers of the hypercube, which is harder, as this graph has weaker isoperimetric properties. To overcome this difficulty, we introduce a modified version of Sapozhenko's container method, along with vertex and edge isoperimetry results for the graphs induced by two consecutive layers in the Boolean lattice.

The proof of Theorem~\ref{thm:Kahn-Park mis} relies on the particular structure of $\mathcal{Q}_{n}$. The absence of some of the properties of $\mathcal{Q}_{n}$ in the middle layers increases the difficulty of asymptotic enumeration of maximal independent sets. For example, the precise asymptotics for the number of $q$-colorings for the middle two layers have been determined only for even values of $q$, see \cite{li2025number} for further details, while for $\mathcal{Q}_{n}$ the number of proper $q$-colorings has been estimated for all $q$, see~\cite{Jenssen2023Homorphisms}. 

For arbitrary positive integers $n$ and $k\leq n$, let $\L_{k}^{n}\coloneqq\{x\subseteq [n]: |x|=k\}$ and $B(n,k)$ be the bipartite graph with vertex set $V(B(n,k))\coloneqq\L_{k}^{n}\cup \L_{k-1}^{n}$ such that $xy\in E(B(n,k))$ if and only if $x\subsetneq y$. For a fixed positive integer $d$, the graph $B(2d-1, d)$ is a $d$-regular bipartite graph. We fix $n=2d-1$, $k=d$ and write $\mathcal{L}_{d}$ and $\mathcal{L}_{d-1}$ instead of $\mathcal{L}_{d}^{2d-1}$ and $\mathcal{L}_{d-1}^{2d-1}$ for ease of notation. We are interested in determining the precise asymptotics of $\textup{mis}(B(2d-1,d))$. We will denote $\mathcal{I}(B(2d-1,d))$ by $\mathcal{I}$. 

Duffus, Frankl and R\"odl~\cite{duffus2011maximal} started the study of maximal independent sets in $B(n,k)$, showing $\binom{n-1}{k-1}\leq \log_{2} \textup{mis}(B(n,k))\leq (1.3563+o(1))\binom{n-1}{k-1}$ for $1\leq k\leq n$ as $n\rightarrow \infty $. Ilinca and Kahn~\cite{ilinca2013counting} obtained the logarithmic asymptotics for $\textup{mis}(B(n,k))$. 

\begin{thm}[Ilinca-Kahn~\cite{ilinca2013counting}]
    $$\log_2 \textup{mis}(B(n,k))=(1+o(1))\binom{n-1}{k-1}.$$
\end{thm}

Ilinca and Kahn~\cite{ilinca2013counting} also raised the question of finding the precise assymptotics of $\textup{mis}(B(n,k))$, and proposed the following conjecture. 
\begin{conj}[Ilinca-Kahn~\cite{ilinca2013counting}]\label{conjIK}
    \[
    \textup{mis}(B(n,k))=(1+o(1))\cdot n\cdot 2^{\binom{n-1}{k-1}}.
    \]
\end{conj}

They~\cite{ilinca2013counting} wrote ``This does not look easy, but we can at least guess what the truth
should be.'' Somewhat surprisingly, Balogh, Treglown and Wagner~\cite{BTW} disproved Conjecture~\ref{conjIK} by providing a construction showing that there are more maximal independent sets in $B(n,k)$ than those conjectured by Ilinca and Kahn~\cite{ilinca2013counting} for certain values of $n$ and $k$.  

\begin{thm}[Balogh-Treglown-Wagner~\cite{BTW}]
    If $|n/2-k|\leq \sqrt{n}$, then 
    \[
    \log_2 \textup{mis}(B(n,k))\geq \binom{n-1}{k-1}+\Omega\left(n^{3/2}\right).
    \]
\end{thm}
The authors of~\cite{BTW} conjectured only after publishing their paper that their construction (in soul) indeed provides the ``typical structure" for most of the maximal independent sets in $B(2d-1,d)$.

For an integer $k \in [2d-1]$, define the edge set $M_{k}\coloneqq \{xy\in E(B(2d-1,d)): k\in x\in \L_d, y=x\setminus\{k\}\}$. Note that $M_{k}$ is an induced matching in $B(2d-1, d)$ for every $k$, and we call it a \textit{canonical matching}. Let $G$ be a simple graph. A set $A\subseteq V(G)$ is \textit{$k$-linked} if $A$ is connected in $G^k$, where $G^k$ is the simple graph defined on $V(G)$, in which two vertices are adjacent if their distance in $G$ is at most $k$. A \textit{$k$-linked component} of $G$ is a maximal $k$-linked subset of $V(G)$. Our main result shows that the conjecture from~\cite{BTW} is true. 

\begin{thm}\label{thm:TypicalStructure}
Almost all maximal independent sets $I\in \textup{MIS}(B(2d-1,d))$ have the following property\footnote{That is, the proportion of maximal independent sets that do not have this property goes to zero as $d \rightarrow \infty$.}. There is a $k\in [2d-1]$ such that every $3$-linked component of $I\setminus V(M_k)$ has size $1$ or $2$.
\end{thm}

From this we conclude the following.

\begin{thm}\label{thm:main}
\begin{equation}\label{eq:mainthm}
    \textup{mis}(B(2d-1,d))=(1+o(1))(2d-1)\exp\left(\frac{(d-1)^{2}}{2^{2d-1}}\binom{2d-2}{d-1}\right)\cdot 2^{\binom{2d-2}{d-1}}.
\end{equation}
\end{thm}
By Stirling's approximation, it follows that $\frac{(d-1)^2}{2^{2d-1}}\binom{2d-2}{d-1}=\frac{(d-1)^{3/2}}{2\sqrt{\pi}}-\frac{(d-1)^{1/2}}{24\sqrt{\pi}}+O(d^{-1/2})$. Therefore, the right hand side of~\eqref{eq:mainthm} is \[
\left(1+o(1)\right)\cdot2d\cdot\exp\left(\frac{(d-1)^{3/2}}{2\sqrt{\pi}}-\frac{(d-1)^{1/2}}{24\sqrt{\pi}}\right)\cdot2^{\binom{2d-2}{d-1}}.\] 

To contextualize Theorem~\ref{thm:TypicalStructure}, Kahn and Park \cite{Kahn2022} introduced a notion for canonical matching in the graph $\mathcal{Q}_{n}$ similarly. They proved that for almost all maximal independent sets $I\in\textup{MIS}(\mathcal{Q}_{n})$, there is a canonical matching $M$ such that $I\setminus V(M)=\emptyset$. Conjecture~\ref{conjIK} was formulated under the analogous assumption that, for almost all maximal independent sets $I\in\textup{MIS}(B(2d-1,d))$, there exists a canonical matching $M$ such that $I\setminus V(M)=\emptyset$. However, this assumption was later shown to be false.

The structure for the rest of the paper is the following. In Section~\ref{sec:tools} we introduce the main tools we use. In particular, to each maximal independent set $I$, we assign an induced matching $M(I)$. In Section~\ref{sec:lowerbound}, we provide a lower bound on $\textup{mis}(B(2d-1,d))$ by using the construction from~\cite{BTW}. In Section~\ref{sec:containers}, we prove that most maximal independent sets $I$ get assigned a large induced matching $M(I)$. For doing so, we use an application of the standard graph container algorithm first, and then a variant of Sapozhenko's graph container algorithm. In Section~\ref{sec:stability}, by using edge isoperimetry, we prove that if an induced matching is large, then most of its edges use the same direction. In Section~\ref{sec:counting}, we complete the proof of Theorem~\ref{thm:main} and prove Theorem~\ref{thm:TypicalStructure}.


\section{Tools}\label{sec:tools}

We use $\log$ for $\log_2$. For a set $A\subseteq V(G)$ let $d_{A}(v)\coloneqq|N(v)\cap A|$ and $N(A)\coloneqq\cup_{v\in A}N(v)$. We write $\binom{n}{\leq b}$ for $\sum_{i\leq b}\binom{n}{i}$.
We use ``cost of $X$" for $\log_2$ of the number of choices for $X$. Floor and ceiling signs will be omitted if they are not essential in the computation.

\begin{prop}\label{prop:Entropy}
    For every $\alpha\in [0,1/2]$ and $n\in \mathbb{Z}^{+}$,
    \[
    \sum_{i\leq \alpha n}\binom{n}{i}\leq 2^{H(\alpha) n},
    \]
    where $H(p)=-p\log p-(1-p)\log (1-p)$ is the binary entropy function. 
\end{prop}

\begin{prop}\label{prop:UpperBoundZ}
    Let $0\leq d'\leq d$ and $Z\subseteq W\subseteq V(B(2d-1,d))$ be such that $d_{Z}(v)\leq d'$ for every $v\in Z$. Let $L\coloneqq e(W,W^{c})$, then
    \[
    |Z|\leq \frac{d|W|+L}{2d-d'}.
    \]
\end{prop}
\begin{proof}
Simple double-counting implies that
\[
d|W\setminus Z|\geq e(Z,W\setminus Z)\geq (d-d') |Z|-L,
\]
which gives
\[
d|W|+L\geq (2d-d')|Z|.
\]
\end{proof}

To obtain an upper bound on the number of maximal independent sets, our starting point is the following theorem.

\begin{thm}[Hujter-Tuza~\cite{HujterTuza}]\label{thm:HujterTuza}
    For every $m$-vertex, triangle-free graph $G$,
    \[
\log \textup{mis}(G)\leq m/2,
    \]
    with equality if and only if $G$ is a perfect matching.
\end{thm}

We use a stability version of Theorem~\ref{thm:HujterTuza}, proved by Kahn and Park~\cite{kahn2020stability}. Before stating it, we need a couple definitions.

\begin{defi}\label{def:M(I)}
    An \textit{induced matching} $M$ in a graph $G$ is a set of disjoint edges $e\in E(G)$ such that for every pair of edges $xy,uv\in M$, we have $xu, xv, yu, yv\notin E(G)$. For tiebreaking purposes only, we fix an arbitrary linear ordering $\prec$ on the set of induced matchings in $G$ such that $M\prec M'$ whenever $|M|>|M'|$.

Given a graph $G$ and an independent set $I\subseteq V(G)$, let $M_{G}(I)$ be the induced matching of maximum size such that if $e=xy\in M_{G}(I)$ then $|\{x,y\}\cap I|=1$ and if $v\in I\setminus V(M_{G}(I))$ then $N(v)\cap V(M_{G}(I))=\emptyset$. If there is more than one such induced matching, then we choose the first one under $\prec$. We consider the empty set to be an induced matching, so $M_{G}(I)$ always exists.

For a set $Z\subseteq \mathcal{L}_{d}\cup \mathcal{L}_{d-1}$, we consider the restriction of $\prec$ to the collection of induced matchings in $B(2d-1,d)[Z]$ and for $I\in\mathcal{I}(B(2d-1,d)[Z])$, define $M_{Z}(I)\coloneqq M_{B(2d-1,d)[Z]}(I)$.
\end{defi}

\begin{thm}[Kahn-Park~\cite{kahn2020stability}, Theorem 14]\label{thm:Hujter-TuzaStab}
    There is a constant $c>0$ such that for every $\varepsilon>0$ and $m$-vertex, triangle-free graph $G$,
    \[
\log |\{I\in \mathcal{I}(G):|M_G(I)|<(1-\varepsilon)m/2\}|< (1-c\varepsilon)m/2.
    \]
\end{thm}

Notice that applying Theorem~\ref{thm:HujterTuza} directly to $B(2d-1,d)$ only gives $\textup{mis}(B(2d-1,d))\leq 2^{\binom{2d-1}{d}}$, which is very far from what we want. Hence, we need to reduce the size of the graph we are working with.

First, we  show that most maximal independent sets $I\in \mathcal{I}$ get assigned a large $M(I)$, where $M(I)\coloneqq M_{B(2d-1, d)}(I)$ for ease of notation. This motivates us to define the set
\[
\J_1\coloneqq\left\{I\in \I: |M(I)|> \left(1-\frac{2\log^3 d}{d}\right)\binom{2d-2}{d-1}\right\}.
\]
Applying the basic graph container algorithm, we will obtain the following result.
\begin{thm}\label{thm:FirstPhase}
    \[
    |\mathcal{I}\setminus \J_1|= 2^{\left(1-\Omega\left(\frac{\log^3 d}{d}\right)\right)\binom{2d-2}{d-1}}=o(|\mathcal{I}|).
    \]
\end{thm}

\noindent For $I\in\J_1$, the size of the matching $M(I)$ is not guaranteed to be large enough for our purposes. To be able to show that most matchings $M(I)$ are not too far away from a canonical matching, we run a second container argument, which shows that most maximal independent sets $I$ satisfy $|M(I)|\geq (1-\beta)\binom{2d-2}{d-1}$ for some $\beta\ll 1/d$. This motivates the following definition.

\begin{equation}\label{defi:J2}
\J_2\coloneqq\left\{I\in \I: |M(I)|> \left(1-\frac{2\log^5 d}{d^{3/2}}\right)\binom{2d-2}{d-1}\right\}.
\end{equation}
Using Sapozhenko's graph container method, we will show the following.
\begin{thm}\label{thm:SecondPhase}
    $$|\mathcal{J}_1\setminus \J_2|= 2^{\left(1-\Omega\left(\frac{\log^5 d}{d^{3/2}}\right)\right)\binom{2d-2}{d-1}}=o(|\mathcal{I}|).$$ 
\end{thm}

Using edge isoperimetry, we will show that if $M(I)$ is large enough, then there is a unique positive integer $k$ such that $|M(I)\cap M_{k}|$ is large. 
This will be helpful to estimate $|\J_2|$.
\begin{thm}\label{thm:stability}
    For every $I\in \J_2$ there exists a canonical matching $M_{k}$ such that 
    \[
    |M(I)\cap M_{k}|\geq \left(1-\frac{25 \log^5 d}{\sqrt{d}}\right)\binom{2d-2}{d-1}.
    \]
\end{thm}


\subsection{Isoperimetry tools}

We will use the following  version of the Kruskal-Katona Theorem~\cite{katona2009theorem, kruskal1963number}.

\begin{thm}[Lov{\' a}sz~\cite{lovasz2007combinatorial}]\label{thm:lovasz}
Let $\mathcal{A}$ be a family of $m$-element subsets of a fixed set U and $\mathcal{B}$ be the family of all $(m-q)$-element subsets of the sets in $\mathcal{A}$. If $|\mathcal{A}|=\binom{x}{m}$ for some real number $x$, then
$|\mathcal{B}|\geq \binom{x}{m-q}$.
\end{thm}

We use an isoperimetry result similar to the one used in \cite{balogh2021independent}, with a small modification of $d^4$ to $d^6$ in (ii).

\begin{lemma}\label{lem:isoperimetry}
There is a $d_0>0$ such that for every integer $d>d_0$ and vertex set $A\subseteq \mathcal{L}_{d}$:
\begin{itemize}
\itemsep-0.3em
\item[\rm (i)] If $|A|\leq d$, then $|N(A)|\geq d|A| - |A|^2/2$. In particular, if $|A|\leq d/4$ then $|N(A)|\geq 7d|A|/8$.
\item[\rm (ii)] If $|A|\leq d^6$, then $|N(A)|\geq d|A|/9$.
\item[\rm (iii)] If $|A|\leq \binom{2d-2}{d}$, then $|N(A)|\geq \left(1+\frac{1}{d-1}\right)|A|$.
\end{itemize}
\end{lemma}
\begin{proof}
    For (i), notice that any two vertices have at most one common neighbor in $B(2d-1,d)$ and therefore $|N(A)|\geq d|A|-|A|^2/2$ for every $A\subseteq \mathcal{L}_{d}\cup \mathcal{L}_{d-1}$.

    For (ii), let $x$ be a real number such that $|A|=\binom{x}{d}$. Suppose for a contradiction that $\lfloor x\rfloor \geq d+7$, then $\binom{x}{d}\geq \binom{\lfloor x\rfloor }{d}\geq \binom{d+7}{d}\geq \frac{d^{7}}{7!}$, which is a contradiction with $\binom{x}{d}\leq d^{6}$ for sufficiently large $d$. Therefore, we can assume $x\leq d+8$. By Theorem~\ref{thm:lovasz} we obtain
    
    \[
    |N(A)|\geq \binom{x}{d-1}=\frac{d}{x-d+1}\binom{x}{d}=\frac{d}{x-d+1}|A|\geq \frac{d}{9}|A|.
    \]
    For (iii), let $x$ be a real number such that $|A|=\binom{x}{d}$, with $x\leq 2d-2$. By Theorem~\ref{thm:lovasz} we have      
    \[
    |N(A)|\geq \binom{x}{d-1}=\frac{d}{x-d+1}\binom{x}{d}=\frac{d}{x-d+1}|A|\geq \left(1+\frac{1}{d-1}\right)|A|.
    \]

\end{proof}

\begin{lemma}\label{lem:Iso3}
    Let $A\subseteq \mathcal{L}_{d-1}^{2d-2}$ such that $|A|\leq \frac{\log^6 d}{\sqrt{d}}\binom{2d-2}{d-1}$. Then, $|N_{B(2d-2,d-1)}(A)|\geq \left(1+\frac{\log d}{5d}\right)|A|$.
\end{lemma}

\begin{proof}
Let $x$ be a real number such that $|A|=\binom{2d-2-x}{d-1}$. Since $|A|\leq \frac{\log^6 d}{\sqrt{d}}\binom{2d-2}{d-1}$, it follows that
    \[
    \binom{2d-2-x}{d-1}\leq \frac{\log^6 d}{\sqrt{d}}\binom{2d-2}{d-1},
    \]
which is equivalent to
    \[
    \prod_{i=0}^{d-2}(2d-2-x-i)\leq \frac{\log^6 d}{\sqrt{d}}\prod_{i=0}^{d-2} (2d-2-i).
    \]
Therefore, we have
    \[
    \left(1-\frac{x}{d-1}\right)^{d-1}\leq\prod_{i=0}^{d-2}\left(1-\frac{x}{2d-2-i}\right)\leq d^{-1/2}\cdot  \log^6 d\leq d^{-1/3}.
    \]
Since $d$ is sufficiently large, we have $\exp\left(-x\right)\approx\left(1-\frac{x}{d-1}\right)^{d-1}\leq\exp\left(-\frac{\ln d}{3}\right)$, thus $x\geq (\ln d)/3$.
Therefore, by Theorem~\ref{thm:lovasz} we have      
    \begin{equation}
        \begin{aligned}
             |N_{B(2d-2,d-1)}(A)|&\geq \binom{2d-2-x}{d-2}=\frac{d-1}{2d-2-x-d+2}\binom{2d-2-x}{d-1}=\frac{(d-1)|A|}{d-x}\\
             &\geq \frac{(d-1)|A|}{d-(\ln d)/3}=\left(1+\frac{(\ln d)/3-1}{d-(\ln d)/3}\right)|A|\geq \left(1+\frac{\log d}{5d}\right)|A|,
        \end{aligned}
    \end{equation}
    where the last inequality holds, since $d$ is sufficiently large.
\end{proof}

Our induced matchings will be typically big, which will imply that most edges in $B(2d-1,d)$ will have at least one endpoint in $V(M(I))$. With this in mind, we use edge isoperimetry inequalities. 
We will need the following theorem from~\cite{IsoBey} on the sum of squares of degrees in a hypergraph.

\begin{thm}[Bey~\cite{IsoBey}, Theorem 1]\label{thm:edgeIso-bey}
    Let $G$ be a $k$-uniform hypergraph on $n$ vertices and $\ell$ be an integer with $0\leq \ell\leq k$, then \[
    \sum_{v\in \binom{V(G)}{\ell}}d_{G}(v)^2\leq \frac{\binom{k}{\ell}\binom{k-1}{\ell}}{\binom{n-1}{\ell}}|E(G)|^{2}+\binom{k-1}{\ell-1}\binom{n-\ell-1}{k-\ell}|E(G)|,
    \]
    where $d_G(v)$ is the number of hyperedges $e\in E(G)$ such that $v\subseteq e$.
\end{thm}

\noindent We instantly obtain the following by applying Theorem~\ref{thm:edgeIso-bey} with $k=d, \ell=d-1$ and $G=A\subseteq \mathcal{L}_{d}^{n}$.

\begin{thm}\label{thm:edgeIso}
For every $A\subseteq \mathcal{L}_{d}^{n}$ we have
\[
\sum_{v\in \L_{d-1}^{n}}d_{A}(v)^2\leq \frac{d}{\binom{n-1}{d-1}}|A|^{2}+(n-d)(d-1)|A|.
\]
\end{thm}

We call an ordered triplet of vertices $(x,y,z)\in (\mathcal{L}_{d}^{n}\cup \mathcal{L}_{d-1}^{n})^{3}$ an \textit{adjacent triplet} if $xy,yz\in E(B(n,d))$. Notice that in this definition, it is allowed that $x=z$. For $A\subseteq \mathcal{L}^n_{d}$, let $\mathcal{B}(A)$ be the set of adjacent triplets $(x,y,z)$ such that $x\in A$ and $z\notin A$.

\begin{lemma}\label{lem:lowerBadTriplets}
    Let $A\subseteq \mathcal{L}_{d}^{n}$, then
    \[
    |\mathcal{B}(A)|\geq n|A|-d\frac{|A|^2}{\binom{n-1}{d-1}}.
    \]
\end{lemma}

\begin{proof} 
   Notice that $\sum_{v\in \mathcal{L}_{d-1}^{n}}d_{A}(v)^2$ counts the number of adjacent triplets $(x,v,z)$ with $x,z\in A$, where $x$ and $z$ are potentially the same vertex. On the other hand, we can construct such a triplet by choosing a vertex $u\in A$, a vertex $v\in N(u)$, and a vertex $w\in N(v)$ but we will have to exclude the triplets in $\mathcal{B}(A)$. This double-counting argument gives
\begin{equation}\label{badtripletintro}
\sum_{v\in \mathcal{L}_{d-1}^{n}} d_{A}(v)^2= (n-d+1)d|A|-|\mathcal{B}(A)|.
\end{equation}

\noindent Using Theorem~\ref{thm:edgeIso} we get
\[
\frac{d}{\binom{n-1}{d-1}}|A|^{2}+(n-d)(d-1)|A|\geq \sum_{v\in \mathcal{L}_{d-1}} d_{A}(v)^2 = d(n-d+1)|A|-|\mathcal{B}(A)|,
\]
which implies our desired inequality.
\end{proof}

\begin{lemma}\label{lem:edgeIso}
    Let $A\subseteq \mathcal{L}_{d}^{n}$ and $M$ be an induced matching in $B(n,d)$. Set $A'\coloneqq N(A)\cap V(M), H\coloneqq N(A)\setminus A'$ and suppose that for every $x\in A\setminus V(M)$, $x$ has no neighbors in $V(M)$. Then,
    \begin{equation}\label{ineq:edgeIso}
    e(H, \mathcal{L}_{d}^{n}\setminus A)\geq\frac{d}{n-d}\left( |A|-\frac{|A|^2}{\binom{n-1}{d-1}}\right).
    \end{equation}
\end{lemma}

\begin{proof}
    If $|A|>\binom{n-1}{d-1}$, then the right-hand side of~\eqref{ineq:edgeIso} is negative, so we may assume $|A|\leq \binom{n-1}{d-1}$. For every vertex $x\in A\cap V(M)$, since $M$ is an induced matching, there is a unique $y\in A'\cap N(x)$. The vertex $y$ has $n-d$ neighbors outside $A$. Therefore, the number of adjacent triplets $(x,y,z)$ with $x\in A$, $y\in A'$ and $z\notin A$ is $(n-d)|A\cap V(M)|\leq (n-d)|A|$.
    
    Let $\mathcal{B}_{0}(A)$ be the set of adjacent triplets $(x,y,z)$ such that $x\in A$, $y\notin A'$ and $z\notin A$. We have that $|\mathcal{B}_{0}(A)|\geq |\mathcal{B}(A)|-(n-d)|A|$. Using Lemma~\ref{lem:lowerBadTriplets} we get
    \[
    |\mathcal{B}_{0}(A)|\geq n|A|-d\frac{|A|^2}{\binom{n-1}{d-1}}-(n-d)|A|= d|A|-d\frac{|A|^2}{\binom{n-1}{d-1}}.
    \]
    Each triplet $(x,y,z)\in\mathcal{B}_{0}(A)$ has $yz\in E(H,\mathcal{L}_{d}^{n}\setminus A)$. Each edge $e\in E(H,\mathcal{L}_d^{n}\setminus A)$ appears in at most $n-d$ adjacent triplets $(x,y,z)$ of $\mathcal{B}_{0}(A)$ as the edge $yz$. We obtain 
    \[
    (n-d)\cdot e(H,\mathcal{L}_{d}^{n}\setminus A)\geq |\mathcal{B}_{0}(A)|\geq d|A|-d\frac{|A|^2}{\binom{n-1}{d-1}}.\qedhere
    \]
    
\end{proof}

\subsection{Sapozhenko containers tools}\label{subsec: Sapozhenko tools}

\begin{lemma}[Sapozhenko~\cite{sapozhenko1987}]\label{lem:linked}
    If $A$ is $m$-linked, and $T\subseteq V(G)$ is such that $d(u,T)\leq r$ for each $u\in A$ and $d(u, A)\leq r$ for each $u\in T$, then $T$ is $(m+2r)$-linked.
\end{lemma}

\begin{lemma}[Folklore, see e.g.~\cite{galvin2011threshold}]\label{lem:numcomp}
Let $\Sigma=(V, E)$ be a $d$-regular graph with $d\geq 2$. The number of $k$-linked subsets of $V$ of size $t$ which contain a fixed vertex $v$ is at most $d^{3kt}$.
\end{lemma}

Recall that our graph $B(2d-1, d)$ is $d$-regular. As the number of choices for $v$ is $\binom{2d-1}{d}$  in $\L_d$, we have the following.

\begin{cor}\label{cor:klinkednumber}
The number of choices for $k$-linked $T\subseteq \mathcal{L}_d$ of size $|T|=t$ is at most $\binom{2d-1}{d}d^{3kt}$.
\end{cor}

The following is a special case of a fundamental result due to Lov\'asz~\cite{lovasz1975ratio} and Stein~\cite{stein1974two}.

\begin{thm}\label{thm:cover}
   Let $\Gamma$ be a bipartite graph with parts $P\cup Q$. If $|N(p)|\geq a$ for every $p\in P$ and $|N(q)|\leq b$ for every $q\in Q$, then there is a set $Q'\subseteq Q$ such that every $p\in P$ has a neighbor in $Q'$ and 
   $$|Q'|\leq \frac{|Q|}{a}\left(1+\ln b\right).$$
\end{thm}

\begin{prop}\label{prop:decom}
The  number  of  compositions  of $n$, which is the number of ways to express $n$ as the sum of an ordered sequence of positive integers, is $2^{n-1}$ and  the  number of compositions with  at  most $b$ parts is $\sum_{i<b}\binom{n-1}{i}<2^{b\log(en/b)}$, when $b< n/2$.
\end{prop}

\subsection{A new isoperimetric lemma}

For a $d$-uniform set system $\mathcal{F}\subseteq \binom{[2d-1]}{d}$, the \textit{shadow} $\partial \mathcal{F}$ is
\[
\partial \mathcal{F}\coloneqq \{E: E=F\setminus \{x\} \text{ for some } F\in \mathcal{F} \text{ and } x\in F\}.
\]
For a $d$-uniform set system $\mathcal{F}\subseteq \binom{[2d-1]}{d}$ and an integer $i\in [2d-1]\setminus\{1\}$, the \textit{shifting} function $S_i:\binom{[2d-1]}{d}\rightarrow \binom{[2d-1]}{d}$ is
\[
S_i(F)\coloneqq\begin{cases}
    F\setminus \{i\}\cup \{1\} & \quad \quad \quad \text{ if } i\in F,\text{ } 1\notin F \text{ and } F\setminus \{i\}\cup \{1\}\notin \mathcal{F},\\
    F & \quad \quad \quad \text{ otherwise.}
\end{cases}
\]
Set $\mathcal{S}_{i}(\mathcal{F})\coloneqq\{S_i(F): F\in \mathcal{F}\}$. The following properties of the shifting function are useful.

\begin{prop}[Das~\cite{kk-das}]\label{shiftprop}
(i) $|\mathcal{S}_i(\mathcal{F})|=|\mathcal{F}|$.\\
    (ii)  $\partial \mathcal{S}_i (\mathcal{F})\subseteq \mathcal{S}_i(\partial\mathcal{F}).$
\end{prop}

\noindent In particular, we have
\begin{equation}\label{shadow}
|\partial \mathcal{S}_i (\mathcal{F})|\leq |\mathcal{ S}_i(\partial\mathcal{F})|=|\partial \mathcal{F}|.    
\end{equation}

\begin{lemma}\label{lem:lastIso}
    Let $M$ be an induced matching in $B(2d-1,d)$ with $|M|\leq \frac{\log^6 d}{\sqrt{d}}\binom{2d-2}{d-1}$. Set $B\coloneqq V(M)\cap \mathcal{L}_{d}$, then
    \begin{equation*}
        |N(B)|\geq \left(2+\frac{\log d}{5d}\right)|B|.
    \end{equation*}
\end{lemma}

\begin{proof}
    We apply the shifting function $\mathcal{S}_i$ one by one for every $2\leq i\leq 2d-1$ to $B$. Denote $B'\coloneqq \mathcal{S}_{2d-1}\circ \ldots \circ \mathcal{S}_{2}(B)$ the final set.

    \begin{claim}\label{clm:contain1} For every $v'\in B'$, we have $1\in v'$.
    \end{claim}

    \begin{proof}
    If every $v\in B$ contains 1, then the claim trivially holds. Otherwise, there is some $v$ satisfying that $1\notin v\in B$, $k\in v$, and the edge formed by $v$ and $v\setminus \{k\}$ is in $M$ for some $k$. Then, $v\setminus \{k\} \cup \{1\}\notin B$ since $M$ is an induced matching. By the definition of the shifting function, if $S_k(v)\neq v\setminus \{k\}\cup\{1\}$, then either it is already shifted, or there exists some $j<k$ and $j\in w\in B$ such that $S_j(w)=w\setminus\{j\}\cup \{1\}=v\setminus \{k\} \cup \{1\}$. In the latter case, we have $w\setminus\{j\}=v\setminus \{k\}$, which implies $w\in N(v\setminus \{k\})$, a contradiction to that $M$ is an induced matching. Therefore, if $v$ is not shifted by $S_k$, then it is already shifted by some $S_j$ with $j<k$, i.e.,~there exists some $\ell\leq k$ such that $S_\ell(v)=v\setminus \{\ell\}\cup \{1\}=v'$, where $v'$ denotes the corresponding vertex of $v$ in $B'$.
    \end{proof}

    Let $B_1\subseteq \partial B'\subseteq \L_{d-1}$ consist of vertices $u$ with $u\cup \{1\}\in B'$. By Claim~\ref{clm:contain1}, we have $|B_1|=|B'|$. By Proposition~\ref{shiftprop} (i), we have \[|B_1|=|B'|=|B|=|M|\leq \frac{\log^6 d}{\sqrt{d}}\binom{2d-2}{d-1}.\]
    Observing that no element of $B_1$ contains $1$, hence $B_1$ could be viewed as a $(d-1)$-uniform set system on $2d-2$ elements. Applying Lemma~\ref{lem:Iso3} to $B_1$, we have \[
    |\partial B_1|\geq \left(1+\frac{\log d}{5d}\right)|B_1|.
    \]
    Let $B_2\coloneqq \partial B'\setminus B_1$. Note that for every $u\in \partial B_1$, there is an integer $2\leq i\leq 2d-1$ such that $u\cup \{i\}\in B_1$. Thus, $u\cup \{1,i\}\in B'$, which implies $u\cup \{1\}\in B_2$. Therefore, we have $|B_2|\geq |\partial B_1|.$

    By~\eqref{shadow}, we conclude that 
    \[
    |N(B)|=|\partial B|\geq |\partial B'|=|B_1|+|B_2|\geq \left(2+\frac{\log d}{5d}\right)|B_1|=\left(2+\frac{\log d}{5d}\right)|B|.\qedhere
    \]
\end{proof}



\section{Lower bound}\label{sec:lowerbound}
In this section, we provide a lower bound for $\textup{mis}(B(2d-1,d))$ by describing a procedure that generates maximal independent sets using the construction from \cite{BTW}. Note that~\cite{BTW} has not provided detailed description. The subsequent sections will be dedicated to show that this lower bound is tight.

\begin{prop}\label{prop:LowerBoundMIS}
\vglue -0.2cm
\begin{equation*}
    \textup{mis}(B(2d-1,d))\geq (2d-1)(1+o(1))\exp\left(\frac{(d-1)^{2}}{2^{2d-1}}\binom{2d-2}{d-1}\right)\cdot 2^{\binom{2d-2}{d-1}}.
\end{equation*}
\end{prop}
\begin{proof}

Recall that for $k\in [2d-1]$, the induced matching in direction $k$ is
\[
M_{k}=\{uv: k\in u\in \mathcal{L}_d, v=u\setminus \{k\}\}.
\]
To build a maximal independent set, fix an integer $k\in[2d-1]$. For each edge $e\in M_k$ we choose one of its endpoints. Then, there is a unique maximal independent set containing the vertices we choose. There might be some overlap among the independent sets generated by different values of $k$, but this simple construction can be used to show $\textup{mis}(B(2d-1,d))\geq (1+o(1))(2d-1)2^{\binom{2d-2}{d}}$.

To get a better lower bound, we use the following improvement. For fixed integers $k\in [2d-1]$ and $0\leq m\leq  d^2$, choose $m$ vertices $v_1,\ldots, v_m\in \mathcal{L}_{d}$ such that $k\notin v_i$ for every $i\in [m]$ and $d(v_i,v_j)\geq 10$ for any distinct $i,j\in [m]$. Fix vertices $v_1',\ldots, v_m'\in \mathcal{L}_{d-1}$ such that $d(v_i,v_i')=3$ and $k\in v_i'$ for each $i\in [m]$. By the definition of $v_i'$, there are two distinct integers $a,b\in v_i$ such that  $v_i'=(v_i\cup\{k\})\setminus\{a,b\}$. Therefore, there are exactly $2d-2$ edges of $M_k$ such that there is an edge from one of their endpoints to at least one of the vertices $v_i,v_i'$. Since $d(v_i,v_j)\geq 10$ for $i\neq j$, it follows that there are exactly $m(2d-2)$ edges adjacent to at least one of the vertices in $\{v_1,\ldots, v_m, v_1',\ldots, v_m'\}$. 

We choose the vertex set $\{v_1,\ldots, v_m, v_1',\ldots, v_m'\}$ as described above and add it to $I$. From each of the remaining $\binom{2d-2}{d-1}-m(2d-2)$ edges in $M_k$ which are not adjacent to a vertex in $\{v_1,\ldots, v_m, v_1',\ldots, v_m'\}$, we choose one of the endpoints and add it to $I$. There is a maximal independent set containing the chosen vertices, which is different with different choices of $\{v_1,\ldots, v_m,$ $ v_1',\ldots, v_m'\}$. There are at least $$\frac{1}{m!}\left(\prod_{i=1}^{m}\left(\binom{2d-2}{d}-(i-1)d^{10}\right)\frac{d(d-1)}{2}\right)$$ choices for $\{v_1,\ldots, v_m, v_1',\ldots, v_m'\}$, and at least $2^{\binom{2d-2}{d-1}-m(2d-2)}$ choices for $I\cap V(M_k)$. Let $\mathcal{I}_{k}^{m}$ be the set of maximal independent sets obtained from $M_k$ following this procedure, then
\[
\begin{split}
|\mathcal{I}_{k}^{m}|&\geq \frac{1}{m!}\left(\left(\binom{2d-2}{d}-md^{10}\right)\frac{d(d-1)}{2}\right)^m 2^{\binom{2d-2}{d-1}-m(2d-2)}\\
&=  \frac{1}{m!}\left(\frac{\binom{2d-2}{d}(d(d-1))} {2^{2d-1} }\right)^m \left(1-\frac{md^{10}}{\binom{2d-2}{d}}\right)^{m}  2^{\binom{2d-2}{d-1}}\\
&\geq  \frac{1}{m!}\left(\frac{\binom{2d-2}{d-1}(d-1)^{2}} {2^{2d-1}} \right)^m  \left(1-\frac{d^{12}}{\binom{2d-2}{d}}\right)^{d^2} 2^{\binom{2d-2}{d-1}}.
\end{split}
\]

\noindent Since $\lim_{x\rightarrow \infty}(1-1/x)^{x}=e^{-1}$, there is a $d_0$ such that for $d\geq d_0$ we have 
\[
\left(1-\frac{d^{12}}{\binom{2d-2}{d}}\right)^{\binom{2d-2}{d}/d^{12}}\geq e^{-2},\]
which implies that for $d\geq d_0$ we have
\[
\left(1-\frac{d^{12}}{\binom{2d-2}{d}}\right)^{d^2}=\left(1-\frac{d^{12}}{\binom{2d-2}{d}}\right)^{\left(\binom{2d-2}{d}/d^{12}\right)\left(d^{14}/\binom{2d-2}{d}\right)}\geq \exp\left(-\frac{2d^{14}}{\binom{2d-2}{d}}\right)=1+o(1).
\]

\noindent Therefore,
\[
\begin{split}
|\mathcal{I}_{k}^{m}|&\geq (1+o(1))\frac{1}{m!}\left(\frac{(d-1)^2}{2^{2d-1}}\binom{2d-2}{d-1}\right)^{m}\cdot 2^{\binom{2d-2}{d-1}}.
\end{split}
\]
Set $y\coloneqq \frac{(d-1)^2}{2^{2d-1}}\binom{2d-2}{d-1}< d^{3/2}$. By Taylor's Theorem, there is a real number $z\in [0,y]$ such that 
\[
\left|e^{y}-\sum_{m=0}^{d^2}\frac{y^m}{m!}\right|\leq \frac{e^{z}}{(d^2+1)!}z^{d^2+1}.
\]
Using that $z< d^{3/2}$ and $(d^2+1)!> \left(d^2/e\right)^{d^2}$, we obtain
\[
\left|e^{y}-\sum_{m=0}^{d^2}\frac{y^m}{m!}\right|\leq \frac{e^{d^{3/2}}}{(d^2+1)!}(d^{3/2})^{d^2+1}\leq \exp\left(2d^{2}+\frac{3}{2}(d^2+1)\ln d-2d^2\ln d\right)=o(1),
\]
implying
\begin{equation}\label{ineq:lowerInclusionFirstTerm}
\sum_{m=0}^{d^{2}} |\mathcal{I}_{k}^{m}|\geq (1+o(1)) \exp\left(\frac{(d-1)^2}{2^{2d-1}}\binom{2d-2}{d-1}\right)\cdot 2^{\binom{2d-2}{d-1}}.
\end{equation}

\noindent By the inclusion-exclusion principle, we obtain
\begin{equation}\label{ineq:lowerBoundInclusion-exclusion}
    \textup{mis}(B(2d-1,d))\geq \left|\bigcup_{k\in [2d-1]}\text{ } \bigcup_{m=0}^{d^2}\mathcal{I}_{k}^{m}\right|\geq \sum_{k\in [2d-1]}\text{ }\sum_{m=0}^{d^2} |\mathcal{I}_{k}^{m}|-\sum_{k_1, k_2\in [2d-1]}\text{ }
    \sum_{m_1,m_2}|\mathcal{I}_{k_1}^{m_1}\cap\mathcal{I}_{k_2}^{m_2}|,
\end{equation}
 where the last summation is over all ordered pairs of non-negative integers $m_1,m_2\leq d^{2}$ such that $(k_1,m_1)\neq (k_2,m_2)$. We proceed to obtain an upper bound for $|\mathcal{I}_{k_1}^{m_1}\cap\mathcal{I}_{k_2}^{m_2}|$.

\begin{claim} \label{claim:Upper bound Intersection Ik}
For every pair of integers $k_1,k_2\in [2d-1]$ and non-negative integers $m_1,m_2\leq d^2$ such that $(k_1,m_1)\neq (k_2,m_2)$ we have 
\[
|\mathcal{I}_{k_1}^{m_1}\cap \mathcal{I}_{k_2}^{m_2}|\leq 3^{\binom{2d-4}{d-2}}\cdot 2^{2\binom{2d-4}{d}+d^5}.
\]
\end{claim}
 \begin{proof}
For a maximal independent set $I$, an edge $e=uv$ is \textit{occupied} if $I\cap \{u,v\}\neq\emptyset$. For each $I\in \mathcal{I}^{m}_{k}$, there are exactly $\binom{2d-2}{d-1}-(2d-2)m$ occupied edges in $M_k$. Therefore, for $m_1\neq m_2$ and $k\in[2d-1]$, we have $\mathcal{I}_{k}^{m_1}\cap \mathcal{I}_{k}^{m_2}=\emptyset$. We proceed to obtain an upper bound for $|\mathcal{I}_{k_1}^{m_1}\cap \mathcal{I}_{k_2}^{m_2}|$ assuming $k_1\neq k_2$.  We also fix an auxiliary integer $k_3\in [2d-1]$ distinct from $k_1$ and $k_2$.

For distinct integers $k_1,k_2,k_3\in [2d-1]$ and $v\in \mathcal{L}_{d-2}^{[2d-1]\setminus\{k_1,k_2,k_3\}}$, let     
\[
U_{v}\coloneqq \{v\cup\{k_1\}, v\cup\{k_2\},v\cup\{k_3\}, v\cup\{k_1,k_2\},v\cup\{k_1,k_3\}, v\cup\{k_2,k_3\}\}.
\]
Each $B(2d-1,d)[U_v]$ is isomorphic to an induced $C_6$ and contains two edges from each of $M_{k_1}$, $M_{k_2}$ and $M_{k_3}$. We say that $U_v$ is \emph{occupied} if the four edges from $M_{k_1}$ and $M_{k_2}$ spanned by $U_v$ are occupied. Notice that if $U_v$ is occupied, then there are $3$ possible choices for $I\cap U_v$, namely

\begin{equation*}
I\cap U_{v}=\begin{cases}\{v\cup\{k_1\}, v\cup\{k_2\}, v\cup\{k_3\}\}, \\\{v\cup\{k_1,k_2\}, v\cup\{k_2,k_3\}, v\cup\{k_1,k_3\}\}, \text{ or}\\\{v\cup\{k_1,k_2\}, v\cup\{k_3\}\}.
\end{cases}
\end{equation*}

The information obtained from the edges $M_{k_2}$ and $M_{k_3}$ in each occupied $U_v$ provides us with fewer possible choices for $I\cap U_v$ than if we only choose the intersection of $I$ with the occupied edges from $M_{k_1}$ induced by $U_v$. That is, ignoring the edges from $M_{k_2}$ and $M_{k_3}$, there are $4$ choices instead of $3$ for $I\cap U_v \cap V(M_{k_1})$. If $U_v$ is not occupied, then there are $9$ possible choices for $I\cap U_v\cap V(M_{k_1})$. The number of unoccupied $U_v$'s will be negligible.

We describe a procedure that generates every maximal independent set $I\in \mathcal{I}_{k_1}^{m_1}\cap \mathcal{I}_{k_2}^{m_2}$ and we use it to provide an upper bound for $|\mathcal{I}_{k_1}^{m_1}\cap \mathcal{I}_{k_2}^{m_2}|$. 

First we choose the set of vertices $\{v_1,\ldots, v_{m_1}, v_1',\ldots, v_{m_1}'\}$. This choice determines which edges from $M_{k_1}$ are occupied. Then, we choose the set of vertices $\{u_1,\ldots, u_{m_2}, u_{1}',\ldots,u_{m_2}'\}$. This second choice determines which edges from $M_{k_2}$ are occupied. The total number of choices up to now is at most $\left(\binom{2d-2}{d}d^2\right)^{m_1+m_2}$. 

Since $\binom{2d-2}{d-1}-2dm_1$ edges from $M_{k_1}$ are occupied and $\binom{2d-2}{d-1}-2dm_2$ edges from $M_{k_2}$ are occupied, it follows that at least $\binom{2d-4}{d-2}-2(m_1+m_2)d$ copies of $C_6$ of the form $B(2d-1,d)[U_{v}]$ are occupied. Therefore, there are at most  $3^{\binom{2d-4}{d-2}}$ choices for the occupied $U_v$'s and $9^{2(m_1+m_2)d}$ choices for the unoccupied $U_{v}$'s. The number of edges from $M_{k_1}$ that are occupied and not specified yet is at most $2\binom{2d-4}{d}$. For each of these edges we choose one of the endpoints to be in $I$. This procedure generates all the possible independent sets $I\in \mathcal{I}_{k_1}^{m_1}\cap \mathcal{I}_{k_2}^{m_2}$. It follows
\begin{equation*}
|\mathcal{I}_{k_1}^{m_1}\cap \mathcal{I}_{k_2}^{m_2}|\leq \left(\binom{2d-2}{d}d^2\right)^{m_1+m_2} 2^{2 \binom{2d-4}{d}}\cdot 3^{\binom{2d-4}{d-2}}\cdot 9^{2(m_1+m_2)d} < 3^{\binom{2d-4}{d-2}}\cdot 2^{2\binom{2d-4}{d}+d^5}.
\end{equation*}
\end{proof}
Claim~\ref{claim:Upper bound Intersection Ik} implies 
\begin{equation}\label{ineq:upperSecondInclusionExclusionTerms}
\sum_{k_1, k_2\in [2d-1]}\text{ }\sum_{m_1,m_2}|\mathcal{I}_{k_1}^{m_1} \cap \mathcal{I}_{k_2}^{m_2}|=o\left(2^{\binom{2d-2}{d-1}}\right),\end{equation} 
where the last sum is over all ordered pairs of non-negative integers $m_1,m_2\leq d^{2}$ such that $(k_1,m_1)\neq (k_2,m_2)$. From~\eqref{ineq:lowerInclusionFirstTerm},~\eqref{ineq:lowerBoundInclusion-exclusion} and~\eqref{ineq:upperSecondInclusionExclusionTerms} we conclude
\[
\textup{mis}(B(2d-1,d))\geq (1+o(1))(2d-1) \exp\left(\frac{(d-1)^2}{2^{2d-1}}\binom{2d-2}{d-1}\right) \cdot 2^{\binom{2d-2}{d-1}}. \qedhere
\]
\end{proof}


\section{Most MIS have a large induced matching} \label{sec:containers}

In this section we prove Theorems~\ref{thm:FirstPhase} and~\ref{thm:SecondPhase}. The proof of Theorem~\ref{thm:FirstPhase} is similar to the proof of Lemma 1.2 in~\cite{Kahn2022}.

\subsection{First stage containers}\label{Sec:FirstStage}
Fix an arbitrary linear ordering $\prec_1$ on $V(B(2d-1,d))$ and a maximal independent set $I$. This linear ordering is for tiebreaking purposes only. Set $Z_{0}\coloneqq \L_{d-1}\cup \L_{d}$. For every $i\geq 0$, we repeat the following procedure: Pick a vertex $v\in Z_i$ with maximum degree in $B(2d-1,d)[Z_i]$. If there is more than one such vertex, then we choose the first under $\prec_1$.
\begin{itemize}
\itemsep-0.3em
    \item[(i)]If $v\in I$, then set $Z_{i+1}\coloneqq Z_i\setminus (N(v)\cup \{v\})$ and $C_{i+1}\coloneqq C_{i}\cup\{v\}$.
    \item[(ii)] If $v\notin I$, then set $Z_{i+1}\coloneqq Z_i\setminus \{v\}$ and $C_{i+1}\coloneqq C_{i}$. 
\end{itemize}
We stop the algorithm when at least one of the two conditions hold: 
\begin{itemize}
\itemsep-0.3em
\item[(a)] $|Z_i|\leq 2\left(1-\frac{\log^3 d}{d}\right)\binom{2d-2}{d-1}$;
\item[(b)]$|C_{i}|\geq \frac{\log d}{d}\binom{2d-1}{d}$.
\end{itemize}
We create a vector $\xi\in \{0,1\}^{|Z_{0}|}$ such that $\xi_{i}=1$ if case (i) holds at step $i$ and $\xi_{i}=0$ otherwise. Once the algorithm stops, we fill the remaining entries of $\xi$ with $0$'s. We will refer to the above process as \emph{running the [Algorithm] for $Z_0$ with stopping conditions (a) and (b).}

Notice that at each iteration, regardless of whether case (i) or (ii) holds, we remove at least one element from $Z_i$. Thus,  condition (a) guarantees that the algorithm stops no matter which maximal independent set $I$ we provide as input.

\begin{proof}[Proof of Theorem~\ref{thm:FirstPhase}]
For each $I\in \mathcal{I}\setminus \mathcal{J}_1$ as input, we run the algorithm. Given $\xi\in \{0,1\}^{|Z_0|}$, our goal is to obtain an upper bound on the number of ways to recover $I\in \mathcal{I}\setminus \mathcal{J}_1$.

If we are provided $\xi$, then we know $j$, the last index for which we added a vertex to the certificate $C_{j}$.

\vspace{0.1in}
\noindent\textbf{Case: }The algorithm stops by condition (a).

\noindent Then $|C_j|=|\text{supp} (\xi)|\leq \frac{\log d}{d}\binom{2d-1}{d}$ and the number of ways to perform this algorithm is at most the number of ways to choose the $1$-entries of $\xi$, the cost of which is  
\[
\log \binom{2\binom{2d-1}{d}}{\leq \frac{\log d}{d}\binom{2d-1}{d}}=O\left(\frac{\log^2 d}{d}\binom{2d-1}{d}\right)=O\left(\frac{\log^2 d}{d}\binom{2d-2}{d-1}\right),
\]
where the first equality follows from Proposition~\ref{prop:Entropy}.

Notice that $I\cap Z_j$ is a maximal independent set in $B(2d-1,d)[Z_j]$. By Theorem~\ref{thm:HujterTuza}, we have that the cost for choosing $I\cap Z_j$ is at most $|Z_j|/2\leq \left(1-\frac{\log^3d}{d}\right)\binom{2d-2}{d-1}$. The total cost is the cost of running the algorithm plus the cost of obtaining $I\cap Z_j$, which is at most
\[
\left(1-\frac{\log^3 d}{d}+O
\left(\frac{\log^2 d}{d}\right)\right)\binom{2d-2}{d-1}
=\left(1-\Omega\left(\frac{\log^3 d}{d}\right)\right)\binom{2d-2}{d-1}.
\]

We conclude that the number of independent sets $I\in \mathcal{I}\setminus \mathcal{J}_1$ for which the algorithm stops because of condition (a) is $o(\mathcal{|I|})$.

\vspace{0.1in}
\noindent \textbf{Case: } The algorithm stops with $|C_{j}|\geq \frac{\log d}{d}\binom{2d-1}{d}$ and 
\begin{equation}\label{ineq:lowerZj}
    |Z_j| > 2\left(1-\frac{\log^3 d}{d}\right)\binom{2d-2}{d-1}.
\end{equation}

\begin{claim}\label{claim:smallZi} For each $0\leq i\leq j$ let $d_i\coloneqq\max_{v\in Z_i} d_{Z_i}(v)$, then we have
\[
|Z_i|\leq \left(1+\frac{d_i}{d}\right)\binom{2d-1}{d}.
\]
\end{claim}

\begin{proof}
By Proposition~\ref{prop:UpperBoundZ} with $Z=Z_i, W=\mathcal{L}_d\cup \mathcal{L}_{d-1}$, $d'=d_i$ and $L=0$, we have

\[
|Z_i|\leq \frac{2d\binom{2d-1}{d}}{2d-d_i}=\left(1+\frac{d_i}{2d-d_i}\right)\binom{2d-1}{d}.
\]
The claim  follows from the fact that $d_i\leq d$ implies $d_i/(2d-d_i)\leq d_i/d$.
\end{proof}

For each $0\leq i\leq j$ let $\alpha_i$ be a real number such that $|Z_i|=(1+\alpha_i)\binom{2d-1}{d}$. In particular, $\alpha_0=1$.  By Claim~\ref{claim:smallZi} we have $\alpha_i \leq d_i/d$.

\begin{claim}\label{claim:UpperAlpha}
    If $\xi_i=1$, i.e.,~we add a vertex to $C_i$, then $\alpha_i\leq \left(1-\frac{d}{\binom{2d-1}{d}}\right)\alpha_{i-1}$.
\end{claim}
\begin{proof}
Using that $\xi_i=1$ we obtain that $|Z_i|=|Z_{i-1}|-d_{i-1}-1$. It follows
    \[
    (1+\alpha_i)\binom{2d-1}{d}=|Z_i|=|Z_{i-1}|-d_{i-1}-1<(1+\alpha_{i-1})\binom{2d-1}{d}-\alpha_{i-1}d,
    \]
    which implies
    \[
    \alpha_i\binom{2d-1}{d}\leq \alpha_{i-1}\binom{2d-1}{d}-\alpha_{i-1}d.
    \]
\end{proof}

\noindent Using the fact that $|\text{supp}( \xi)|\geq \frac{\log d}{d}\binom{2d-1}{d}$ and Claim~\ref{claim:UpperAlpha} we obtain an upper bound on $|Z_j|$.

\begin{align*}
|Z_j|=(1+\alpha_j)\binom{2d-1}{d} &\leq \left(1+\left(1-\frac{d}{\binom{2d-1}{d}}\right)^{\frac{\log d}{d}\binom{2d-1}{d}}\right)\binom{2d-1}{d}.\\
\end{align*}

\noindent From the fact that $(1-x)\leq e^{-x}$ for every real number $x$, it follows
\begin{equation}\label{ineq:upperZj}
|Z_j|\leq \left(1+e^{\log (1/d)}\right)\binom{2d-1}{d}\leq\left(1+\frac{1}{d}\right)\binom{2d-1}{d}<2\left(1+\frac{1}{d}\right)\binom{2d-2}{d-1}.
\end{equation}

We construct an auxiliary graph, for which we apply Theorem~\ref{thm:Hujter-TuzaStab}. Let $Z^{*}\coloneqq Z_{j}\cup C_{j}$. Given a maximal independent set $I$ in $B(2d-1,d)$, we have that $I\subseteq Z^*$ is also a maximal independent set in $G'\coloneqq B(2d-1,d)[Z^{*}]$. 

Notice that $\prec$, the linear ordering for the induced matchings in $B(2d-1,d)$, induces a linear ordering for the induced matchings in $G'$. Therefore, we have $|M_{G'}(I)|\leq |M(I)|$.

From~\eqref{ineq:upperZj} and $|C_j|\leq \frac{2\log d}{d}\binom{2d-1}{d}\leq \frac{4\log d}{d}\binom{2d-2}{d-1}$, we obtain
\begin{equation}\label{ineq:UpperZstar}
    |Z^{*}|\leq 2\left(1+\frac{3\log d}{d}\right)\binom{2d-2}{d-1}.
\end{equation}

\noindent Since $|Z^*|\geq |Z_j|$, by~\eqref{ineq:lowerZj} we obtain
\begin{equation}\label{ineq:lowerXY}
|Z^{*}|\left(1-\frac{\log^3 d}{d}\right)\geq 2\left(1-\frac{\log^3 d}{d}\right)^2 \binom{2d-2}{d-1}\geq 2\left(1-\frac{2\log^3 d}{d}\right)\binom{2d-2}{d-1}.
\end{equation}

\noindent If $I\in \mathcal{I}\setminus \mathcal{J}_1$, then by~\eqref{ineq:lowerXY} and the definition of $\mathcal{J}_1$, 
\begin{equation}\label{ineq:UpperMG'}
|M_{G'}(I)|\leq |M(I)|\leq \left(1-\frac{2\log^3 d}{d}\right)\binom{2d-2}{d-1}\leq \left(1-\frac{\log^3 d}{d}\right)\frac{|Z^{*}|}{2}.
\end{equation}
Applying Theorem~\ref{thm:Hujter-TuzaStab} on $G'$ with $\varepsilon=\log^3 d/d$, we obtain
\[
\log \left|\left\{I\in \I(G'): |M_{G'}(I)|\leq \left(1-\frac{\log^3 d}{d}\right)\frac{|Z^{*}|}{2}\right\}\right|\leq \left(1-\frac{c\log^3 d}{d}\right)\frac{|Z^{*}|}{2},
\]
where $c$ is the constant in Theorem \ref{thm:Hujter-TuzaStab}. By~\eqref{ineq:UpperZstar},
\begin{equation*}
\left(1-\frac{c\log^3 d}{d} \right)\frac{|Z^{*}|}{2}\leq \left(1-\frac{c\log^3 d}{d}\right)\left(1+\frac{3\log d}{d}\right)\binom{2d-2}{d-1}=\left(1-\Omega\left(\frac{\log^3 d}{d}\right)\right)\binom{2d-2}{d-1}.
\end{equation*}
Inequality~\eqref{ineq:UpperMG'} shows that $I\in \mathcal{I}\setminus \J_{1}$ implies $|M_{G'}(I)|\leq (1-\log^3 d/d)|Z^{*}|/2$. Therefore, we can specify $I\in \mathcal{I}\setminus \J_{1}$ by choosing a vector $\xi$ with $|\text{supp}(\xi)|\leq \frac{2\log d}{d}\binom{2d-1}{d}$; with the algorithm we recover $Z_j$ and $C_j$, thus we recover $Z^{*}$; finally, we choose $I\subseteq V(G')$ from $\{I\in \I({G'}): M_{G'}(I)\leq (1-\log^3 d/d)|Z^{*}|/2\}$. It follows that the cost to specify $I\in \mathcal{I}\setminus \mathcal{J}_1$ for which the algorithm stops by condition (b) is
\[
\left(O\left(\frac{\log^2 d}{d}\right)+1-\Omega\left(\frac{\log^3 d}{d}\right)\right)\binom{2d-2}{d-1}=\left(1-\Omega\left(\frac{\log^3 d}{d}\right)\right)\binom{2d-2}{d-1}.
\]
Since $2^{\left(1-\Omega\left(\frac{\log^3 d}{d}\right)\right)\binom{2d-2}{d-1}}=o(|\mathcal{I}|)$, we conclude
\[
|\mathcal{I}\setminus\mathcal{J}_1|=o(|\mathcal{I}|).\qedhere
\]
\end{proof}


\subsection{Second stage containers}\label{Subsec:SecondStage}
In this subsection we prove Theorem~\ref{thm:SecondPhase}. Through this subsection we will assume $I\in \mathcal{J}_1\setminus\mathcal{J}_2 $, which means we can use that $M(I)$ is considerably large. We will use a modified version of Sapozhenko's graph container lemma to exploit the fact that most edges in $B(2d-1,d)$ have at least one endpoint in $V(M(I))$.

Fix a set $I\in \mathcal{J}_1\setminus \mathcal{J}_2 $ and choose $\beta$ so that $|M(I)|=(1-\beta)\binom{2d-2}{d-1}$. By the definitions of $\mathcal{J}_1$ and $\mathcal{J}_2$,
\begin{equation}\label{eq:betabounds}
\frac{2\log^5 d}{d^{3/2}}\leq \beta< \frac{2\log^3 d}{d}.
\end{equation}
Let $E_0$ be the set of edges with at least one endpoint in $V(M(I))$. Each edge $xy\in M(I)$ contributes exactly $2d-1$ edges to $E_0$ with at least one endpoint in $\{x,y\}$, hence  
\begin{equation}\label{eq:E0bounds}
|E_0|=(2d-1)(1-\beta)\binom{2d-2}{d-1}=d\binom{2d-1}{d}-(2d-1)\beta\binom{2d-2}{d-1}.
\end{equation}
No edge with one endpoint in $I\setminus V(M(I))$ belongs to $E_0$, so
\begin{equation*}
\begin{split}
|I\setminus V(M(I))|&\leq \frac{1}{d}|E(B(2d-1,d))\setminus E_0|=\frac{2d-1}{d}\beta \binom{2d-2}{d-1}=2\beta \binom{2d-2}{d-1}-\frac{\beta}{d}\binom{2d-2}{d-1}.
\end{split}
\end{equation*}

Set $A_{1}\coloneqq (V(M(I))\cup I)\cap \mathcal{L}_{d}$ and $A_{2}\coloneqq (V(M(I))\cup I)\cap \mathcal{L}_{d-1}$. Without loss of generality we may assume that $|A_1|\leq |A_2|$. Thus, $|A_1|=|M(I)|+|(I\setminus V(M(I)))\cap \mathcal{L}_d|$ and 

\begin{equation}\label{ineq:sizeA}
|A_1|\leq \left(1-\beta+\beta-\frac{\beta}{2d}\right)\binom{2d-2}{d-1}\leq  \left(1-\frac{\log^5 d}{d^{5/2}}\right)\binom{2d-2}{d-1}.
\end{equation}

\noindent Let $H_1\coloneqq N(A_1)\setminus A_2$ and $H_2\coloneqq N(A_2)\setminus A_1$. Finally let
\begin{equation}
r_1\coloneqq e(H_1,\mathcal{L}_{d}\setminus A_1). \end{equation}

\noindent From~\eqref{eq:E0bounds} and~\eqref{eq:betabounds} we get
\begin{equation}\label{ineq:upperBoundr}
r_1=e(H_1,\mathcal{L}_d\setminus A_1)\leq d\binom{2d-1}{d}-|E_0|= (2d-1)\beta \binom{2d-2}{d-1}\leq 2\log^3 d \binom{2d-1}{d}.
\end{equation}

Let $\{A_1^{1},A_2^{1},\ldots\}$ be the collection of $2$-linked components of $A_1=(V(M(I))\cup I)\cap \mathcal{L}_{d}$. For each $i$, set $H_{i}^{1}\coloneqq N(A_{i}^{1})\setminus A_{2}$ and $r_{i}^{1}\coloneqq e(H_{i}^{1},\mathcal{L}_{d}\setminus A_{i}^{1})$. By the definition of $2$-linked component, it follows that the sets $N(A_i^{1})$ are mutually disjoint, thus the sets $H_{i}^{1}$ are mutually disjoint as well. We also have $E(H_{i}^{1},A_1\setminus A_i^1)=\emptyset$ for every $i$, and $\sum_{i}r_{i}^{1}=r_{1}$.

Let $L\subseteq \mathbb{N}$ be the set of integers $i$ for which $|A_{i}^1|\leq d^6$. By Lemma~\ref{lem:isoperimetry}, for every $i\in L$ we have $|N(A_i^{1})|\geq d|A_i^{1}|/9$, which implies $|H_{i}^1|\geq (d/9)|A_{i}^1|-|A_i^1|$. It follows that
\begin{equation}\label{ineq:edgeGi}
\begin{aligned}
    r_{i}^{1}&=e(H_{i}^{1},\mathcal{L}_{d}\setminus A_i^1)=d|H_{i}^{1}|-e(H_{i}^{1}, A_i^1)\geq d|H_{i}^{1}|-d|A_i^{1}|\geq \left(\frac{d^2}{9}-2d\right)|A_i^1|\geq \frac{d^2}{10}|A_i^1|
\end{aligned}
\end{equation}
for sufficiently large $d$. Summing up~\eqref{ineq:edgeGi} for all $i\in L$ and using~\eqref{ineq:upperBoundr}, we get
\begin{equation}\label{ineq:upperSmallAi}
|\cup_{i\in L}A_i^1|=\sum_{i\in L}|A_i^1|\leq \frac{20\log^3 d}{d^2}\binom{2d-1}{d}.
\end{equation}
Since $|\cup_{i\in L}A_i^{1}|$ is relatively small, we will specify it by simply choosing it from $\mathcal{L}_{d}$.

Our goal is to use a variation of Sapozhenko's graph container lemma for the $2$-linked components $A_{i}^{1}\subseteq A_1\setminus \cup_{i\in L} A_{i}^{1}$. We will obtain sets $F_i, S_i$ such that $F_i\subseteq H_i^1$ and $A_i^1\subseteq S_i$ for $i\notin L$, which will be good approximations for $H_{i}^{1}$ and $A_i^1$, respectively.

Our final goal for this subsection is to set $F_1\coloneqq \cup_{i\notin L} F_i$ and $S_1\coloneqq \cup_{i\notin L}S_i$, and reproduce the method at the end of Section~\ref{Sec:FirstStage}, using the graph induced by $S_1\cup (\mathcal{L}_{d-1}\setminus F_1)\cup(\cup_{i\in L}A_{i}^{1})$ this time.

Fix a set $A\in \{A_{1}^1,A_{2}^{1},\ldots\}$ such that $|A|>d^6$ and let $H\coloneqq N(A)\setminus A_2$. Since $A\subseteq A_1$, by~\eqref{ineq:sizeA} we have $|A|\leq (1-\log^5 d/d^{5/2})\binom{2d-2}{d-1}$. Thus, from Lemma~\ref{lem:edgeIso} we obtain
\begin{equation}\label{ineq:lowert}
r\coloneqq e(H,\mathcal{L}_d\setminus A)\geq |A|\left(1-\frac{|A|}{\binom{2d-2}{d-1}}\right)\geq |A|\left(1-\left(1-\frac{\log^5 d}{d^{5/2}}\right)\right)=\frac{|A|\log^ 5 d}{d^{5/2}}.
\end{equation}
Since we also have $h\coloneqq |H|\leq |N(A)|\leq d|A|$, from~\eqref{ineq:lowert} we obtain
\begin{equation}\label{ineq:Uppergt}
    \frac{h}{r}\leq \frac{d^{7/2}}{\log^5 d}<d^4.
\end{equation}

\noindent Now we start building our containers for the pair $(A,H)$.

\begin{defi}\label{defi:G(A,H)}
For arbitrary positive integers $a$, $h$ and $r$, 
let $\mathcal{G}_1(a,h,r)\subseteq 2^{\mathcal{L}_d}\times 2^{\mathcal{L}_{d-1}}$ denote the collection of ordered pairs of sets $(A,H)$ such that:
\begin{itemize}
\itemsep-0.3em
\item[i)]  $A\subseteq \mathcal{L}_{d}$, $A$ is $2$-linked and $|A|=a$.
\item[ii)] $H\subseteq N(A)$, $|H|=h$ and if $u\in \mathcal{L}_{d}$ has $N(u)\subseteq H$ then $u\in A$. 
\item[iii)] $e(H,\mathcal{L}_{d}\setminus A)=r$.
\item[iv)] Set $A'\coloneqq N(A)\setminus H$. There is an induced matching $M\subseteq E(A, A')$ in $B(2d-1,d)[A,N(A)]$, such that $M$ covers $A'$ and $N(V(M))\cap (A\setminus V(M))=\emptyset$. 
\end{itemize}
\end{defi}
\begin{rem}
    For every set $A$, there might be more than one set $H$ such that $(A,H)\in \mathcal{G}_1(a,h,r)$. Our goal is not to bound $|\mathcal{G}_1(a,h,r)|$ directly and we do not have to address this potential problem.

In contrast with other applications of Sapozhenko's container lemma, in Definition~\ref{defi:G(A,H)}, $A$ is not closed but contains all vertices $u$ with $N(u)\subseteq H$. 
\end{rem}

The first step is obtaining a $\varphi$-approximation. Informally, that means we first create a set $F'\subseteq H$ covering $A$ such that $F'$ contains all vertices in $H$ with high degree towards $A$.

\begin{defi}
    For every $A\subseteq \mathcal{L}_d, H\subseteq N(A)$ and $1\leq \varphi\leq d$, let 
\[
H^{\varphi}\coloneqq\{u\in H: d_{A}(u)\geq \varphi\}.
\]
\end{defi}
In our applications, we will set $\varphi=d/2$. We state our results in general form.
\begin{defi}
    For arbitrary positive integers $h, a, r$ and an arbitrary pair $(A,H)\in \mathcal{G}_1(a,h,r)$, a $\varphi$-approximation of $(A,H)$ is a set $F'$ such that $H^{\varphi}\subseteq F' \subseteq H$ and $A\subseteq N(F')$.
\end{defi}

\begin{lemma}\label{lem:phiaprox}
    For arbitrary positive integers $h, a\geq d^6, r\geq h\log^5 d/d^{7/2}$ and $100\leq \varphi\leq d$, there is a family $\mathcal{V}=\mathcal{V}(a,h,r,\varphi)\subseteq 2^{\mathcal{L}_{d-1}}$ such that every $(A,H)\in \mathcal{G}_1(a,h,r)$ has a $\varphi$-approximation $F'\in \mathcal{V}$ and 

    \[
    \log|\mathcal{V}|= \begin{cases}
    O\left(\frac{h\log^2 d}{\varphi d}\right)+O\left(\frac{r\log^2 d}{d \varphi }\right)&\text{ if } r\leq \frac{h(d-\varphi)}{\varphi },\\
    O\left(\frac{r\log^2 d}{d(d-\varphi)}\right)+O\left(\frac{r\log^2 d}{d\varphi} \right)&\text{ if } r> \frac{h(d-\varphi)}{\varphi }.
    \end{cases}
    \]
\end{lemma}
\begin{proof}
    Let $p=40\frac{\log d}{\varphi d}$. Let $X$ be a random subset of $H$, where each vertex in $H$ belongs to $X$ independently with probability $p$. Let $\Omega_{X}\coloneqq E(X, N(X)\setminus A)$. Notice that $\mathbb{E}(|X|)=ph$ and $\mathbb{E}(|\Omega_X|)=pr$. Set \[
    T_{X}\coloneqq\{v\in \mathcal{L}_{d-1}: |N_{A}(v)\cap N_{A}(X)|\geq 2\},
    \]
    then we have $T_{X}\subseteq N(A)$. Let $M\subseteq E(A, N(A)\setminus H)$ be the induced matching covering $N(A)\setminus H$ guaranteed by the definition of $\mathcal{G}_1(a,h,r)$. Then, for each $v\in V(M)\cap\mathcal{L}_{d-1}$, we have $|N_A(v)|=1$, which implies that $T_{X}\cap V(M)=\emptyset$. We conclude that $T_{X}\subseteq N(A)\setminus V(M)=H$.

    Since $|N_A(y)|\geq \varphi$ for every $y\in H^{\varphi}$, by Lemma~\ref{lem:isoperimetry} we have  $|N(N_{A}(y))|\geq \varphi d/9$, which implies 
    \[
    |N(N_{A}(y))\cap H|\geq \varphi d/9-|N_A(y)|\geq \varphi d/9-d\geq \varphi d/10.
    \]
    Furthermore, if $N_{A}(y)\cap N_A(X)=\{w\}$ for $y\in H^{\varphi}$, then $|N_A(y)\setminus\{w\}|\geq \varphi-1$ and $N_A(y)\setminus\{w\}$ has no neighbors in $X$. Similarly, we obtain \[
    |N(N_{A}(y)\setminus\{w\})\cap H|\geq (\varphi-1)d/9-|N_{A}(y)\setminus\{w\}|\geq \varphi d/10.\]
    It follows that for every $y\in H^{\varphi}$, we have
    \begin{align*}
    \mathbb{P}(y\notin T_X)&\leq \mathbb{P}(|N_{A}(y)\cap N_A(X)|=0)+\mathbb{P}(|N_{A}(y)\cap N_A(X)|=1)\\
    &\leq (1-p)^{\varphi d/10}+d(1-p)^{\varphi d/10}\leq  (d+1)\exp (-40\log d /10)\leq 2d^{-3},
    \end{align*}
    which implies
    \[
    \mathbb{E}(|H^{\varphi}\setminus T_{X}|)\leq  \frac{2h}{d^{3}}.
    \]

    \noindent By Markov's inequality there are vertex sets $T_0\coloneqq X\subseteq H$, $T_1\coloneqq T_X\subseteq H$ and $T_2\coloneqq  H^{\varphi}\setminus (T_0\cup T_1)$ satisfying
    \begin{equation}\label{uppert0}
    |T_0|\leq \frac{200h\log d}{\varphi d},
    \end{equation}
    \begin{equation}\label{upperOmegat0'}
    e(T_0, N(T_0)\setminus A)\leq \frac{200r\log d}{\varphi d},
    \end{equation} and
    \begin{equation}\label{uppert0'}
    |T_2|\leq |H^{\varphi}\setminus T_1|\leq \frac{10h}{d^{3}}.
    \end{equation}

Let $T^*\coloneqq T_0\cup T_1\cup T_2\supseteq H^{\varphi}$. Up to now, if we are provided $T_0$, $\Omega_{T_0}$ and $T_2$, then we can recover $T_1$ and build $T^*$ such that $H^{\varphi}\subseteq T^*\subseteq H$. This is almost our $\varphi$-approximation, the set $F'$. We are missing the property $A\subseteq N(F')$.

Let $T_3\subseteq H\setminus T^*$ be a minimal cover of $A\setminus N(T^*)$ and set 
\[
F'\coloneqq T^*\cup T_{3}=T_0\cup T_1\cup T_2\cup T_3,
\] 
which is our $\varphi$-approximation. Notice that $T_{i}\subseteq H$ for every $i$, thus $F'\subseteq H$.

Now we proceed to get an upper bound on the number of possible choices of $F'$.
Let 
\[
    T\coloneqq T_0\cup T_2\cup T_3\subseteq N(A).
\]
\noindent We have $d(u,T)\leq 3$ for every $u\in A$ by the definitions of $T_3$ and $T_1$, and $d(u,A)=1$ for every $u\in T$. Using Lemma~\ref{lem:linked}, it follows that $T$ is $8$-linked.

For every $T_0$ and $\Omega_{T_0}=E(T_0,N(T_0)\setminus A)$, we can determine $N_{A}(T_0)$ and $T_1$. Hence, if we are provided with $T$ and $\Omega_{T_0}$, then we can determine $F'$. 

To get an upper bound on $|T_3|$, notice that $T_3\subseteq H\setminus T^* \subseteq H\setminus H^{\varphi}$. We have 
\[
|H\setminus T^*|\leq \frac{e(H\setminus T^*,\mathcal{L}_{d}\setminus A )}{d-\varphi}\leq \frac{e(H,\mathcal{L}_{d}\setminus A)}{d-\varphi}=\frac{r}{d-\varphi}.
\] 
We also have $d_{H\setminus T^*}(u)\geq d-1$ for every $u\in A\setminus N(T^*)$ and trivially $d_{A\setminus N(T^*)}(u)\leq d$ for every $u\in H\setminus T^*$. Applying Theorem~\ref{thm:cover} with $P=A\setminus N(T^*)$, $Q=H\setminus T^*$, $a=d-1$ and $b=d$, we get 
    \begin{equation}\label{uppert1}
    |T_3|\leq  \frac{|H\setminus T^*|}{d-1}(1+\ln d) \leq  \frac{2r\log d}{d(d-\varphi)}.
    \end{equation}

\noindent From~\eqref{uppert0},~\eqref{uppert0'} and~\eqref{uppert1}, it follows that
\begin{equation}\label{upperT}
|T|\leq |T_0|+|T_2|+|T_3|\leq \frac{200 h\log d}{\varphi d}+\frac{10h}{d^{3}}+\frac{2r\log d}{d(d-\varphi)}.
\end{equation}

\noindent With a fixed $T_0$, since $\Omega_{T_0}\subseteq E(T_0,\mathcal{L}_d)$, from~\eqref{uppert0} and~\eqref{upperOmegat0'} it follows that the number of choices of $\Omega_{T_0}$ is at most 
$$\binom{d|T_0|}{\leq \frac{200r\log d}{\varphi d}}\leq\binom{\frac{200h\log d}{\varphi}}{\leq \frac{200r\log d}{\varphi d}}= 2^{O\left(\frac{r\log^2 d}{\varphi d}\right)},$$ 
where we used~\eqref{ineq:Uppergt} and Proposition~\ref{prop:Entropy} to obtain the last equality.

To bound the cost of specifying $T$, we split the proof into two cases.
    
\noindent $\bullet$ If $r\leq \frac{h(d-\varphi)}{\varphi }$: From~\eqref{upperT} we get $|T|\leq\frac{202h\log d}{\varphi d}$. Recall that $T$ is $8$-linked. By Corollary~\ref{cor:klinkednumber}, the number of choices of $T$ is $\binom{2d-1}{d}\cdot 2^{O\left(\frac{h\log^2 d}{\varphi d}\right)}=2^{O\left(\frac{h\log^2 d}{\varphi d}\right)}$. Since $T=T_0\cup T_2\cup T_3$, we have at most $3^{|T|}< 2^{\frac{406h\log d}{\varphi d}}$ choices of $T_0,T_2$ and $T_3$. Given $T_0$, $T_2$, $T_3$ and $\Omega_{T_0}$, we can reconstruct $F'$. Thus, the cost of specifying $F'$ is
$$O\left(\frac{h\log^2 d}{\varphi d}\right)+O\left(\frac{r\log^2 d}{d \varphi }\right).$$

\noindent$\bullet$ If $r>\frac{h(d-\varphi)}{\varphi }$: From~\eqref{upperT} we get $|T|\leq \frac{202r\log d}{d(d-\varphi)}$. Recall that $T$ is $8$-linked. By Corollary~\ref{cor:klinkednumber}, the number of choices of $T$ is $\binom{2d-1}{d}\cdot 2^{O\left(\frac{r\log^2 d}{d(d-\varphi)}\right)}=2^{O\left(\frac{r\log^2 d}{d(d-\varphi)}\right)}$. Since $T=T_0\cup T_2\cup T_3$, we have at most $3^{|T|}<2^{\frac{406r\log d}{d(d-\varphi)}}$ choices of $T_0,T_2$ and $T_3$.  As in the previous case, the cost of choosing $\Omega_{T_0}$ is $O\left(\frac{r\log^2 d}{\varphi d}\right)$. Thus, in this case, the cost of specifying $F'$ is
\[O\left(\frac{r\log^2 d}{d(d-\varphi)}\right)+O\left(\frac{r\log^2 d}{d\varphi} \right).\qedhere\]

\end{proof}
\noindent Now that we determined $F'$, we start the second phase of the container algorithm. 

\begin{defi}\label{defi:psiAprox} For every $1\leq\psi\leq d$, a $\psi$-approximation for $(A,H)\in \mathcal{G}_{1}(a,h,r)$ is a pair $(S,F)$ such that $F\subseteq H$, $A\subseteq S$ and:
\begin{itemize}
\itemsep-0.3em
    \item [i)] $d_{F}(u)\geq d-\psi$ for every $u\in S$.
    \item [ii)] $d_{\mathcal{L}_d\setminus S}(v)\geq d-\psi$ for every $v\in \mathcal{L}_{d-1}\setminus F$.
\end{itemize}
\end{defi}

\noindent We will set $\psi= \sqrt{d} \log d$.
\begin{lemma}\label{lem:psiapprox}
    For arbitrary integers $h,a\geq d^6, r\geq h\log^5 d/d^{7/2}, 100\leq \varphi\leq d$ and $2\leq\psi\leq d$, there is a family $\mathcal{W}=\mathcal{W}(a,h,r,\varphi,\psi)$ with
\[
    \log|\mathcal{W}|= O\left(\frac{r\log d}{\psi (d-\varphi)}\right)+O\left(\frac{r\log d}{\psi(d-\psi)}\right)
\]
    satisfying that if $F'\in \mathcal{V}(a,h,r,\varphi)$ is a $\varphi$-approximation of $(A,H)\in \mathcal{G}_1(a,h,r)$, then there is $(S,F)\in\mathcal{W}$,  a $\psi$-approximation of $(A,H)$, such that $F'\subseteq F\subseteq H$ and $A\subseteq S$.
\end{lemma}

\begin{proof}

We split the proof into two steps. In the first step we make sure that condition (i) in Definition~\ref{defi:psiAprox} is satisfied, while in the second step we make sure that condition (ii) is satisfied.

\textbf{Step 1:} Fix an arbitrary linear ordering $\prec_2$ on $\mathcal{L}_{d}\cup \mathcal{L}_{d-1}$. Following the ordering $\prec_2$, for each $u\in A$, if $d_{H\setminus F'}(u)\geq \psi-1$, then update $F'$ to $F'\cup N(u)$. We stop when there is no $u\in A$ with $d_{H\setminus F'}(u)\geq \psi-1$, i.e., for every $u\in A$ we have $d_{H\setminus F'}(u)<\psi-1$. In the end, set $F^{*}\coloneqq F'$.

Notice that each $v\in H\setminus F'\subseteq H\setminus H^{\varphi}$ satisfies $d_{\mathcal{L}_{d}\setminus A}(v)> d-\varphi$, thus 
\[
(d-\varphi)|H\setminus F'|\leq e(H\setminus F', \mathcal{L}_{d}\setminus A)\leq r,\]
which implies $|H\setminus F'|\leq r/(d-\varphi)$. In each step where we add vertices to $F'$, we remove at least $\psi-1$ vertices from $H\setminus F'$, hence the number of times we could add vertices to $F'$ is at most $r/((\psi-1)(d-\varphi))$. Since we are choosing vertices from $A\subseteq N(F')$ and $|N(F')|\leq dh$, the number of possible outcomes of the algorithm described above is at most
$$\binom{dh}{\leq r/((\psi-1)(d-\varphi))}= 2^{O\left(\frac{r\log d}{\psi(d-\varphi)}\right)},$$
where the equality comes from Proposition~\ref{prop:Entropy} together with~\eqref{ineq:Uppergt}.

Set $\hat{F}^{*}\coloneqq V(M)\cap F^{*}$. Notice that in each step where we add vertices to $F'$, we add at most one vertex from $\mathcal{L}_{d-1}\setminus H$. So, for each $u\in A$, either $N(u)\cap V(M)=\emptyset$ or we have $d$ choices for $N(u)\cap V(M)$. Therefore, the number of ways to specify $\hat{F}^{*}$ is at most
\[
(d+1)^{r/((\psi-1)(d-\varphi))}=2^{O\left(\frac{r\log d}{\psi(d-\varphi)}\right)}.
\]
After specifying $\hat{F}^{*}$, we set $F''\coloneqq F^{*}\setminus \hat{F}^{*}\subseteq H$. Since $\hat{F}^{*}\cap H\subseteq V(M)\cap H=\emptyset$ and $d_{H\setminus F^{*}}(u)<\psi-1$ for every $u\in A$, we have $d_{F''}(u)+d_{\hat{F}^{*}}(u)> d-\psi+1$. Set $S''\coloneqq\{u\in \mathcal{L}_d: d_{F''}(u)\geq d-\psi\}$. Since $u\in A$ implies $d_{F''}(u)\geq d-\psi$, we have $A\subseteq S''$.

Now we obtain the desired pair $(F'',S'')$ in Step 1 such that $d_{F''}(u)\geq d-\psi$ for every $u\in S''$ and $A\subseteq S''$, $F''\subseteq H$. 

\textbf{Step 2:} Again, we use the linear ordering $\prec_2$ on $\mathcal{L}_{d}\cup \mathcal{L}_{d-1}$. Following $\prec_2$, for each $w\in\mathcal{L}_{d-1}\setminus H $, if $d_{S''}(w)>\psi$, then we update $S''$ with $S''\setminus N(w)$. Note that such $w$ has exactly one neighbor in $A$ when $w\in V(M)$, and no neighbors in $A$ otherwise. So, when $w\in V(M)$, in order to satisfy $A\subseteq S''$, we need to add the vertex in $N(w)\cap V(M)$ back to $S''$. Therefore, for each $w\in \mathcal{L}_{d-1}\setminus H$ we have $d+1$ options: add back either one of its $d$ neighbors or none of them. In the end, 
set $S\coloneqq S''$ and $F\coloneqq F''\cup \{u\in \mathcal{L}_{d-1}: d_{S}(u)>\psi\}$. Notice that for each $u\in S$ we have $d_{F}(u)\geq d_{F''}(u)\geq d-\psi$.

Since $d_{F''}(v)\geq d-\psi$ for every $v\in S''$, we have 
\[
|S''\setminus A|\cdot (d-\psi)\leq e(S''\setminus A, F'')\leq r.
\]
Each time when we remove vertices from $S''$, we remove at least $\psi-1$ vertices from $S''\setminus A$, that is, the guaranteed $\psi$ vertices from $N(w)$ minus the potential vertex from $V(M)$. Therefore, the number of times we could delete vertices from $S''\setminus A$ is at most $r/((\psi-1)(d-\psi))$.

Note that each vertex $w$ is chosen from $N(S'')\subseteq N^2(H)\subseteq N^{3}(A)\subseteq N^4(F'')$. Since $|F''|\leq h$, the number of ways to run Step 2 is at most
\[\binom{d^4 h}{\leq r/((\psi-1)(d-\varphi))}(d+1)^{r/((\psi-1)(d-\varphi))}= 2^{O\left(\frac{r\log d}{\psi(d-\psi)}\right)}.\]\qedhere
\end{proof}

\begin{lemma}\label{lem:boundSF}
   If $(S,F)$ is a $\psi$-approximation for $(A,H)$, then \[|S|\leq |F|+\frac{\psi h-(\psi-1)a}{d-\psi}.\]
\end{lemma}
\begin{proof}
We double-count edges to get
    \[
    (d-\psi)|S\setminus A|+(d-1)|A|\leq e(S,H)\leq d|F|+\psi |H\setminus F|.
    \]
Then,
\[
(d-\psi)|S|+(\psi-1)|A|\leq (d-\psi)|F|+\psi|H|,
\]
which implies the desired result.
\end{proof}
 
\begin{proof}[Proof of Theorem~\ref{thm:SecondPhase}]
Let $I\in \mathcal{J}_1\setminus \mathcal{J}_2$. Recall that $A_{1}=\mathcal{L}_d\cap (I\cup V(M(I)))$, $A_{2}=\mathcal{L}_{d-1}\cap (I\cup V(M(I)))$ and $H_1=N(A_1)\setminus A_2$. Let $a_1\coloneqq |A_1|$ and $h_1\coloneqq |H_1|$, then $h_1\leq |N(A_1)|\leq da_1$. By~\eqref{ineq:sizeA},~\eqref{ineq:upperBoundr} and setting $n=2d-1$ in Lemma~\ref{lem:edgeIso}, we have
\[
2\log^3d\cdot\binom{2d-1}{d}\geq r_1\geq a_1\left(1-\frac{a_1}{\binom{2d-2}{d-1}}\right)\geq \frac{a_1\log^5 d}{d^{5/2}}.
\]
Therefore, we obtain an upper bound on the number of ways of choosing the integers $a_1, h_1$ and $r_1$ by choosing $a_{1}$ with $\left(1-\frac{2\log^3 d}{d}\right)\binom{2d-2}{d-1}\leq a_1\leq \left(1-\frac{\log^5 d}{d^{5/2}}\right)\binom{2d-2}{d-1}$, choosing $h_1$ with $h_1\leq d a_1$ and choosing $r_1$ with $a_1\cdot\log^5 d/d^{5/2}\leq r_1\leq 2\log^3d\cdot\binom{2d-1}{d}$. The number of ways of choosing these integers is at most $2^{O(d)}$, which will be negligible for our purposes.

Recall that $A_1=\cup A_{i}^{1}$, where each $A_{i}^{1}$ is a $2$-linked component of $A_1$, and $L$ is the set of integers $i$ with $|A_{i}^{1}|\leq d^6$. Set $A_{s}\coloneqq\cup_{i\in L}A_{i}^{1}$ and $H_s\coloneqq\cup_{i\in L}H_{i}^{1}$. By~$\eqref{ineq:upperSmallAi}$, to specify $A_s$ we can choose a subset of $\mathcal{L}_d$ of size at most $\frac{20\log^3 d}{d^2}\binom{2d-1}{d}$. To specify $H_s$, we specify $N(A_{s})\setminus H_s$ by choosing a subset of $N(A_s)$ of size at most $|A_s|\leq \frac{20\log^3 d}{d^2}\binom{2d-1}{d}$. Therefore, the number of ways to specify $A_s$ and $H_s$ is at most
\[
\binom{\binom{2d-1}{d}}{\leq \frac{20\log^3 d}{d^2}\binom{2d-1}{d}}^2=2^{O\left(\frac{\log^4 d}{d^2}\binom{2d-2}{d-1}\right)}.
\]
This process also determines $r_{s}\coloneqq|\cup_{i\in L}E(H_{i}^{1},\mathcal{L}_{d}\setminus A_{i}^1)|$.

With the small components determined, we start the process of choosing the large components. Set $a_{m}\coloneqq a_{1}-a_s$, $h_{m}\coloneqq h_{1}-h_s$ and $r_{m}\coloneqq r_{1}-r_{s}$.
Recall that when $a_i^1\geq d^6$, we have $h_i^1\geq r_i^1/d\geq a_i^1\log^5d/d^{7/2}\geq d^{5/2}\log^5d$ by~\eqref{ineq:lowert}. We choose positive integers $a_{i}^{1}$, $h_{i}^{1}$ and $r_{i}^{1}$ such that $$a_{i}^{1}\geq d^6,\quad\quad h_{i}^{1}\geq d^{5/2}\log^5d,\quad\quad r_{i}^{1}\geq \frac{a_i^1\log^5 d}{d^{5/2}}\geq d^{7/2}\log^5 d,$$ $$\sum a_{i}^{1}=a_m,\quad\quad \sum h_{i}^{1}=h_m,\quad\quad \sum r_{i}^{1}=r_m.$$
By Proposition~\ref{prop:decom}, the cost of determining the compositions for each of $a_m$, $h_m$ and $r_m$ is at most 
\[ \frac{a_m}{d^{6}}\log(ed^{6})=O\left(\frac{\log d}{d^{6}}\binom{2d-2}{d-1}\right),\quad \quad
\frac{h_m}{d^{5/2}\log^5d}\log\left(ed^{5/2}\log^5d\right)=O\left(\frac{1}{d^{3/2}\log^4d}\binom{2d-2}{d-1}\right),\]
and\[
\frac{r_m}{d^{7/2}\log^5 d}\log\left(e d^{7/2}\log^5 d\right)=O\left(\frac{1}{d^{7/2}\log d}\binom{2d-2}{d-1}\right),\]
respectively. Thus, the cost of choosing the integers $a_{i}^{1}$, $h_{i}^{1}$ and $r_{i}^{1}$ is $o\left(\frac{1}{d^{3/2}}\binom{2d-2}{d-1}\right)$, which is again negligible for our purposes.

Inequality~\eqref{ineq:upperBoundr} implies $r_1\leq 4\log ^3 d\cdot\binom{2d-2}{d-1}$ and trivially we have $h_1\leq \binom{2d-1}{d}$. For every $i\notin L$, setting $\varphi=d/2$ and using Lemma~\ref{lem:phiaprox}, we obtain a set $F_i'\in \mathcal{V}(a_i,h_i, r_i, \varphi)$ such that $F_{i}'\subseteq H_i^1$; setting $\psi=\sqrt{d} \log d$ and using Lemma~\ref{lem:psiapprox}, we obtain a pair $(S_i,F_i)\in \mathcal{W}(a_i, h_i, r_i, \varphi,\psi)$.

To summarize, up to now, we chose a layer for the larger value between $|A_1|$ and $|A_2|$. We chose integers $a_1,h_1$ and $r_1$, then we chose $A_s$ and $H_s$. For each of $a_m$, $h_m$ and $r_m$, we chose a composition which determines the integers $a_i^{1},h_i^{1}$ and $r_i^{1}$ for $i\notin L$. Finally, for $i\notin L$ we chose a pair $(S_{i},F_i)$. We refer to this process as running the container algorithm. The total cost of running the container algorithm is at most
\[
O\left(\frac{\log^4 d}{d^2}\binom{2d-2}{d-1}\right)+\sum_{i\notin L}\left(O\left(\frac{\log^{2} d}{d^{2}}h_i^{1}\right)+O\left(\frac{\log^{2} d}{d^{2}}r_i^{1}\right)+O\left( \frac{1}{d^{3/2}}r_i^{1}\right)\right)
\]

\[
=O\left(\frac{\log^4 d}{d^{2}}\binom{2d-2}{d-1}\right)+O\left(\frac{\log^{2} d}{d^{2}}h_m\right)+O\left( \frac{1}{d^{3/2}}r_m\right)=O\left(\frac{\log^3 d}{d^{3/2}}\binom{2d-1}{d}\right),
\]
where the last equality uses~\eqref{ineq:upperBoundr}. Let $F\coloneqq\cup_{i\notin L} F_{i}$ and $S\coloneqq\cup_{i\notin L} S_i$. By Lemma~\ref{lem:boundSF}, for every $i\notin L$ we have

\begin{equation}\label{ineq:upperSTi}
|S_i|\leq |F_{i}|+\frac{\psi h_i^1-(\psi-1)a_i^1}{d-\psi}.
\end{equation}
Since the sets $N(A_{i}^{1})$ and the sets $A_i^{1}$ are mutually disjoint and $F_i\subseteq N(A_i^{1})$, summing up~\eqref{ineq:upperSTi} for all $i\notin L$ we get
\begin{equation}\label{ineq:upperSTF}
|S|\leq \sum_{i\notin L} |S_i|\leq \sum_{i\notin L} |F_i|+\frac{\psi \sum_{i\notin L}h_i^1-(\psi-1)\sum_{i\notin L}a_i^1}{d-\psi}=|F|+\frac{\psi h_m-(\psi-1)a_m}{d-\psi}.
\end{equation}

\noindent By~\eqref{ineq:upperBoundr}, since $dh_m-da_m\leq r_m \leq r_1\leq 2\log^3 d\cdot \binom{2d-1}{d}$, we get 
\begin{equation*}
h_m-a_m\leq \frac{2\log^3 d}{d}\binom{2d-1}{d},
\end{equation*} which implies that
\begin{equation*}
\frac{\psi h_m-(\psi-1)a_m}{d-\psi}\leq \frac{a_m}{d-\psi}+\frac{2\psi \log^3 d}{d(d-\psi)}\binom{2d-1}{d}=\frac{a_m}{d}+\frac{\psi a_m}{d(d-\psi)}+\frac{2\psi \log^3 d}{d(d-\psi)}\binom{2d-1}{d}.
\end{equation*}
From~\eqref{ineq:upperSTF}, the fact that $a_m\leq \binom{2d-2}{d-1}$ and $\psi=\sqrt{d} \log d$, we get
\begin{equation*}
    |S|-|F|\leq \frac{a_m}{d}+\frac{\psi a_m}{d(d-\psi)}+\frac{2\psi \log^3 d}{d(d-\psi)}\binom{2d-1}{d}\leq\frac{1}{d}\binom{2d-2}{d-1}+\frac{4\log^{4} d}{d^{3/2}}\binom{2d-1}{d}.
\end{equation*}

\noindent Therefore,
\begin{equation}\label{ineq:upperX*Y*}
|S \cup (\mathcal{L}_{d-1}\setminus F )|=\binom{2d-1}{d}+|S|-|F|\leq 2\binom{2d-2}{d-1}+\frac{8\log^{4} d}{d^{3/2}}\binom{2d-2}{d-1}.
\end{equation}
We split the proof into $3$ cases.

\textbf{Case 1:} If $|S|\leq \left(1-\frac{\log^5 d}{d^{3/2}}\right)\binom{2d-2}{d-1}$, then we choose a subset of $S\cup A_s$ to specify $I\cap \mathcal{L}_{d}$. Since $I$ is a maximal independent set, we can uniquely recover $I\cap \mathcal{L}_{d-1}$ and completely determine $I$. The total cost of choosing the subset and running the container algorithm is at most
\[
\left(1-\frac{\log^5 d}{d^{3/2}}\right)\binom{2d-2}{d-1}+O\left(\frac{\log^3 d}{d^{3/2}}\binom{2d-1}{d}\right)=\left(1-\Omega\left(\frac{\log^ 5 d}{d^{3/2}}\right)\right)\binom{2d-2}{d-1}.
\]

\textbf{Case 2:} If $|\mathcal{L}_{d-1}\setminus F|\leq \left(1-\frac{\log^5 d}{d^{3/2}}\right)\binom{2d-2}{d-1}$, then we can specify $I\cap \mathcal{L}_{d-1}$ by choosing a subset of $\mathcal{L}_{d-1}\setminus F$ and this again uniquely determines $I$. Similarly to the previous case, the total cost is at most $\left(1-\Omega\left(\frac{\log^5 d}{d^{3/2}}\right)\right)\binom{2d-2}{d-1}$.

\textbf{Case 3:}  $|S|, |\mathcal{L}_{d-1}\setminus F|>\left(1-\frac{\log^5 d}{d^{3/2}}\right)\binom{2d-2}{d-1}$. Consider the induced bipartite graph $G'\coloneqq X^{*}\cup Y^{*}$ where $X^{*}\coloneqq S\cup A_s$ and $Y^{*}\coloneqq\mathcal{L}_{d-1}\setminus F$. Let $\prec'$ be the linear ordering on the set of induced matchings in $G'$ induced by the linear ordering $\prec$ on the induced matchings in $B(2d-1,d)$.

Notice that $V(M(I))\cup I\subseteq V(G')$. Therefore, $I$ is a maximal independent set in $G'$ and $M_{G'}(I)=M(I)$. We also have
\begin{align*}\label{ineq:lowerXYPhase2}
\left(1-\frac{\log^5 d}{d^{3/2}}\right) \frac{|X^*\cup Y^*|}{2}\geq \left(1-\frac{\log^5 d}{d^{3/2}}\right)^2 \binom{2d-2}{d-1}\geq \left(1-\frac{2\log^5 d}{d^{3/2}}\right)\binom{2d-2}{d-1}\geq |M_{G'}(I)|.
\end{align*}
Applying Theorem~\ref{thm:Hujter-TuzaStab} to $G'$, we obtain \[
\log \left| \left\{I\in \I(G'): |M_{G'}(I)|\leq \left(1-\frac{\log^5 d}{d^{3/2}}\right) \frac{|X^*\cup Y^*|}{2}\right\} \right|\leq \left(1-\frac{c\log^5 d}{d^{3/2}}\right)\frac{|X^*\cup Y^*|}{2},\]
where $c$ is the constant in Theorem \ref{thm:Hujter-TuzaStab}. By~\eqref{ineq:upperSmallAi} and~\eqref{ineq:upperX*Y*},
\begin{equation*}
\begin{aligned}
    \left(1-\frac{c\log^5 d}{d^{3/2}}\right)\frac{|X^*\cup Y^*|}{2}&\leq  \left(1-\frac{c\log^5 d}{d^{3/2}}\right )\left(\frac{10\log^3d}{d^2}\binom{2d-1}{d}+\left(1+\frac{4\log^4d}{d^{3/2}}\right)\binom{2d-2}{d-1}\right)\\
    &\leq \left(1-\frac{c\log^5 d}{2d^{3/2}}\right)\binom{2d-2}{d-1}.
\end{aligned}
\end{equation*}
The total cost of determining $I$ in this case is the cost of running the container algorithm plus the cost of choosing a maximal independent set from $\{I\in \I(G'): |M_{G'}(I)|\leq \left(1-\log^5 d/d^{3/2}\right) |X^*\cup Y^*|/2\}$, which is at most
\[
 \left(1-\frac{c\log^5 d}{2d^{3/2}}\right)\binom{2d-2}{d-1}+O\left(\frac{\log^3 d}{d^{3/2}}\binom{2d-1}{d}\right)= \left(1-\Omega\left(\frac{\log^5 d}{d^{3/2}}\right)\right)\binom{2d-2}{d-1}.
\]

The combination of these three cases implies that $\log |\mathcal{J}_1\setminus \mathcal{J}_2|=\left(1-\Omega\left(\frac{\log^5 d}{d^{3/2}}\right)\right)\binom{2d-2}{d-1}$, completing the proof of Theorem~\ref{thm:SecondPhase}. 

\end{proof}


\section{For every large induced matching, most of the edges are in one direction}\label{sec:stability}

For this section only, let $$B_1=V(M(I))\cap \mathcal{L}_{d},\quad B_2=V(M(I))\cap \mathcal{L}_{d-1},\quad H_1=N(B_1)\setminus B_2, \quad \text{and}\quad H_2=N(B_2)\setminus B_1.$$ 
Assume $I\in \mathcal{J}_{2}$, then $|B_1|=|B_2|\geq (1-2\log^5 d/d^{3/2})\binom{2d-2}{d-1}$. Recall that in Section~\ref{sec:tools} we defined an adjacent triplet as an ordered triplet of vertices $(x,y,z)\in (\mathcal{L}_{d}\cup \mathcal{L}_{d-1})^{3}$ such that $xy,yz\in E(B(2d-1,d))$. We say an adjacent triplet $(x,y,z)$ is \textit{nice} if $x\in B_1$, $y\in H_1$, $z\in B_1$ or $x\in B_2$, $y\in H_2$, $z\in B_2$, with $x\neq z$ in either case. Denote by $\mathcal{T}_{nice}$ the set of nice triplets.

\begin{claim}\label{claim:lowerNiceTriplets}
For every $I\in \mathcal{J}_2$, we have
    \[
    \begin{split}
    |\mathcal{T}_{\textup{nice}}|&\geq \left(2d(d-1)-24\sqrt{d}\log^5 d\right)\binom{2d-2}{d}.
    \end{split}
    \]
\end{claim}

\begin{proof}   
Choosing a vertex in $H_1$  and two distinct vertices from its neighbors in $B_1$, we can build two nice triplets. We can also choose one vertex in $H_2$ and two of its neighbors from $B_2$ to construct two nice triplets. Therefore, using Jensen's inequality, we have
\begin{equation}\label{eq:lowerTriplets1}
|\mathcal{T}_{\textup{nice}}|=\sum_{v\in H_1}2\binom{d_{B_{1}}(v)}{2}+\sum_{v\in H_2}2\binom{d_{B_{2}}(v)}{2}\geq 2|H_1|\binom{\frac{1}{|H_1|}\sum\limits_{v\in H_1}d_{B_{1}}(v)}{2}+2|H_2|\binom{\frac{1}{|H_2|}\sum\limits_{v\in H_2}d_{B_{2}}(v)}{2}.\end{equation}
Each vertex in $B_1$ has exactly $d-1$ neighbors in $H_1$, so $e(B_1,H_1)=(d-1)|B_1|$. From~\eqref{eq:lowerTriplets1} we obtain
\begin{equation}\label{eq:lowerTriplets}
|\mathcal{T}_{\textup{nice}}|\geq  2|H_1|\binom{\frac{e(B_1,H_1)}{|H_1|}}{2}+2|H_2|\binom{\frac{e(B_2,H_2)}{|H_2|}}{2}= 2|H_1|\binom{\frac{(d-1)|B_1|}{|H_1|}}{2}+2|H_2|\binom{\frac{(d-1)|B_2|}{|H_2|}}{2}.
\end{equation}
Using $\binom{2d-1}{d}=\frac{2d-1}{d-1}\binom{2d-2}{d}$ and $\binom{2d-2}{d-1}=\frac{d}{d-1}\binom{2d-2}{d}$, we have 
\[
|H_1|\leq |\mathcal{L}_{d-1}\setminus B_2|\leq \binom{2d-1}{d}-\left(1-\frac{2\log^5 d}{d^{3/2}}\right)\binom{2d-2}{d-1}\leq \left(1+\frac{3\log^5 d}{d^{3/2}}\right)\binom{2d-2}{d},
\]
which implies
\begin{equation}\label{eq:lowertriplets2}
\frac{(d-1)|B_1|}{|H_1|}\geq \frac{(d-1)(1-2\log^5 d/d^{3/2})\binom{2d-2}{d-1}}{(1+3\log^5 d/d^{3/2})\binom{2d-2}{d}}=\frac{d(1-2\log^5 d/d^{3/2})}{1+3\log^5 d/d^{3/2}}\geq d\left(1-\frac{5\log^5 d}{d^{3/2}}\right).
\end{equation}

We have $d|H_1|\geq e(H_1,B_1)= (d-1)|B_1|$, so $|H_1|\geq (1-2\log^5 d/d^{3/2})\binom{2d-2}{d}$. Similarly, $|H_2|\geq (1-2\log^5 d/d^{3/2})\binom{2d-2}{d}$. We also have the corresponding version of~\eqref{eq:lowertriplets2} for $|H_2|$ and $|B_2|$. Therefore, from~\eqref{eq:lowerTriplets} we conclude
\[
\begin{split}
|\mathcal{T}_{\textup{nice}}|&\geq 4\left(1-\frac{2\log^{5}d}{d^{3/2}}\right)\binom{d(1-5\log^5 d/d^{3/2})}{2}\binom{2d-2}{d}\\
&=\left(2-\frac{4\log^{5}d}{d^{3/2}}\right)\left(d-\frac{5\log^5 d}{\sqrt{d}}\right)\left(d-1-\frac{5\log^5 d}{\sqrt{d}}\right)\binom{2d-2}{d}\\
&= \left(2d-\frac{4\log^5 d}{\sqrt{d}}-\frac{10\log^5 d}{\sqrt{d}}+\frac{20\log^{10}d}{d^2}\right)\left(d-1-\frac{5\log^5 d}{\sqrt{d}}\right)\binom{2d-2}{d}\\
&\geq \left(2d-\frac{14\log^5 d}{\sqrt{d}}\right)\left(d-1-\frac{5\log^5 d}{\sqrt{d}}\right)\binom{2d-2}{d}\\
&= \left(2d(d-1)-\frac{14(d-1)\log^5 d}{\sqrt{d}}-10\sqrt{d}\log^5 d+\frac{70\log^{10}d}{d}\right) \binom{2d-2}{d}\\
&\geq \left(2d(d-1)-24\sqrt{d}\log^5 d\right)\binom{2d-2}{d},
\end{split}
\]
for large enough $d$.\qedhere

\end{proof}

For each $i\in[2d-1]$, let $B_i^{1}\subseteq \mathcal{L}_d$ be the set of vertices $v\in V(M(I))\cap \L_{d}$ such that the edge formed by $v$ and $v\setminus\{i\}$ is in $M(I)$ and $B_i^{2}\subseteq \mathcal{L}_{d-1}$ be the set of vertices $v\in V(M(I))\cap \L_{d-1}$ such that the edge formed by $v$ and $v\cup\{i\}$ is in $M(I)$. Our goal is to show that $|B_i^{1}|=(1-o(1))\binom{2d-2}{d-1}$ for some $i\in [2d-1]$, hence most edges of $M(I)$ are in one fixed direction.


\begin{claim}\label{claim:LowerBadTripletsL}
Let $L\subseteq [2d-1]$. Set $B_{L}^{1}\coloneqq \cup_{i\in L}B_i^{1}$ and
$B_{L}^{2}\coloneqq \cup_{i\in L}B_{i}^{2}$. The number of adjacent triplets $(x,y,z)$ such that $x\in B_{L}^{1}\cup B_{L}^{2}$ and $y,z\in (\mathcal{L}_{d}\cup\mathcal{L}_{d-1})\setminus V(M(I))$ is at least
\[
d|B_{L}^1|-\frac{d}{\binom{2d-2}{d-1}}|B_{L}^1|^2.
\]

\end{claim}

\begin{proof}
By Lemma~\ref{lem:lowerBadTriplets}, we get
\[
|\mathcal{B}(B_{L}^{1})|\geq (2d-1)|B_{L}^{1}|-\frac{d}{\binom{2d-2}{d-1}}|B_{L}^{1}|^2.
\]
From $|\mathcal{B}(B_{L}^{1})|$ we have to subtract the number of triplets $(u,v,w)$ such that $u\in B_{L}^{1}$, $uv\in M(I)$ and $vw\in E(B(2d-1,d))\setminus M(I)$, which is $(d-1)|B_{L}^{1}|$. 

It follows that the number of adjacent triplets $(x,y,z)$ with $x\in B_{L}^{1}, y\notin V(M(I))$ and $z\in \mathcal{L}_{d}\setminus B_{L}^{1}$ is at least
\[
d|B_{L}^{1}|-\frac{d}{\binom{2d-2}{d-1}}|B_{L}^{1}|^2.
\]
We could still have $z\in V(M(I))\setminus B_{L}^{1}$. For $x\in B_{L}^{1}$, let $j\in L$ be the integer with $x(x-j)\in M(I)$ and assume that for some $k\in (x-j)$ and $l\in [2d-1]\setminus x$, the adjacent triplet $(x,x-k,x-k+l)$ satisfies $x\in B_L^1$ and $x-k+l\in V(M(I))\setminus B_L^1$. Then, $(x-k+l)(x-k+l-j)\notin M(I)$. Since $M(I)$ is an induced matching, we have $x-k+l-j\notin V(M(I))\cap \mathcal{L}_{d-1}$.

Therefore, for every adjacent triplet $(x,y,z)$ with $x\in B_{L}^{1}, y\notin V(M(I))$ and $z\in \mathcal{L}_{d}\setminus B_{L}^{1}$, either $x\in B_{L}^{1}$ and $y, z\in (\mathcal{L}_d\cup \mathcal{L}_{d-1})\setminus V(M(I))$, or we injectively obtain the triplet $(x',y',z')=(x-j, x-j+l, x-j+l-k)$ where $x'=x-j\in B_{L}^{2}$ and both $y'=x-j+l$ and $z'=x-j+k-l$ are in $(\mathcal{L}_d\cup\mathcal{L}_{d-1})\setminus V(M(I))$. Note that $(x',y',z')$ is not considered in the previous case. We conclude that the number of adjacent triplets $(x,y,z)$ with $x\in B_{L}^{1}\cup B_{L}^{2}$ and $y,z\in (\mathcal{L}_d\cup \mathcal{L}_{d-1})\setminus V(M(I))$ is at least
\[
d|B_{L}^{1}|-\frac{d}{\binom{2d-2}{d-1}}|B_{L}^{1}|^2,
\]
as desired.
\end{proof}

\begin{claim}\label{claim:UpperNiceTriplets}
For every $I\in \mathcal{J}_2$, 
    \[
    |\mathcal{T}_{\textup{nice}}|\leq 2d(d-1)\binom{2d-2}{d}-\sum_{i\in [2d-1]}d|B_i^{1}|\left(1-\frac{|B_i^{1}|}{\binom{2d-2}{d-1}}\right).
    \]
\end{claim}
\begin{proof}
  Notice that the number of nice triplets is at most the number of adjacent triplets $(x,y,z)$ such that $x\in V(M(I)), y\notin V(M(I))$ and $z\neq x$, excluding the ones such that $z\in (\mathcal{L}_{d}\cup\mathcal{L}_{d-1})\setminus V(M(I))$. Applying Claim~\ref{claim:LowerBadTripletsL} with $B_{L}^{1}=B_i^{1}$ for each $i\in [2d-1]$ and using that $|M(I)|\leq \binom{2d-2}{d-1}$, we get that the number of nice triplets is at most
    \[
    2(d-1)^2\binom{2d-2}{d-1}-\sum_{i\in [2d-1]}d|B_i^{1}|\left(1-\frac{|B_i^{1}|}{\binom{2d-2}{d-1}}\right)=2d(d-1)\binom{2d-2}{d}-\sum_{i\in [2d-1]}d|B_i^{1}|\left(1-\frac{|B_i^{1}|}{\binom{2d-2}{d-1}}\right).
    \]\qedhere
    
\end{proof}

\begin{proof} [Proof of Theorem~\ref{thm:stability}]
Set $b_i^{1}=|B_i^{1}|/\binom{2d-2}{d-1}$. By Claims~\ref{claim:lowerNiceTriplets} and~\ref{claim:UpperNiceTriplets},
\[
\left(2d(d-1)-24\sqrt{d}\log^5 d\right)\binom{2d-2}{d}\leq |\mathcal{T}_{\textup{nice}}|\leq 2d(d-1)\binom{2d-2}{d}-\sum_{i\in [2d-1]}d|B_i^{1}|\left(1-\frac{|B_i^{1}|}{\binom{2d-2}{d-1}}\right).
\]
Hence, using $\binom{2d-2}{d}<\binom{2d-2}{d-1}$, we have
\begin{equation*}
\sum_{i\in[2d-1]}b_i^{1}(1-b_i^{1})\leq \frac{24\log^5 d}{\sqrt{d}}.
\end{equation*}

Notice that $x(1-x)\leq 24\log^5 d/\sqrt{d}$ implies that either $x\leq 25\log^5 d/\sqrt{d}$ or $1-x\leq 25\log^5 d/\sqrt{d}$. Assume for a contradiction that there is no $i\in [2d-1]$ with $b_i^{1}\geq 1-25\log^5 d/\sqrt{d}$. Then, $b_i^{1}\leq 25\log^5 d/\sqrt{d}$ for every $i\in[2d-1]$. Since $I\in \mathcal{J}_2$, we have \begin{equation}\label{eq:sumbi}
    \sum_{i\in[2d-1]} b_{i}^{1}=\frac{|B_1|}{\binom{2d-2}{d-1}}>1-\frac{2\log^5 d}{d^{3/2}}.
\end{equation}
The assumption that $b_i^1\leq 25\log^5 d/\sqrt{d}$ for every $i\in[2d-1]$ together with~\eqref{eq:sumbi} implies that there exists a nonempty set $L\subseteq [2d-1]$ with
\begin{equation}\label{sumbi}
    \sum_{i\in L}b_i^{1}=\frac{1}{2}-O\left(\frac{\log^5 d}{\sqrt{d}}\right).
\end{equation}
By the definition of nice triplets, $|\mathcal{T}_{\text{nice}}|$ is the number of adjacent triplets $(x,y,z)$ with $x\in B_1\cup B_2, y\in (\mathcal{L}_{d}\cup\mathcal{L}_{d-1})\setminus V(M(I))$ and $z\neq x$, excluding the ones with $x\in B_1\cup B_2\supseteq B_{L}^{1}\cup B_{L}^{2}$ and $y,z\in (\mathcal{L}_{d}\cup\mathcal{L}_{d-1})\setminus V(M(I))$. Hence, by Claim~\ref{claim:LowerBadTripletsL}, \[
|\mathcal{T}_{\text{nice}}|\leq 2(d-1)^2\binom{2d-2}{d-1}-d|B_{L}^{1}|\left(1-\frac{|B_{L}^{1}|}{\binom{2d-2}{d-1}}\right)=2d(d-1)\binom{2d-2}{d}-d|B_{L}^{1}|\left(1-\frac{|B_{L}^{1}|}{\binom{2d-2}{d-1}}\right).
\]
Combining with Claim~\ref{claim:lowerNiceTriplets}, we have
\[
\left(2d(d-1)-24\sqrt{d}\log^5 d\right)\binom{2d-2}{d}\leq |\mathcal{T}_{\text{nice}}|\leq 2d(d-1)\binom{2d-2}{d}-d|B_{L}^{1}|\left(1-\frac{|B_{L}^{1}|}{\binom{2d-2}{d-1}}\right).
\]
Consequently,
\[
\frac{|B_{L}^{1}|}{\binom{2d-2}{d}}\left(1-\frac{|B_{L}^1|}{\binom{2d-2}{d-1}}\right)\leq \frac{24\log^5 d}{\sqrt{d}},
\]
which implies
\[
\frac{d}{d-1}\sum_{i\in L}b_i^{1}\left(1-\sum_{i\in L}b_i^{1}\right)\leq \frac{24\log^5 d}{\sqrt{d}},
\]
a contradiction to~\eqref{sumbi}. Therefore, there must exist $i\in [2d-1]$ with $b_i^{1}\geq 1-25\log^ 5 d/\sqrt{d}$. 
\end{proof}


\section{Upper bound in Theorem~\ref{thm:TypicalStructure}}\label{sec:counting}

We will prove our desired upper bound on $|\mathcal{J}_2|$. Recall that for $k\in [2d-1]$, $M_{k}=\{uv: k\in u\in \mathcal{L}_d, v=u\setminus \{k\}\}$ is the canonical matching in direction $k$. By Theorem~\ref{thm:stability}, for every $I\in\mathcal{J}_2$, there is a unique integer $k\in[2d-1]$ such that $I\in \mathcal{U}^{k}$, where $\mathcal{U}^{k}$ denotes the set of independent sets $I\in \mathcal{J}_2$ satisfying $|M(I)\cap M_k|\geq \left(1-\frac{25\log^5 d}{\sqrt{d}}\right)\binom{2d-2}{d-1}$, in particular, $\mathcal{J}_2\subseteq \bigcup_{k} \mathcal{U}^{k}$. We proceed to obtain an upper bound on $|\mathcal{U}^{k}|$, which will be stated in Proposition~\ref{prop:upperUi}. Without loss of generality, we may assume $k=2d-1$. 

Denote $V(B(2d-1,d))\setminus V(M_{2d-1})$ by $V(M_{2d-1})^c$. For $I\in \mathcal{U}^{2d-1}$, let 
\[
\mathcal{A}\coloneqq I\cap V(M_{2d-1})^c= \{v\in I: v\in \mathcal{L}_d \text{ and } 2d-1\notin v; \text{ or } v\in \mathcal{L}_{d-1}\text{ and }2d-1\in v\}.
\]
\noindent Let $\mathcal{A}_{com}\coloneqq  \{A_1,A_2,\ldots \}$ denote the set of $3$-linked components of $\mathcal{A}$ such that either $|A_i|=2$ and $A_i\cap \mathcal{L}_{d}, A_i\cap \mathcal{L}_{d-1}\neq \emptyset$; or $|A_i|\geq 3$.


We provide a procedure to specify every component $A_i\in \mathcal{A}_{\text{com}}$, and use this procedure to generate every maximal independent set $I\in \mathcal{U}^{2d-1}$.

For $D\subseteq V(M_{2d-1})$, define
\[
\widetilde{D}\coloneqq \{v\in V(M_{2d-1}): \{v,w\}\cap D\neq \emptyset, \text{ where }w\text{ is the unique vertex such that }vw\in M_{2d-1}\}.
\]
We will use this $\sim$ notation for any set, even with a different letter. For every component $A_i\in \mathcal{A}_{\text{com}}$, let $A_{i}^{1}\coloneqq A_{i}\cap \mathcal{L}_{d}$ and $A_{i}^{2}\coloneqq A_{i}\cap \mathcal{L}_{d-1}$. Let $G_i\coloneqq N(A_i)$, $G_{i}^{1}\coloneqq N(A_{i}^1)$ and $G_{i}^{2}\coloneqq N(A_{i}^{2})$. Note that $G_i\subseteq V(M_{2d-1})$ since $A_i\subseteq V(M_{2d-1})^c$. Set $\tilde{g}_i\coloneqq |\widetilde{G}_i|$. Then, $\tilde{g}_i/2$ is the number of edges $e\in M_{2d-1}$ such that at least one of the endpoints of $e$ belongs to $N(A_i)$, hence $|M(I)\cap M_{2d-1}|\leq \binom{2d-2}{d-1}-\sum_{i}\tilde{g}_i/2$. The inequality may not be tight since $I$ may have components of size $1$ or $2$ that we ignored. By the definition of $\mathcal{U}^{2d-1}$, we have
\begin{equation}\label{upbd-sumgi}
    \sum_{i}\frac{\tilde{g}_i}{2}\leq \frac{25\log^5d}{\sqrt{d}}\binom{2d-2}{d-1}.
\end{equation}

We say $A_i$ is \textit{small} if $3\leq |A_i|\leq d^{6}$ and \textit{large} if $|A_i|>d^6$.


\subsection{Small components}\label{smallcomp}
We describe a simple procedure to choose the small $3$-linked components of $\mathcal{A}_{\text{com}}$. We establish one more isoperimetric result.

\begin{claim}\label{claim:IsoDoublesided}
     For every $A_i\in \mathcal{A}_{\text{com}}$, we have
    \[
    |\widetilde{G}_i|\geq 2d|A_i|-4|A_i|^2.
    \]
\end{claim}

\begin{proof}
    Let $u,v\in A_i$ be two distinct vertices. If $u,v\in \mathcal{L}_{d}$ or $u,v\in\mathcal{L}_{d-1}$, then $|\widetilde{N(u)}\cap \widetilde{N(v)}|\leq 2$. If $u\in\mathcal{L}_d$ and $v\in\mathcal{L}_{d-1}$, then $|\widetilde{N(u)}\cap \widetilde{N(v)}|\leq 4$. Therefore, $|\widetilde{G}_i|\geq 2d|A_i|-4|A_i|^2$.
\end{proof}

\begin{prop}\label{prop:ChooseSmallComps}
    When $3\leq a_i\leq d^{6}$ and $d$ is sufficiently large, the cost of choosing a $3$-linked set $A_i\subseteq V(M_{2d-1})^c$ with $|A_i|=a_i$ is at most $\tilde{g}_i/2-da_i/19$.
\end{prop}

\begin{proof}
Since $B(2d-1,d)$ is $d$-regular, $A_i$ is $3$-linked, and the number of choices for $v\in V(M_{2d-1})^c$ is $2\binom{2d-2}{d}$, by Lemma~\ref{lem:numcomp}, the number of ways to choose $A_{i}$ is at most 
\[
2\binom{2d-2}{d}d^{9a_i}\leq 2^{2d+ 10 a_i\log d}.
\]

\noindent \textbf{Case:} $3\leq a_i < d/4$. 

By Claim \ref{claim:IsoDoublesided} we have $\tilde{g}_i/2\geq da_i-2a_i^2$. Since $d$ is sufficiently large, the cost of choosing $A_i$ is at most
\begin{equation*}
\frac{\tilde{g}_i}{2}-\frac{\tilde{g}_i}{2}+2d+10a_i\log d\leq \frac{\tilde{g}_i}{2}+2d+10a_i\log d-da_i+2a_i^{2}\leq \frac{\tilde{g}_i}{2}+2d-\frac{3da_i}{4}+2a_i^2\leq \frac{\tilde{g}_i}{2}-\frac{da_i}{15}.
\end{equation*}

\noindent \textbf{Case:} $d/4\leq a_i\leq d^6$. 

Without loss of generality we may assume $|A_i^1|\geq |A_i^2|$. By Lemma~\ref{lem:isoperimetry} we have $|\widetilde{G}_i|/2\geq |N(A_i^1)|\geq d|A_i^1|/9\geq d|A_i|/18$. Since $d$ is sufficiently large, the cost of choosing $A_i$ is at most
\[
\begin{split}\frac{\tilde{g}_i}{2}-\frac{\tilde{g}_i}{2}+2d+ 10a_i\log d&\leq \frac{\tilde{g}_i}{2}+2d+10a_i\log d-\frac{da_i}{18}\leq \frac{\tilde{g}_i}{2}-\frac{da_i}{19}.\qedhere
\end{split}
\]
\end{proof}


\subsection{Large components}\label{largecomp}
Before we outline the procedure for choosing the large components $A_i\in \mathcal{A}_{\text{com}}$, we need to gather several facts about the structure of $A_i$ and $\widetilde{G}_i$, and to stablish another variant of Sapozhenko's container method. Set 
\[
[A_i]\coloneqq \{v\in V(M_{2d-1})^{c}: N(v)\subseteq \widetilde{G}_i\}.
\]
\noindent Since $d|[A_i]|=e([A_i],\widetilde{G}_i)\leq (d-1)|\widetilde{G}_i|$, we have \begin{equation}\label{ineq:lowertildeg_i}
    \tilde{g}_i\geq \frac{d|[A_i]|}{d-1}\geq \frac{d|A_i|}{d-1}>\frac{d^7}{d-1}>d^6.
\end{equation}
For $I\in \mathcal{I}$ and $A_i\in\mathcal{A}_{\text{com}}$, let $I'\coloneqq I\cap (\widetilde{G}_i\cup [A_i])$.

\begin{claim}\label{claim:I is a MIS in GU[A]}
The vertex set $I'$ is a maximal independent set in $B(2d-1,d)[\widetilde{G}_i\cup [A_i]]$.
\end{claim}
\begin{proof}
    Clearly, $I'$ is an independent set. To prove $I'$ is maximal, it suffices to show that $N(v)\cap I'\neq \emptyset$ for every $v\in (\widetilde{G}_i\cup [A_i])\setminus I$. Note that by the maximality of $I$, there is a vertex $u\in N(v)\cap I$.
    
    If $v\in [A_i]\setminus I$, then $u\in N(v)\cap I\subseteq \widetilde{G}_i\cap I$ by the definition of $[A_i]$, which implies that $N(v)\cap I'=N(v)\cap( \widetilde{G}_i\cup [A_i])\cap I\neq \emptyset$.

    Assume that $v\in \widetilde{G}_i\setminus I$. Let $w\in \widetilde{G}_i$ be the unique vertex such that $vw\in M_{2d-1}$, then $N(\{v,w\})\cap A_i\neq \emptyset$ by the definition of $\widetilde{G}_i$. If $N(v)\cap A_i\neq \emptyset$, then $A_i\subseteq [A_i]\cap I$ implies that $N(v)\cap I'=N(v)\cap( \widetilde{G}_i\cup [A_i])\cap I\neq \emptyset$. If there exists $w'\in N(w)\cap A_i$, then $d(w',u)=3$, implying that $u\in A_i$ since $A_i$ is $3$-linked. We conclude again that $N(v)\cap I'=N(v)\cap(\widetilde{G}_i\cup [A_i])\cap I\neq \emptyset$.
\end{proof}

Recall Definition~\ref{def:M(I)}. 
 Fix a linear ordering $\prec'$ on the induced matchings in $B(2d-1,d)[\widetilde{G}_i\cup [A_i]]$ in which $M\prec' M'$ whenever $|M|>|M'|$. When $|M|=|M'|$, then we set $\prec'$ so that $M\prec' M'$ if $|M\cap M_{2d-1}|>|M'\cap M_{2d-1}|$.   

Our objective is to apply a modified version of Sapozhenko's container method on the graph $B(2d-1,d)[\widetilde{G}_i\cup [A_i]]$. To facilitate the process, we establish the following definitions.

Let $B_i\coloneqq[A_i]\cap V(M_{\widetilde{G}_i\cup [A_i]}(I'))$, $B^{1}_i\coloneqq B_i\cap \mathcal{L}_{d}$ and $B^{2}_i\coloneqq B_i\cap \mathcal{L}_{d-1}$. Set $b_i\coloneqq|B_i|$, $b_i^1\coloneqq|B_i^1|$ and $b_i^2\coloneqq|B_i^2|$. Let $Q_i\coloneqq A_i\cup B_i$, $Q_i^1\coloneqq Q_i\cap \mathcal{L}_d$ and $Q_i^2\coloneqq Q_i\cap \mathcal{L}_{d-1}$. Let $H_i^{1}\coloneqq N(Q_i^1)\setminus V(M_{\widetilde{G}_i\cup [A_i]}(I'))$, $H_i^{2}\coloneqq N(Q_i^2)\setminus V(M_{\widetilde{G}_i\cup [A_i]}(I'))$ and $H_i\coloneqq H_i^1\cup H_i^2$. Set $h_i\coloneqq|H_i|$, $h_i^1\coloneqq|H_i^1|$ and $h_i^2\coloneqq|H_i^2|$.

Each edge $e\in E(B(2d-1,d))$ has at least one endpoint in $V(M_{2d-1})$, thus \begin{equation}\label{ineq:simple uppper b_i}
b_i=|M_{\widetilde{G}_i\cup [A_i]}(I')|\leq \frac{\tilde{g}_i}{2}.
\end{equation}

\noindent The value $\tilde{g}_i/2-b_i$ measures how many edges from the maximum of $\tilde{g}_i/2$ we are ``losing" in $M_{\widetilde{G}_i\cup [A_i]}(I')$ by the choice of $A_i$. Note that
\begin{equation*}
    \frac{\tilde{g}_i}{2}-b_i
    =\frac{\tilde{g}_i-2b_i}{2}
    =\frac{\tilde{g}_i/2-2b_i^{1}+\tilde{g}_i/2-2b_i^2}{2}.
\end{equation*}
Since $H_i^{1}=N(Q_i^1)\setminus V(M_{\widetilde{G}_i\cup [A_i]}(I'))$, $H_i^{2}=N(Q_i^2)\setminus V(M_{\widetilde{G}_i\cup [A_i]}(I'))$ and $N(Q_i)\subseteq N(A_i)\cup N([A_i])\subseteq \widetilde{G}_i$, we have
\begin{equation}\label{ineq:upper h-b}
\frac{\tilde{g}_i}{2}-b_i\geq \frac{|N(Q_i^{1})|-2b_i^1+|N(Q_i^{2})|-2b_i^2}{2}
    =\frac{h_i^1-b_i^1+h_i^2-b_i^2}{2}= \frac{h_i-b_i}{2}.
\end{equation}
By ~\eqref{upbd-sumgi},~\eqref{ineq:simple uppper b_i} and Lemma~\ref{lem:lastIso},
\begin{equation*}
|N(Q_i^1)|\geq |N(B_i^1)|\geq \left(2+\frac{\log d}{5d}\right)b_i^1 \quad\text{ and }\quad |N(Q_i^2)|\geq|N(B_i^2)|\geq \left(2+\frac{\log d}{5d}\right)b_i^2,
\end{equation*}
hence
\begin{equation}\label{ineq:lowerNAi}
\left(|N(Q_i^1)|+|N(Q_i^2)|\right)\left(1-\frac{\log d}{11d}\right)\geq  2b_i.
\end{equation}
Together with $h_i^{1}+b_i^{1}=|N(Q_i^{1})|, h_i^{2}+b_i^{2}= |N(Q_i^2)|$ and $|G_i^1|+|G_i^2|\geq \tilde{g}_i/2$, it follows from~\eqref{ineq:lowerNAi} that
\begin{equation}\label{ineq:lower h-b}
\begin{aligned}
   h_i-b_i&=h_i^1-b_i^1+h_i^2-b_i^2=|N(Q_i^1)|-2b_i^1+|N(Q_i^2)|-2b_i^2\\
   &\geq \frac{\log d}{11 d}\left(|N(Q_i^1)|+|N(Q_i^2)|\right) \geq \frac{\log d}{11 d}\left(|N(A_i^1)|+|N(A_i^2)|\right)\geq  \frac{\log d}{22d}\tilde{g}_i.
\end{aligned}
\end{equation}
Let $$r_i^1\coloneqq e(H_i^1,V(M_{2d-1})^c\setminus Q_i^1),\quad \quad r_i^2\coloneqq e(H_i^2,V(M_{2d-1})^c\setminus Q_i^2),\quad \quad r_i\coloneqq e(H_i,V(M_{2d-1})^c\setminus Q_i).$$ From $|Q_i^1|\geq |B_i^1|$, it follows that \begin{equation*}
    r_i^1=(d-1)|H_i^1|-e(H_i^1, Q_i^1)=(d-1)h_i^1-(d|Q_i^1|-b_i^1)\leq (d-1)(h_i^1-b_i^1).
\end{equation*}
Similarly, $r_i^2\leq (d-1)(h_i^2-b_i^2)$, implying
\begin{equation}\label{eq:upboundr_i}
r_i=r_i^1+r_i^2\leq (d-1)(h_i-b_i).
\end{equation}

\begin{claim}\label{claim:quotientrh2}
For every $A_i\in\mathcal{A}_{\text{com}}$ with $|A_i|>d^6$, we have $$\frac{h_i}{r_i}=\frac{|H_i|}{e(H_i,V(M_{2d-1})^c\setminus Q_i)} \leq 2d.$$
\end{claim}

\begin{proof}
     Since $Q^1_i\subseteq \mathcal{L}_{d}$, $H^1_i\subseteq \mathcal{L}_{d-1}$ and $2d-1\notin v$ for every $v\in Q^1_i\cup H_{i}^1$, we may view $Q^1_i$ as a subset of $\mathcal{L}_{d}^{2d-2}$ and $H^1_i$ as a subset of $\mathcal{L}_{d-1}^{2d-2}$. Applying Lemma \ref{lem:edgeIso} with $n=2d-2$, we have
    \[
    r_i^{1}=e_{B(2d-2,d)}(H^{1}_i, \mathcal{L}_{d}^{2d-2}\setminus Q^{1}_i)\geq \frac{d}{d-2}\left(|Q^{1}_i|-\frac{|Q^{1}_i|^2}{\binom{2d-3}{d-1}}\right).
    \]
    Since $N(Q_i^1)\subseteq \widetilde{G}_i$, implying that $|N(Q_i^1)|\leq \tilde{g}_i/2$, by Lemma~\ref{lem:isoperimetry} (iii) we have \begin{equation}\label{upboundq_i1}
        |Q_i^1|\leq \frac{d-1}{d}|N(Q_i^1)|< \frac{\tilde{g}_i}{2}.
    \end{equation}
    Combining with~\eqref{upbd-sumgi} we get \[
    |Q_i^1|\leq \frac{25\log^5d}{\sqrt{d}}\binom{2d-2}{d-1}=\frac{50\log^5d}{\sqrt{d}}\binom{2d-3}{d-1}\leq \frac{d+2}{2d}\binom{2d-3}{d-1},\] which implies that
    \[
    r_i^{1}\geq \frac{d}{d-2}\left(|Q^{1}_i|-\frac{|Q^{1}_i|^2}{\binom{2d-3}{d-1}}\right)\geq\frac{d|Q_i^1|}{d-2}\left(1-\frac{d+2}{2d}\right)=\frac{|Q^{1}_i|}{2}.\]
    Consider the graph induced by taking the complement of every vertex in $\L_d^{2d-1}\cup \L_{d-1}^{2d-1}$, which preserves incidence. Applying the same arguments as above, we have $r_{i}^{2}\geq |Q^{2}_i|/2$, so $r_i=r_i^{1}+r^{2}_{i}\geq |Q_i|/2$. Note that $H_i\subseteq N(Q_i)$, so $|H_i|\leq d|Q_i|$, implying that $h_i/d \leq |Q_i|\leq 2r_i$. We conclude that $h_i/r_i\leq 2d$. 
\end{proof}

Now, we proceed to build the containers $(S_i,F_i)$ for the pairs $(Q_i,H_i)$. The following definition of $\mathcal{G}_{2}(q,h,r,\tilde{g})$ mirrors Definition~\ref{defi:G(A,H)}: Instead of using $A\subseteq \mathcal{L}_d$, we use $Q\subseteq V(M_{2d-1})^c$ now. The key property of $Q$ is $3$-linkedness, which means that we could still use Lemma~\ref{lem:numcomp} when obtaining the $\varphi$-approximation.

\begin{defi}\label{defi:G2(A,H)}
For arbitrary positive integers $q$, $h$, $r$ and $\tilde{g}$,
let $\mathcal{G}_{2}(q,h,r,\tilde{g})\subseteq 2^{V(M_{2d-1})^c}\times 2^{V(M_{2d-1})}$ denote the collection of ordered pairs of sets $(Q,H)$ such that:
\begin{itemize}
\itemsep-0.3em
\item[i)]  $Q\subseteq V(M_{2d-1})^c $, $Q$ is $3$-linked and $|Q|=q$.
\item[ii)] $H\subseteq N(Q)$, $|H|=h$ and if $u\in V(M_{2d-1})^c$ has $N(u)\subseteq H$ then $u\in Q$. 
\item[iii)] $e(H,V(M_{2d-1})^c\setminus Q)=r$.
\item[iv)] Set $Q'\coloneqq N(Q)\setminus H$. There is an induced matching $M\subseteq E(Q, Q')$ which is maximal in $B(2d-1,d)[Q,N(Q)]$, such that $M$ covers $Q'$ and $N(V(M))\cap (Q\setminus V(M))=\emptyset$. 

\item[v)] $|\widetilde{N(Q)}|=\tilde{g}$. 
\end{itemize}
\end{defi}


\begin{defi}
For every $Q\subseteq \mathcal{L}_{d}\cup \mathcal{L}_{d-1}$, $H\subseteq N(Q)$ and $1\leq \varphi\leq d$ , let \[
H^{\varphi}\coloneqq \{v\in S: d_{Q}(v)\geq \varphi\}.\]
\end{defi}

\begin{defi}
For every $(Q,H)\in \mathcal{G}_{2}(q,h,r,\tilde{g})$, a $\varphi$-approximation of $(Q,H)$ is a vertex set $F'$ such that $H^{\varphi}\subseteq F'\subseteq H$ and $Q\subseteq N(F')$.    
\end{defi}

\begin{lemma}\label{lem:phiaprox2}
    For arbitrary positive integers $q, h, r\geq h/(2d), \tilde{g}\geq d^6$ and $100\leq \varphi\leq d$, there is a family $\mathcal{V}_2=\mathcal{V}_2 (q,h,r,\tilde{g},\varphi)\subseteq 2^{V(M_{2d-1})}$ such that every $(Q,H)\in \mathcal{G}_{2}(q,h,r,\tilde{g})$  has a $\varphi$-approximation $F'\in \mathcal{V}_2$ and 

    \[
    \log|\mathcal{V}_2|\leq \begin{cases}
    O\left(\frac{h\log^2 d}{\varphi d}\right)+O\left(\frac{r\log^2 d}{d \varphi }\right)&\text{ if } r\leq \frac{h(d-\varphi)}{\varphi },\\
    O\left(\frac{r\log^2 d}{d(d-\varphi)}\right)+O\left(\frac{r\log^2 d}{d\varphi} \right)&\text{ if } r> \frac{h(d-\varphi)}{\varphi }.
    \end{cases}
    \]
\end{lemma}

\begin{defi}\label{defpsi2} For every $1\leq\psi\leq d$, a $\psi$-approximation for $(Q,H)\in \mathcal{G}_{2}(a,h,r,\tilde{g})$ is a pair $(S,F)$ such that $F\subseteq H$, $Q\subseteq S$ and:
\begin{itemize}
\itemsep-0.3em
    \item[i)] $d_{F}(u)\geq d-\psi$ for every $u\in S$.
    \item[ii)] $d_{\mathcal{L}_d\setminus S}(v)\geq d-\psi$ for every $v\in \mathcal{L}_{d-1}\setminus F$.
\end{itemize}
\end{defi}

\begin{lemma}\label{lem:psiapprox2}
    For arbitrary positive integers $q, h, r\geq h/(2d), \tilde{g}\geq d^6, 100\leq \varphi\leq d$ and $2\leq \psi\leq d$, there is a family $\mathcal{W}_2=\mathcal{W}_2(q,h,r,\tilde{g},\varphi,\psi)$ with
\[
    \log|\mathcal{W}_2|= O\left(\frac{r\log d}{\psi (d-\varphi)}\right)+O\left(\frac{r\log d}{\psi(d-\psi)}\right)
\] satisfying that if $F'\in \mathcal{V}_2(q,h,r,\tilde{g},\varphi)$ is a $\varphi$-approximation of $(Q,H)\in \mathcal{G}_{2}(q,h,r,\tilde{g})$, then there is $(S,F)\in\mathcal{W}_2$ such that $(S,F)$ is a $\psi$-approximation of $(Q,H)$.
\end{lemma}

The proofs for Lemmas~\ref{lem:phiaprox2} and~\ref{lem:psiapprox2} are similar to those of Lemmas~\ref{lem:phiaprox} and~\ref{lem:psiapprox}, respectively, which are provided in Appendix. We also have the following analogue of Lemma~\ref{lem:boundSF}.

\begin{lemma}\label{lem:boundSF2}
   If $(S,F)$ is a $\psi$-approximation for $(Q,H)\in \mathcal{G}_{2}(q,h,r,\tilde{g})$, then 
   \[
   |S|\leq |F|+\frac{\psi (h-q)+q}{d-\psi}.
   \]
\end{lemma}
\begin{proof}
We double-count edges to get
    \[
    (d-\psi)|S\setminus Q|+(d-1)|Q|\leq e(S,H)\leq d|F|+\psi |H\setminus F|.
    \]
Then, 
\[
(d-\psi)|S|+(\psi-1)|Q|\leq (d-\psi)|F|+\psi|H|,
\]
which implies the desired result.
\end{proof}

Now we have all the needed tools to outline our procedure of constructing every large component $A_i\in \mathcal{A}_{\text{com}}$. We choose an approximation $(S_i,F_i)$ for each $(Q_i, H_i)$, with $\varphi_i$ and $\psi_i$ to be determined. Denote $|\widetilde{F}_i|$ by $\tilde{f}_i$. We will split the proof into cases based on the values of $h_i-b_i$, which determines the cost of running the modified Sapozhenko's container method, and $\tilde{g}_i-\tilde{f}_i$, which determines whether it is too costly to recover $\widetilde{G}_i\cup [A_i]$ from $(S_i,F_i)$. Specifically, when $\tilde{g}_i-\tilde{f}_i$ is small, recovering $\widetilde{G}_i\cup [A_i]$ is cheap and we do so. With $\widetilde{G}_i\cup [A_i]$ built, the proof will mirror Sections~\ref{Sec:FirstStage} and~\ref{Subsec:SecondStage}, that is, if $h_i-b_i$ is small, then we use the basic graph container algorithm and apply stability; if $h_i-b_i$ is large, then we take the graph $(\widetilde{G}_i\setminus F_i)\cup ([A_i]\cap S_i)$ and apply stability. On the other hand, when $\tilde{g}_i-\tilde{f}_i$ is large, the cost of recovering $\widetilde{G}_i\cup [A_i]$ is high, in which case our procedure is similar to the proof of Lemma 6.6 in~\cite{Kahn2022}. We will show that there are multiple vertices that are not covered by the induced matching, which will help to reduce the size of the matching enough to apply stability. 

Let $\varepsilon=\varepsilon(d)>0$ be a real number with
\begin{equation}\label{eq:defiEpsilon}
\frac{1}{\log d}\ll \varepsilon \ll 1.
\end{equation}


\noindent \textbf{Case 1: $h_i-b_i\geq \tilde{g}_i/\sqrt{d}$.} Let $\varphi_i=d/2$ and $\psi_i=\sqrt{d}/\log d$. By Lemmas~\ref{lem:phiaprox2} and~\ref{lem:psiapprox2}, we obtain an approximation $(S_i,F_i)$ of $(Q_i, H_i)$, and the cost of choosing $F_i'$ plus the cost of choosing $(S_i,F_i)$ with a fixed $F_i'$ is at most
\begin{equation}\label{costofpaircase1}
\begin{aligned}
     O\left(\frac{h_i\log^2 d}{d^{2}}+\frac{r_i\log^{2} d}{d^2}+\frac{r_i\log^2 d}{d^{3/2}}\right)
     &=O\left(\frac{\tilde{g}_i\log^2d}{d^2}+\frac{(d-1)(h_i-b_i)\log^2 d}{d^{3/2}}\right)\\
     &=O\left(\frac{(h_i-b_i)\log^2 d}{\sqrt{d}}\right),
\end{aligned}
\end{equation}
where the last two equalities follow from that $h_i\leq \tilde{g}_i$,~\eqref{eq:upboundr_i} and the assumption that $h_i-b_i\geq \tilde{g}_i/\sqrt{d}$.
Recalling $\psi_i=\sqrt{d}/\log d$ and $b_i\leq |Q_i|\leq \tilde{g}_i$, by Lemma~\ref{lem:boundSF2} we obtain
\begin{equation}\label{ineq: Case 1 bound SF}
|S_i|-|F_i|\leq \frac{\sqrt{d}(h_i-b_i)}{\log d\cdot(d-\sqrt{d}/\log d)}+\frac{\tilde{g}_i}{d-\sqrt{d}/\log d}\leq \frac{2(h_i-b_i)}{\sqrt{d}\log d}+\frac{2\tilde{g}_i}{d}.
\end{equation}

We split the remaining proof of this case into two subcases.\hfill\break

\noindent \textbf{Case 1.1: }$\tilde{g}_i-\tilde{f}_i\leq 2\varepsilon (h_i-b_i)$. \hfill\break

In this case, we will prove the following upper bound on the cost of specifying $A_i$.

\begin{prop}\label{prop:Case1.1}
    The cost of choosing a large component $A_i$ such that $(Q_i,H_i)\in \mathcal{G}_2(q_i,h_i,r_i,\tilde{g}_i)$, $h_i-b_i\geq \tilde{g}_i/\sqrt{d}$ and $\tilde{g}_i-\tilde{f}_i\leq 2\varepsilon(h_i-b_i)$ is at most $\tilde{g}_i/2-\Omega\left(\tilde{g}_i/2-b_i\right)$. 
\end{prop}
\begin{proof} We will recover $\widetilde{G}_i\cup [A_i]$ from $(S_i,F_i)$, then apply the basic container algorithm in $\widetilde{G}_i\cup [A_i]$ twice, and recover $I'$ using Theorem~\ref{thm:Hujter-TuzaStab}.
\begin{claim}\label{claim: Case 1.1 reconstruction}
    The cost of choosing $\widetilde{G}_i\cup[A_i]$ with a fixed pair $(S_i,F_i)$ is at most $o\left(h_i-b_i\right)$.
\end{claim}
\begin{proof}
    We choose an arbitrary pair $(Q^*, H^{*})\in \mathcal{G}_2(q_i,h_i,r_i,\tilde{g}_i)$ such that $(S_i,F_i)$ is a $\psi_i$-approximation of $(Q^*,H^*)$ at no cost. Note that $F_i\subseteq H^*$ and $Q^*\subseteq S_i$. Denote $\widetilde{G}^*\coloneqq\widetilde{N(Q^*)}$. To specify $\widetilde{G}_i$, we choose the sets $\widetilde{G}_i\setminus \widetilde{G}^{*}$ and $\widetilde{G}^{*}\setminus \widetilde{G}_i$ separately.
    
    Since $\widetilde{G}^*\setminus \widetilde{G}_i\subseteq \widetilde{G}^{*}\setminus \widetilde{F}_i$, the number of ways to choose $\widetilde{G}^{*}\setminus \widetilde{G}_i$ is at most the number of ways to choose a subset of $\widetilde{G}^{*}\setminus \widetilde{F}_i$, the cost of which is at most $\tilde{g}_i-\tilde{f}_i\leq 2\varepsilon(h_i-b_i)=o\left(h_i-b_i\right)$. 

    To specify $\widetilde{G}_i\setminus \widetilde{G}^*$, we choose a vertex set from $S_i\setminus Q^*$ covering $G_i\setminus \widetilde{G}^{*}$. By~\eqref{ineq: Case 1 bound SF} and the assumption that $h_i-b_i\geq \tilde{g}_i/\sqrt{d}$,
    \[
    |S_i\setminus Q^*|= |S_i|-|F_i|+|F_i|-|Q^*|\leq \frac{2(h_i-b_i)}{\sqrt{d}\log d}+\frac{2\tilde{g}_i}{d}+h_i-b_i\leq 2(h_i-b_i).
    \]
    Since $G_i\setminus \widetilde{G}^*\subseteq \widetilde{G}_i\setminus \widetilde{F}_i$, we have $|G_i\setminus \widetilde{G}^*|\leq \tilde{g}_i-\tilde{f}_i\leq 2\varepsilon(h_i-b_i)$. Therefore, by Proposition~\ref{prop:Entropy}, the cost of specifying $\widetilde{G}_i\setminus \widetilde{G}^*$ is at most
    \[
    \log \binom{|S_i\setminus Q^*|}{\leq |G_i\setminus \widetilde{G}^*|}\leq 
    \log \binom{2(h_i-b_i)}{\leq 2\varepsilon (h_i-b_i)}=O\left(H(\varepsilon)(h_i-b_i)\right)=o(h_i-b_i).
    \]
    
    Since from the vertex set $\widetilde{G}_i\subseteq V(M_{2d-1})$ we uniquely recover $[A_i]$, we conclude that the cost of specifying $\widetilde{G}_i\cup[A_i]$ is the cost of choosing $\widetilde{G}_i\setminus \widetilde{G}^{*}$ plus the cost of choosing $\widetilde{G}^{*}\setminus \widetilde{G}_i$, which is at most $o(h_i-b_i)$, as desired.
\end{proof}

Similarly to Section~\ref{Sec:FirstStage}, we will run the basic graph container algorithm in $\widetilde{G}_i\cup [A_i]$. Let $X_0=\widetilde{G}_i\cup [A_i]$. We first run the [Algorithm] for $X_0$ with the stopping condition
\begin{itemize}
    \item $d_{X_j}\leq d^{2/3}$,
\end{itemize}
where $d_{X_j}$ is the maximum degree in $B(2d-1,d)[X_j]$.
Let $Z_1$ be the final $X_j$ and $C_1$ be the certificate. Since each time when we add a vertex to the certificate, we remove at least $d^{2/3}$ vertices from $X_j$,
\begin{equation}\label{upbd:c1}
    |C_1|=|\text{supp}(\xi)|\leq \frac{|X_0|}{d^{2/3}}=\frac{\tilde{g}_i+|[A_i]|}{d^{2/3}}\leq \frac{2\tilde{g}_i}{d^{2/3}},
\end{equation}
where the last inequality uses~\eqref{ineq:lowertildeg_i}. Therefore, the cost of running the [Algorithm] for the first time is at most
\[
\log\binom{|X_0|}{\leq |C_1|}\leq \log\binom{2\tilde{g}_i}{\leq 2\tilde{g}_i/d^{2/3}}=O\left(\frac{\tilde{g}_i\cdot\log d}{d^{2/3}}\right).
\]
Note that the number of edges between $\widetilde{G}_i\cup [A_i]$ and $(\L_d\cup \L_{d-1})\setminus(\widetilde{G}_i\cup [A_i])$ is $(d-1)\tilde{g}_i-d|[A_i]|$.
By Proposition \ref{prop:UpperBoundZ},
\begin{equation}\label{case1.1:ineq:upperZ1}
|Z_1|\leq \frac{d(\tilde{g}_i+|[A_i]|)+(d-1)\tilde{g}_i-d|[A_i]|}{2d-d^{2/3}}\leq \left(1+\frac{d^{2/3}}{2d-d^{2/3}}\right)\tilde{g}_i\leq\left(1+\frac{1}{d^{1/3}}\right)\tilde{g}_i.
\end{equation}

Now we run the basic graph container algorithm for the second time, that is, run the [Algorithm] for $Z_1$ with the following stopping conditions, i.e.~stop if one of them holds:
\begin{itemize}
\itemsep-0.3em
    \item[(a)] $|X_j|\leq 2b_i$;
    \item[(b)] $d_{X_j}\leq d^{1/3}$.
\end{itemize}
Let $Z_2$ be the final $X_j$ and $C_2$ be the certificate. Notice that no matter which stopping condition is applied, we may assume that each time when we add a vertex to the certificate, at least $d^{1/3}$ vertices are removed from $X_j$.

Suppose that the algorithm stops because of condition $(a)$. Then, $|Z_2|\leq 2b_i$ and we may also assume that $|Z_2|\geq2b_i-d$ since $B(2d-1,d)$ is $d$-regular. By~\eqref{case1.1:ineq:upperZ1} and~\eqref{ineq:lowertildeg_i},
\[
|Z_1|-|Z_2|\leq\frac{\tilde{g}_i}{d^{1/3}}+\tilde{g}_i-2b_i+d\leq\frac{2\tilde{g}_i}{d^{1/3}}+\tilde{g}_i-2b_i.
\]
Hence, the number of steps we add vertices to the certificate is
\[
|C_2|=|\text{supp}(\xi)|\leq \frac{|Z_1|-|Z_2|}{d^{1/3}}\leq \frac{2\tilde{g}_i}{d^{2/3}}+\frac{\tilde{g}_i-2b_i}{d^{1/3}},
\]
implying that the cost of running the [Algorithm] for the second time is at most
\[
\log\binom{|Z_1|}{\leq|C_2|}\leq \log \binom{2\tilde{g}_i}{\leq \frac{2\tilde{g}_i}{d^{2/3}}+\frac{\tilde{g}_i-2b_i}{d^{1/3}}}= O\left(\frac{\tilde{g}_i\cdot \log d}{d^{2/3}}\right)+O\left(\frac{(\tilde{g}_i/2-b_i)\log d}{d^{1/3}}\right).
\]
Note that $I'\cap Z_2$ is a maximal independent set in $B(2d-1,d)[Z_2]$. By Theorem~\ref{thm:HujterTuza}, the cost of specifying $I'\cap Z_2$ is at most $|Z_2|/2\leq b_i=\tilde{g}_i/2-(\tilde{g}_i/2-b_i)$. Therefore, the total cost of specifying $I'=I\cap (\widetilde{G}_i\cup [A_i])$ is the sum of the costs of choosing the pair $(S_i,F_i)$, recovering $\widetilde{G}_i\cup [A_i]$, running the basic graph container algorithm twice and recovering $I'$ with a fixed $Z_2$, which is at most
\[
O\left(\frac{(h_i-b_i)\log^2 d}{\sqrt{d}}\right)+o(h_i-b_i)+O\left(\frac{\tilde{g}_i\cdot \log d}{d^{2/3}}\right)+O\left(\frac{(\tilde{g}_i/2-b_i)\log d}{d^{1/3}}\right)+\frac{\tilde{g}_i}{2}-\left(\frac{\tilde{g}_i}{2}-b_i\right)
\]
by~\eqref{costofpaircase1} and Claim~\ref{claim: Case 1.1 reconstruction}. By~\eqref{ineq:upper h-b} and the assumption that $h_i-b_i\geq \tilde{g}_i/\sqrt{d}$, we conclude that the total cost of specifying $I'$ when the algorithm stops because of $(a)$ is at most 
\[\frac{\tilde{g}_i}{2}-\Omega\left(\frac{\tilde{g}_i}{2}-b_i\right).\]

Suppose that the algorithm stops because of condition (b) from now on. We may assume $|Z_2|>2b_i$. By~\eqref{case1.1:ineq:upperZ1},
\[
|Z_1|-|Z_2|\leq \frac{\tilde{g}_i}{d^{1/3}}+\tilde{g}_i-2b_i.
\]
Similarly as above, 
\begin{equation}\label{upbd:c2}
    |C_2|=|\text{supp}(\xi)|\leq \frac{|Z_1|-|Z_2|}{d^{1/3}}\leq \frac{\tilde{g}_i}{d^{2/3}}+\frac{\tilde{g}_i-2b_i}{d^{1/3}},
\end{equation}
hence the cost of running the [Algorithm] for the second time is at most
\[
\log\binom{|Z_1|}{\leq|C_2|}\leq \log \binom{2\tilde{g}_i}{\leq \frac{\tilde{g}_i}{d^{2/3}}+\frac{\tilde{g}_i-2b_i}{d^{1/3}}}=O\left(\frac{\tilde{g}_i\cdot \log d}{d^{2/3}}\right)+O\left(\frac{(\tilde{g}_i/2-b_i)\log d}{d^{1/3}}\right).
\]
Since the algorithm stops because of condition (b), by Proposition~\ref{prop:UpperBoundZ},
\begin{equation}\label{upbd:z2}
    |Z_2|\leq \frac{d(\tilde{g}_i+|[A_i]|)+(d-1)\tilde{g}_i-d|[A_i]|}{2d-d^{1/3}}\leq \left(1+\frac{d^{1/3}}{2d-d^{1/3}}\right)\tilde{g}_i\leq \left(1+\frac{1}{d^{2/3}}\right)\tilde{g}_i.
\end{equation}
Set $Z\coloneqq Z_2\cup C_1\cup C_2$. By~\eqref{upbd:c1},~\eqref{upbd:c2} and~\eqref{upbd:z2},
\begin{equation}\label{ineq:Upper Z case 1.1}
|Z|\leq \tilde{g}_i+\frac{4\tilde{g}_i}{d^{2/3}}+\frac{2(\tilde{g}_i/2-b_i)}{d^{1/3}}.
\end{equation}
Note that $I'\subseteq Z\subseteq \widetilde{G}_i\cup [A_i]$ is a maximal independent set in $B(2d-1,d)[Z]$, so
\[
|M_{Z}(I')|\leq |M_{\widetilde{G}_i\cup [A_i]}(I')|=b_i.
\]

\begin{claim}
The cost of specifying $I'$ from a fixed $Z$ is at most $\tilde{g}_i/2-\Omega\left(\tilde{g}_i/2-b_i\right)$.
\end{claim}

\begin{proof}
If $|Z|\leq \tilde{g}_i-\left(\tilde{g}_i/2-b_i\right)$, then by Theorem~\ref{thm:HujterTuza}, the cost of specifying $I'\subseteq Z$ is at most $|Z|/2\leq \tilde{g}_i/2-\left(\tilde{g}_i/2-b_i\right)/2$, as desired.

We may assume that $|Z|>\tilde{g}_i-\left(\tilde{g}_i/2-b_i\right)$. Let $\varepsilon'\coloneqq(\tilde{g}_i/2-b_i)/\tilde{g}_i$. Then,
\[
\left(1-\varepsilon'\right)\frac{|Z|}{2}>\left(1-\varepsilon'\right)\left(\frac{\tilde{g}_i}{2}-\frac{1}{2}\left(\frac{\tilde{g}_i}{2}-b_i\right)\right)\geq \frac{\tilde{g}_i}{2}-\frac{1}{2}\left(\frac{\tilde{g}_i}{2}-b_i\right)-\frac{1}{2}\left(\frac{\tilde{g}_i}{2}-b_i\right)=b_i\geq |M_{Z}(I')|.
\]
\noindent By Theorem \ref{thm:Hujter-TuzaStab},
\[
\log \left|\left\{I\in \mathcal{I}(Z): |M_{Z}(I)|\leq \left(1-\varepsilon'\right)\frac{|Z|}{2}\right\}\right|\leq \left(1-c\varepsilon'\right)\frac{|Z|}{2}.
\]
Recalling~\eqref{ineq:upper h-b} and the assumption that $h_i-b_i\geq \tilde{g}_i/\sqrt{d}$, we have $\tilde{g}_i/2-b_i\geq (h_i-b_i)/2\geq \tilde{g}_i/(2\sqrt{d})$.
Combining with~\eqref{ineq:Upper Z case 1.1}, we obtain
\[\frac{1}{\log d}\left(\frac{\tilde{g}_i}{2}-b_i\right)\geq\frac{\tilde{g}_i}{4\log d\sqrt{d}}+\frac{\tilde{g}_i/2-b_i}{2\log d}
\geq \frac{2\tilde{g}_i}{d^{2/3}}+\frac{\tilde{g}_i/2-b_i}{d^{1/3}}\geq \frac{|Z|}{2}-\frac{\tilde{g}_i}{2},\]
which implies that (recalling $\varepsilon'=(\tilde{g}_i/2-b_i)/\tilde{g}_i$)
\[
\left(1-c\varepsilon'\right)\frac{|Z|}{2}\leq \left(1-c\varepsilon'\right)\left(\frac{\tilde{g}_i}{2}+\frac{1}{\log d}\left(\frac{\tilde{g}_i}{2}-b_i\right)\right)\leq \frac{\tilde{g}_i}{2}+\frac{1}{\log d}\left(\frac{\tilde{g}_i}{2}-b_i\right)-\frac{c}{2}\left(\frac{\tilde{g}_i}{2}-b_i\right).
\]
Since $d$ is sufficiently large, the cost of choosing $I'\subseteq Z$ is at most
\[
\left(1-c\varepsilon'\right)\frac{|Z|}{2}\leq \frac{\tilde{g}_i}{2}-\frac{c}{4}\left(\frac{\tilde{g}_i}{2}-b_i\right)=\frac{\tilde{g}_i}{2}-\Omega\left(\frac{\tilde{g}_i}{2}-b_i\right).\qedhere
\]
\end{proof}

Similarly to the case where the algorithm stops because of $(a)$, we conclude that when the algorithm stops because of $(b)$, the total cost of specifying $I'=I\cap (\widetilde{G}_i\cup [A_i])$ is the sum of the costs of choosing the pair $(S_i,F_i)$, recovering $\widetilde{G}_i\cup [A_i]$, running the basic graph container algorithm twice and recovering $I'$ with a fixed $Z$, which is at most
    \[O\left(\frac{(h_i-b_i)\log^2 d}{\sqrt{d}}\right)+o(h_i-b_i)+O\left(\frac{\tilde{g}_i\cdot \log d}{d^{2/3}}\right)+O\left(\frac{(\tilde{g}_i/2-b_i)\log d}{d^{1/3}}\right)+\frac{\tilde{g}_i}{2}-\Omega\left(\frac{\tilde{g}_i}{2}-b_i\right)\]
    \[=\frac{\tilde{g}_i}{2}-\Omega\left(\frac{\tilde{g}_i}{2}-b_i\right).\qedhere\]
\end{proof}

\noindent\textbf{Case 1.2:} $\tilde{g}_i-\tilde{f}_i\geq2\varepsilon(h_i-b_i)$.
We postpone and merge it with Case 2.2.\\


\noindent\textbf{Case 2: $h_i-b_i\leq \tilde{g}_i/\sqrt{d}$.}
By~\eqref{ineq:lower h-b}, we have $\tilde{g}_i\log d/(22d)\leq h_i-b_i\leq \tilde{g}_i/\sqrt{d}$, so
\begin{equation}\label{case 2:uppper h-b}
\frac{\log d}{22d}\leq \frac{h_i-b_i}{\tilde{g}_i}\leq \frac{1}{\sqrt{d}}.
\end{equation}
Let $\varphi_i=d/2$ and $\psi_i=\tilde{g}_i/(h_i-b_i)$, then $\psi_i\leq 22d/\log d\leq d/2=\varphi_i$. By Lemmas~\ref{lem:phiaprox2} and~\ref{lem:psiapprox2}, we obtain an approximation $(S_i,F_i)$ of $(Q_i,H_i)$, and the cost of choosing $F'$ plus the cost of choosing $(S_i, F_i)$ with a fixed $F_i'$ is at most
\begin{equation}\label{costofpaircase2}
\begin{aligned}
     O\left(\frac{h_i\log^2 d}{d^{2}}+\frac{r_i\log^{2} d}{d^2}+\frac{r_i\log d(h_i-b_i)}{d\tilde{g}_i}\right)
     &=O\left(\frac{\tilde{g}_i\log^2d}{d^2}+\frac{(d-1)(h_i-b_i)\log d}{d^{3/2}}\right)\\
     &=O\left(\frac{(h_i-b_i)\log d}{\sqrt{d}}\right),
\end{aligned}
\end{equation}
where the last two equalities follow from the definition of $H_i$,~\eqref{eq:upboundr_i} and~\eqref{case 2:uppper h-b}. As $\psi_i=\tilde{g}_i/(h_i-b_i)\leq 22d/\log d$, we have $d-\psi_i\geq d/2$. Together with $b_i\leq |Q_i|\leq \tilde{g}_i$, by Lemma~\ref{lem:boundSF2} we obtain
\begin{equation}\label{ineq:Case 2 upper s-f}
|S_i|-|F_i|\leq \frac{\tilde{g}_i(h_i-b_i)}{(h_i-b_i)(d-\psi_i)}+\frac{\tilde{g}_i}{d-\psi_i}=\frac{2\tilde{g}_i}{d-\psi_i}\leq \frac{4\tilde{g}_i}{d}.
\end{equation}

Similarly to Case 1, we split the remaining proof into two subcases.\\


\noindent \textbf{Case 2.1 :} $\tilde{g}_i-\tilde{f}_i\leq2 \varepsilon (h_i-b_i)$.\\

We will prove the same upper bound on the cost of specifying $A_i$ as in Proposition~\ref{prop:Case1.1}.

\begin{prop}\label{prop:Case2.1}
The cost of choosing a large component $A_i$ such that $(Q_i,H_i)\in \mathcal{G}_2(q_i,h_i,r_i,\tilde{g}_i)$, $h_i-b_i\leq \tilde{g}_i/\sqrt{d}$ and $\tilde{g}_i-\tilde{f}_i\leq 2\varepsilon(h_i-b_i)$ is at most $\tilde{g}_i/2-\Omega\left(\tilde{g}_i/2-b_i\right)$.
\end{prop}

\begin{proof} We will recover $\widetilde{G}_i\cup [A_i]$ from $(S_i,F_i)$ and specify $I'$ by applying Theorem~\ref{thm:Hujter-TuzaStab} in the graph induced by $(\widetilde{G}_i\setminus F_i)\cup(S_i\cap [A_i])$.

\begin{claim}\label{claim: Case 2.1 reconstruction}
    The cost of choosing $\widetilde{G}_i\cup [A_i]$ with a fixed pair $(S_i,F_i)$ is at most $o\left(h_i-b_i\right)$.
\end{claim}

\begin{proof}
    We choose an arbitrary pair $(Q^*, H^{*})\in \mathcal{G}_2(q_i,h_i,r_i,\tilde{g}_i)$ such that $(S_i,F_i)$ is a $\psi_i$-approximation of $(Q^*,H^*)$ at no cost. Note that $F_i\subseteq H^*$ and $Q^*\subseteq S_i$. Denote $\widetilde{G}^*\coloneqq\widetilde{N(Q^*)}$. Similarly to the proof of Claim~\ref{claim: Case 1.1 reconstruction}, to specify $\widetilde{G}_i$, we choose the sets $\widetilde{G}_i\setminus \widetilde{G}^{*}$ and $\widetilde{G}^{*}\setminus \widetilde{G}_i$ separately.

    Since $\widetilde{G}^*\setminus \widetilde{G}_i\subseteq \widetilde{G}^{*}\setminus \widetilde{F}_i$, the number of ways to choose $\widetilde{G}^{*}\setminus \widetilde{G}_i$ is at most the number of ways to choose a subset of $\widetilde{G}^{*}\setminus \widetilde{F}_i$, the cost of which is at most $\tilde{g}_i-\tilde{f}_i\leq 2\varepsilon(h_i-b_i)=o\left(h_i-b_i\right)$. 

    To specify $\widetilde{G}_i\setminus \widetilde{G}^*$, we choose a vertex set from $S_i\setminus Q^*$ covering $G_i\setminus \widetilde{G}^{*}$. By~\eqref{ineq:Case 2 upper s-f} and~\eqref{ineq:lower h-b},
    \[
    |S_i\setminus Q^*|= |S_i|-|F_i|+|F_i|-|Q^*|\leq \frac{4\tilde{g}_i}{d}+h_i-b_i\leq 2(h_i-b_i).
    \]
    Since $G_i\setminus \widetilde{G}^*\subseteq \widetilde{G}_i\setminus \widetilde{F}_i$, we have $|G_i\setminus \widetilde{G}^*|\leq \tilde{g}_i-\tilde{f}_i\leq 2\varepsilon(h_i-b_i)$. Therefore, by Proposition~\ref{prop:Entropy}, the cost of specifying $\widetilde{G}_i\setminus \widetilde{G}^*$ is at most
    \[
    \log \binom{|S_i\setminus Q^*|}{\leq |G_i\setminus \widetilde{G}^*|}\leq 
    \log \binom{2(h_i-b_i)}{\leq 2\varepsilon (h_i-b_i)}=O\left(H(\varepsilon)(h_i-b_i)\right)=o(h_i-b_i).
    \]
    Since from the vertex set $\widetilde{G}_i\subseteq V(M_{2d-1})$ we uniquely recover $[A_i]$, we conclude that the cost of specifying $\widetilde{G}_i\cup[A_i]$ is the cost of choosing $\widetilde{G}_i\setminus \widetilde{G}^{*}$ plus the cost of choosing $\widetilde{G}^{*}\setminus \widetilde{G}_i$, which is at most $o(h_i-b_i)$, as desired.
\end{proof}

Set $Z\coloneqq(\widetilde{G}_i\setminus F_i)\cup ([A_i]\cap S_i)$. Since $F_i\subseteq H_i= N(Q_i)\setminus V(M_{\widetilde{G}\cup [A_i]}(I'))$, we have $F_i\cap I'=\emptyset$. We also have $A_i\subseteq Q_i\subseteq S_i$, implying that $I'=I\cap (\widetilde{G}_i\cup [A_i])\subseteq Z$. Therefore, $I'$ is a maximal independent set in $B(2d-1,d)[Z]$ and
\begin{equation}\label{ineq:upperMICase2.1}
|M_{Z}(I')|\leq |M_{\widetilde{G}_i\cup [A_i]}(I')|=b_i= \frac{\tilde{g}_i}{2}-\left(\frac{\tilde{g}_i}{2}-b_i\right).
\end{equation}
By~\eqref{ineq:Case 2 upper s-f},
\begin{equation}\label{ineq:upperZCase2.1}
|Z|\leq |\widetilde{G}_i|-|F_i|+|S_i|\leq \tilde{g}_i+\frac{4\tilde{g}_i}{d}.
\end{equation}

\begin{claim}\label{chooseI'case2.1}
The cost of choosing $I'$ with a fixed $Z$ is at most $\tilde{g}_i/2-\Omega\left(\tilde{g}_i/2-b_i\right)$.    
\end{claim}

\begin{proof}
If $|Z|\leq \tilde{g}_i-\left(\tilde{g}_i/2-b_i\right)$, then by Theorem \ref{thm:HujterTuza}, the cost of specifying $I'\subseteq Z$ is at most $|Z|/2\leq\tilde{g}_i/2-\left(\tilde{g}_i/2-b_i\right)/2$, as desired.

We may assume that $|Z|>\tilde{g}_i-\left(\tilde{g}_i/2-b_i\right)$.
Let $\varepsilon'\coloneqq(\tilde{g}_i/2-b_i)/\tilde{g}_i$. Then,
\[
\left(1-\varepsilon'\right)\frac{|Z|}{2}>\left(1-\varepsilon'\right)\left(\frac{\tilde{g}_i}{2}-\frac{1}{2}\left(\frac{\tilde{g}_i}{2}-b_i\right)\right)\geq \frac{\tilde{g}_i}{2}-\frac{1}{2}\left(\frac{\tilde{g}_i}{2}-b_i\right)-\frac{1}{2}\left(\frac{\tilde{g}_i}{2}-b_i\right)=b_i\geq |M_{Z}(I')|.
\] By Theorem \ref{thm:Hujter-TuzaStab},
\[
\log \left|\left\{I\in \mathcal{I}(Z): |M_{Z}(I)|\leq \left(1-\varepsilon'\right)\frac{|Z|}{2}\right\}\right|\leq \left(1-c\varepsilon'\right)\frac{|Z|}{2}.
\]
From~\eqref{ineq:upper h-b} and~\eqref{ineq:lower h-b}, it follows that $\tilde{g}_i/2-b_i\geq \tilde{g}_i\log d/(44d)$. Therefore, by~\eqref{ineq:upperZCase2.1}, the cost of specifying $I'\subseteq Z$ from a fixed $Z$ is at most
\begin{equation*}
    \begin{split}
    \left(1-c\varepsilon'\right)\frac{|Z|}{2}
    &\leq\left(\frac{\tilde{g}_i}{2}-\frac{c}{2}\left(\frac{\tilde{g}_i}{2}-b_i\right)\right)\left(1+\frac{4}{d}\right)
    \leq \tilde{g}_i+\frac{4\tilde{g}_i}{d}-\frac{c}{2}\left(\frac{\tilde{g}_i}{2}-b_i\right)\\
    &\leq\tilde{g}_i+\frac{176}{\log d}\left(\frac{\tilde{g}_i}{2}-b_i\right)-\frac{c}{2}\left(\frac{\tilde{g}_i}{2}-b_i\right)=\tilde{g}_i-\Omega\left(\frac{\tilde{g}_i}{2}-b_i\right).\qedhere
    \end{split}
\end{equation*}
\end{proof}

We conclude that the total cost of specifying $A_i$ is the sum of the costs of choosing the pair $(S_i,F_i)$, recovering $\widetilde{G}_i\cup [A_i]$ and choosing $I'$ from $Z=(\widetilde{G}_i\setminus F_i)\cup ([A_i]\cap S_i)$, which is at most
\[
O\left(\frac{(h_i-b_i)\log d}{\sqrt{d}}\right)+o(h_i-b_i)+\frac{\tilde{g}_i}{2}-\Omega\left(\frac{\tilde{g}_i}{2}-b_i\right)=\frac{\tilde{g}_i}{2}-\Omega\left(\frac{\tilde{g}_i}{2}-b_i\right)
\] by~\eqref{costofpaircase2}, Claims~\ref{claim: Case 2.1 reconstruction} and~\ref{chooseI'case2.1}, the assumption that $h_i-b_i\leq \tilde{g}_i/\sqrt{d}$ and~\eqref{ineq:upper h-b}.
\end{proof}

\noindent \textbf{Cases 1.2 and 2.2:} $\tilde{g}_i-\tilde{f}_i\geq 2\varepsilon(h_i-b_i)$.

We now merge the two remaining cases. We specify all the remaining components simultaneously instead of one by one. Let $L_0$ be the set of indices $i$ such that $A_i\in\mathcal{A}_{\text{com}}$ and $|A_i|=2$, and $L_1$ be the set of indices $i$ such that $A_i$ is small. Let $L_2, L_3, L_4$ and $L_5$ be the sets of indices covered by Case~1.1,~1.2,~2.1 and~2.2, respectively. Set $L^*\coloneqq L_0\cup L_1\cup L_2\cup L_4$, then for every $A_i\in\mathcal{A}_{\text{com}}$ with $i\in L^*$, we have provided a procedure to specify $A_i$ in Section~\ref{sec:counting}, thus we may assume that $\cup_{i\in L^*} A_i$ is fixed from now on.

Recall that $(S_i, F_i)$ is a $\psi_i$-approximation of $(Q_i, H_i)$ for every $i\in L_2\cup L_3\cup L_4\cup L_5$. For $I\in \mathcal{J}_2$ set
\[
Z\coloneqq\left(\left(\bigcup_{i\in L_{3}\cup L_{5}}\widetilde{N(S_i)}\setminus F_i\right)\setminus \bigcup_{i\in L^{*}}\widetilde{G}_i\right)\cup \left(\left(\bigcup_{i\in L_3\cup L_5}S_i\right)\setminus\left(\bigcup_{i\in L^*}N(\widetilde{G}_i)\right)\right)
\]
and $I''\coloneqq I\cap Z$. Note that $A_i\subseteq Z$ for $i\in L_3\cup L_5$.
\begin{claim}
The vertex set $I''$ is a maximal independent set in $B(2d-1,d)[Z]$.
\end{claim}
\begin{proof}
As in the proof of Claim~\ref{claim:I is a MIS in GU[A]}, our goal is to show that for every vertex $v\in Z\setminus I''$ we have $N(v)\cap I''\neq \emptyset$.

Note that $N(\widetilde{G}_j)\cap A_i=\emptyset$ when $i\neq j$. We also have $A_i\subseteq Q_i\subseteq S_i$ for every $i\in L_2\cup L_3\cup L_4\cup L_5$. Therefore, for every $j\in L_3\cup L_5$, it follows that $A_j\subseteq \left(\bigcup_{i\in L_3\bigcup L_5}S_i\right)\setminus\left(\bigcup_{i\in L^*}N(\widetilde{G}_i)\right)$.

For every edge $e$ contained in $\left(\cup_{i\in L_3\cup L_5} \widetilde{N(S_i)}\right)\setminus\left(\cup_{i}\widetilde{G}_i\right)$, at most one of its endpoints has a neighbor in $I\cap V(M_{2d-1})^c$. Therefore, one of the endpoints of $e$ belongs to $I$. We obtain that every vertex $v\in \left(\left(\cup_{i\in L_3\cup L_5} \widetilde{N(S_i)}\right)\setminus\left(\cup_{i}\widetilde{G}_i\right)\right)\setminus I$ has a neighbor in $I\cap \left(\left(\cup_{i\in L_3\cup L_5} \widetilde{N(S_i)}\right)\setminus\left(\cup_{i}\widetilde{G}_i\right)\right)$.

For every vertex $v\in\left(\cup_{i\in L_3\cup L_5} \widetilde{G}_i\right)\setminus I$, since $\cup_{i\in L_3\cup L_5} A_i\subseteq (\cup_{i\in L_3\cup L_5}S_i)\setminus \left(\cup_{i\in L^{*}}N(\widetilde{G}_i)\right)$, we have that $N(v)\cap I\cap Z\neq \emptyset$. 

For every vertex $v\in\left(\cup_{i\in L_3\cup L_5} \widetilde{G}_i\right)\setminus I$, assume $v\in \widetilde{G}_{j}$. By the maximality of $I$, there is a vertex $u\in N(v)\cap I$. Let $w\in \widetilde{G}_j$ be the unique vertex such that $vw\in M_{2d-1}$, then $N(\{v,w\})\cap A_j\neq \emptyset$ by the definition of $\widetilde{G}_j$. If $N(v)\cap A_j\neq \emptyset$, then $A_j\subseteq [A_j]\cap I$ implies that $N(v)\cap I''\neq \emptyset$. If there exists $w'\in N(w)\cap A_j$, then $d(w',u)=3$, implying that $u\in A_j$ since $A_j$ is $3$-linked. We conclude again that $N(v)\cap I''\neq \emptyset$.

For every vertex $v\in \left((\cup_{i\in L_3\cup L_5}S_i)\setminus\left(\cup_{i\in L^{*}}N(\widetilde{G}_i)\right)\right)\setminus I
$, by the maximality of $I$, we have $N(v)\cap I\neq \emptyset$. We choose a vertex $v'\in N(v)\cap I$. If $v'\in \cup_{i\in L_3\cup L_5}\widetilde{G}_i$, then, since $F_i\cap I=\emptyset$, it follows that $v'\in \cup_{i\in L_3\cup L_5}\left(\widetilde{G}_i\setminus F_i\right)$ and conclude $N(v)\cap I\cap Z\neq \emptyset$. If $v'\in \cup_{i\in L_3\cup L_5}\left(N(S_i)\setminus \widetilde{G}_i\right)$ then, since we have $v\notin \cup_{i\in L^{*}}N(\widetilde{G}_i)$, it follows that $v'\in Z$  and conclude  $N(v)\cap I\cap Z\neq \emptyset$. 
\end{proof}

Set $\widetilde{X}\coloneqq \cup_{i\in L_3\cup L_5} \widetilde{N(S_i)}\setminus (\cup_{i}\widetilde{G}_i)$ and $\tilde{x}\coloneqq|\widetilde{X}|$. It follows that
\begin{equation}\label{eq:Size of Z Case 2.2}
|Z|\leq\tilde{x}+\sum_{i\in L_3\cup L_5}(\tilde{g}_i-f_i+s_i).
\end{equation}

\begin{prop}\label{prop:Case1.2andCase2.2}
With fixed parameters $(q_i,h_i,r_i,g_i)$ for every $i\in L_{3}\cup L_5$ and fixed compoments $A_i$ for every $i\in L^{*}$, the cost of choosing the $\psi_i$-approximations $(S_i,F_i)$ for $i\in L_{3}\cup L_5$ plus the cost of choosing $I''=I\cap Z$ is at most  
\[
\frac{\tilde{x}}{2}+\sum_{i\in L_3\cup L_5}\left(\frac{\tilde{g}_i}{2}-\Omega\left(\varepsilon(h_i-b_i)\right)\right).
\]
\end{prop}

\begin{proof}
\noindent Recalling the assumption of Case~1, for every $i\in L_3$ we have $\psi_i=\sqrt{d}/\log d$ and $h_i-b_i\geq \tilde{g}_i/\sqrt{d}$. By~\eqref{ineq: Case 1 bound SF},~\eqref{ineq:lower h-b} and~\eqref{eq:defiEpsilon},
\[
s_i-f_i\leq \frac{2(h_i-b_i)}{\sqrt{d}\log d}+\frac{2\tilde{g}_i}{d}\leq \frac{(h_i-b_i)}{\log d}\left(\frac{2}{\sqrt{d}}+44\right)=o(\varepsilon(h_i-b_i)).
\]
Recalling the assumption of Case~2, for every $i\in L_5$ we have $\psi_i=\tilde{g}_i/(h_i-b_i)$. By~\eqref{ineq:Case 2 upper s-f},~\eqref{ineq:lower h-b} and~\eqref{eq:defiEpsilon},
\begin{equation}\label{ineq:Upper s-f case 2.2}
s_i-f_i\leq \frac{4\tilde{g}_i}{d}\leq \frac{88(h_i-b_i)}{\log d}=o\left(\varepsilon(h_i-b_i)\right).
\end{equation}
Therefore, by~\eqref{eq:Size of Z Case 2.2},
\begin{equation}\label{ineq:UpperZ Case 2.2}
|Z|\leq \tilde{x}+\sum_{i\in L_3\cup L_5}\left(\tilde{g}_i+o\left(\varepsilon(h_i-b_i)\right)\right).
\end{equation}
We remark that~\eqref{ineq:Upper s-f case 2.2} is the reason for requiring the lower bound $\varepsilon\gg 1/\log d$ in~\eqref{eq:defiEpsilon}, and all the case analysis and the choices of $\psi_i$ were specifically for~\eqref{ineq:UpperZ Case 2.2} to hold. Now we proceed to establish an upper bound on $|M_{Z}(I'')|$ so that we could apply Theorem~\ref{thm:HujterTuza} as before.

Since $F_i\subseteq H_i\subseteq \widetilde{G}_i$ and all the sets $\widetilde{G}_i$ are pairwise disjoint, we have $\widetilde{G}_i\setminus\widetilde{F}_i= \widetilde{N(Q_i)}\setminus\widetilde{F}_i\subseteq \widetilde{N(S_i)}\setminus\widetilde{F}_i \subseteq Z$ for $i\in L_3\cup L_5$. By the definition of $G_i=N(A_i)$, for each edge $e=uv$ contained in $B(2d-1,d)[\widetilde{G}_i]$, we have $|\{u,v\}\cap N(A_i)|\geq 1$. Since $A_i\subseteq Z$ for $i\in L_3\cup L_5$, no edge contained in $\cup_{i\in L_3\cup L_5} (\widetilde{G}_i\setminus \widetilde{F}_i)\subseteq Z$ can be in $M_{Z}(I'')$.
Hence,
\[
|M_{Z}(I'')|\leq \frac{1}{2}\left(|Z|-\sum_{i\in L_3\cup L_5}\left(\tilde{g}_i-\tilde{f}_i\right)\right).
\]
By~\eqref{ineq:UpperZ Case 2.2} and the assumption of Cases~1.2 and~2.2 that $\tilde{g}_i-\tilde{f}_i\geq2\varepsilon(h_i-b_i)$,
\[
|M_{Z}(I'')|\leq \frac{|Z|}{2}-\sum_{i\in L_3\cup L_5}\varepsilon(h_i-b_i)\leq \frac{\tilde{x}}{2}+\sum_{i\in L_3\cup L_5}\left(\frac{\tilde{g}_i}{2}-\frac{\varepsilon(h_i-b_i)}{2}\right).
\]

\begin{claim}
The cost of choosing $I''$ is at most $\tilde{x}/2+\sum_{i\in L_3\cup L_5}\left(\tilde{g}_i/2-\Omega\left(\varepsilon(h_i-b_i)\right)\right)$.
\end{claim}

\begin{proof}
If $|Z|\leq \tilde{x}+\sum_{i\in L_3\cup L_5}\left(\tilde{g}_i-\varepsilon(h_i-b_i)/2\right)$, then by Theorem~\ref{thm:HujterTuza}, the cost of specifying $I''=I\cap Z$ is at most $|Z|/2\leq\tilde{x}/2+\sum_{i\in L_3\cup L_5}\left(\tilde{g}_i/2-\varepsilon(h_i-b_i)/4\right)$, as desired.

We may assume that $|Z|>\tilde{x}+\sum_{i\in L_3\cup L_5}\left(\tilde{g}_i-\varepsilon(h_i-b_i)/2\right)$. Set \[\varepsilon'\coloneqq\frac{\sum_{i\in L_3\cup L_5}\varepsilon(h_i-b_i)
}{2\left(\tilde{x}+\sum_{i\in L_3\cup L_5}\tilde{g}_i\right)},\] then
\[
(1-\varepsilon')\frac{|Z|}{2}> (1-\varepsilon')\left(\frac{\tilde{x}}{2}+\sum_{i\in L_3\cup L_5}\left(\frac{\tilde{g}_i}{2}-\frac{\varepsilon(h_i-b_i)}{4}\right)\right)\geq \frac{\tilde{x}}{2}+\sum_{i\in L_3\cup L_5}\left(\frac{\tilde{g}_i}{2}-\frac{\varepsilon(h_i-b_i)}{2}\right)\geq |M_{Z}(I'')|.
\]
By Theorem~\ref{thm:Hujter-TuzaStab},
\[
\log \left|\left\{I\in \mathcal{I}(Z): |M_{Z}(I)|\leq \left(1-\varepsilon'\right)\frac{|Z|}{2}\right\}\right|\leq \left(1-c\varepsilon'\right)\frac{|Z|}{2}.
\]
By~\eqref{ineq:UpperZ Case 2.2},
\[
(1-c\varepsilon')\frac{|Z|}{2}\leq (1-c\varepsilon')\left(\frac{\tilde{x}}{2}+\sum_{i\in L_3\cup L_5}\left(\frac{\tilde{g}_i}{2}+o\left(\varepsilon(h_i-b_i)\right)\right)\right) \leq\frac{\tilde{x}}{2}+\sum_{i\in L_3\cup L_5}\left(\frac{\tilde{g}_i}{2}-\frac{c\varepsilon(h_i-b_i)}{2}\right).
\]
Therefore, the cost of specifying $I''=I\cap Z$ is at most $\tilde{x}/2+\sum_{i\in L_3\cup L_5}\left(\tilde{g}_i/2-\Omega\left(\varepsilon(h_i-b_i)\right)\right)$. 
\end{proof}
Recalling~\eqref{costofpaircase1}, ~\eqref{costofpaircase2} and~\eqref{eq:defiEpsilon}, we conclude that the cost of choosing $(S_i,F_i)$ for $i\in L_3\cup L_5$ plus the cost of choosing $I''=I\cap Z$ is at most
\[
O\left(\sum_{i\in L_3}\frac{(h_i-b_i)\log^2 d}{\sqrt{d}}+\sum_{i\in L_5}\frac{(h_i-b_i)\log d}{\sqrt{d}}\right)+\frac{\tilde{x}}{2}+\sum_{i\in L_3\cup L_5}\left(\frac{\tilde{g}_i}{2}-\Omega\left(\varepsilon(h_i-b_i)\right)\right)\]
\[
\leq\frac{\tilde{x}}{2}+\sum_{i\in L_3\cup L_5}\left(\frac{\tilde{g}_i}{2}-\Omega\left(\varepsilon(h_i-b_i)\right)\right).\qedhere
\]
\end{proof}

\subsection{Proof of Theorem~\ref{thm:TypicalStructure}}
Recall that $\mathcal{U}^{k}=\{I\in \mathcal{J}_2: |M(I)\cap M_{k}|\geq \left(1-\frac{25\log^5 d}{\sqrt{d}}\right)\binom{2d-2}{d-1}\}$. Let $\mathcal{U}_{m}^{k}$ be the collection of $I\in \mathcal{U}^{k}$ such that there are exactly $m$ $3$-linked components $A_i\in \mathcal{A}_{\text{com}}$ with $|A_i|=2$. Let $\mathcal{U}_{m,>0}^{k}$ be the collection of $I\in \mathcal{U}_{m}^{k}$ such that there is at least one $A_i\in \mathcal{A}_{\textup{com}}$ with $|A_i|\geq 3$, and $\mathcal{U}_{m,0}^{k}$ be the collection of $I\in \mathcal{U}_{m}^{k}$ such that $|A_i|=2$ for every $A_i\in\mathcal{A}_{\textup{com}}$. Note that~\eqref{upbd-sumgi} implies that $\mathcal{U}_{m}^{k}=\emptyset$ for $m> \frac{25\log^5 d}{(2d-2)\sqrt{d}}\binom{2d-2}{d-1}$. Therefore, to prove Theorem~\ref{thm:main}, it is sufficient to show that $|\mathcal{U}_{m,>0}^{k}|=o\left(|\mathcal{U}_{m,0}^{k}|\right)$ and to provide the right upper bound on $|\mathcal{U}^{k}_{m,0}|$ for each integer $m\leq \frac{25\log^5 d}{(2d-2)\sqrt{d}}\binom{2d-2}{d-1}$.

\begin{prop}\label{prop:upperUi} For every integer $k\in [2d-1]$ we have
\begin{equation}\label{eq:upperUi}
\sum_{m=0}^{\frac{25\log^5 d}{(2d-2)\sqrt{d}}\binom{2d-2}{d-1}}\left|\mathcal{U}_{m,>0}^{k}\right|\leq \exp\left(\frac{(d-1)^2}{2^{2d-1}}\binom{2d-2}{d-1}\right)\cdot 2^{\binom{2d-2}{d-1}-\frac{d}{40}}.
\end{equation}
\end{prop}

\begin{proof}

Without loss of generality, we may assume $k=2d-1$. We will describe a procedure to generate every maximal independent set $I\in \mathcal{U}_{m,>0}^{2d-1}$. First, we choose each component $A_i\in \mathcal{A}_{\text{com}}$ following a procedure according to the cases described in Sections~\ref{smallcomp} and~\ref{largecomp}. 

If the $3$-linked component $A_i\in\mathcal{A}_{\text{com}}$ has $|A_i|=2$, say, $A_i=\{v,u\}$ with $u\in \L_d$ and $v\in \L_{d-1}$, then there are two distinct integers $a,b$ such that $v=(u\cup\{2d-1\})\setminus \{a,b\}$. Therefore, the number of ways to choose all the components $A_i\in\mathcal{A}_{\text{com}}$ with $|A_i|=2$ is at most \[\frac{1}{m!}\left(\binom{2d-2}{d}\frac{d(d-1)}{2}\right)^m=\frac{1}{m!}\left(\binom{2d-2}{d-1}\frac{(d-1)^2}{2}\right)^m.\]
 
Let $a_s\coloneqq\sum_{i}|A_i|$, where the sum is over all integers $i$ such that $A_i\in \mathcal{A}_{\text{com}}$ is small. We specify $a_s$ by choosing an integer between $0$ and $\frac{25\log^5 d}{\sqrt{d}}\binom{2d-2}{d-1}$, after which we choose a composition of it. The number of compositions of $a_s$ is at most $2^{a_s-1}\leq 2^{a_s}$ by Proposition~\ref{prop:decom}. By Proposition~\ref{prop:ChooseSmallComps}, the cost of specifying each of the small components $A_i$ with $|A_i|=a_i$ is at most $\tilde{g}_i/2-da_i/19$.

To specify the large components, let $q_{\ell}\coloneqq\sum_{i}|Q_i|$, where the sum is over all integers $i$ such that $A_i\in \mathcal{A}_{\text{com}}$ is large. Define similarly $h_{\ell}$, $r_{\ell}$ and $\tilde{g}_{\ell}$. We choose a list of integers $(q_{\ell}, h_{\ell}, r_{\ell}, \tilde{g}_{\ell})$ such that
\[
d^6\leq q_\ell\leq \frac{(d-1)\tilde{g}_\ell}{d},\quad
d^5\leq \frac{\tilde{g}_\ell\log d}{22d}\leq h_\ell\leq dq_{\ell}\leq (d-1)\tilde{g}_\ell, \quad
\frac{d^4}{2}\leq \frac{h_\ell}{2d}\leq r_\ell\leq (d-1)h_\ell\leq (d-1)^2\tilde{g}_\ell,
\] and
\[
d^6\leq \tilde{g}_{\ell}\leq \frac{50\log^5 d}{\sqrt{d}}\binom{2d-2}{d-1},
\] where the bounds follow from~\eqref{ineq:lowertildeg_i},~\eqref{ineq:lower h-b},~\eqref{upbd-sumgi} and Claim~\ref{claim:quotientrh2}.
For each integer in the list, we choose a composition $q_{\ell}=\sum_{i}q_i$, $h_\ell=\sum_{i}h_i$, $r_\ell=\sum_i r_{i}$ and $\tilde{g}_\ell=\sum_{i}\tilde{g}_i$ such that the corresponding bounds for $q_i, h_i, r_i$ and $\tilde{g}_i$ hold. By Proposition~\ref{prop:decom}, the cost of determining the composition for each of $q_\ell$, $h_\ell$, $r_\ell$ and $\tilde{g}_\ell$ is at most 
\[ 
\frac{q_\ell}{d^{6}}\log\left(ed^{6}\right)=O\left(\frac{g_\ell\log d}{d^{6}}\right), \quad\quad
\frac{h_\ell}{d^{5}}\log\left(ed^{5}\right)=O\left(\frac{g_\ell\log d}{d^4}\right), \quad\quad
\frac{2r_\ell}{d^{4}}\log\left(\frac{ed^{4}}{2}\right)=O\left(\frac{g_\ell\log d}{d^2}\right),
\] and
\[
\frac{\tilde{g}_\ell}{d^{6}}\log\left(ed^{6}\right)=O\left(\frac{g_\ell\log d}{d^{6}}\right),
\] respectively.
Therefore, for each large component $A_i\in\mathcal{A}_{\text{com}}$, we have chosen an ordered list of integers $(q_i, h_i, r_i, \tilde{g}_i)$. Since $r_i=(d-1)h_i-dq_i+b_i$ by the definition of $r_i$, the value of $b_i$ is uniquely determined. Once $b_i$ and $\tilde{g}_i$ are provided, we obtain $\tilde{g}_i/2-b_i$, which uniquely determines $\varphi_i$ and $\psi_i$. We apply Lemma~\ref{lem:psiapprox2} to get a set $F'_i\in \mathcal{W}$ and Lemma~\ref{lem:phiaprox2} to get a pair $(S_i, F_i)\in \mathcal{V}$, which fixes $|\widetilde{F}_i|=\tilde{f}_i$. Based on the values $h_i-b_i$ and $\tilde{g}_i-\tilde{f}_i$, we obtain the partition $L_0\cup L_1\cup L_2\cup L_3\cup L_4\cup L_5$.  

By Propositions~\ref{prop:Case1.1} and~\ref{prop:Case2.1}, for every $i\in L_2\cup L_4$, the cost of specifying the $3$-linked component $A_i$ is at most $\tilde{g}_i/2-\Omega\left(\tilde{g}_i/2-b_i\right)$. Once we specify all sets $A_i$ for $i\in L_0\cup L_1 \cup L_2 \cup L_4$ and all pairs $(S_i,F_i)$ for $i\in L_3\cup L_5$, by Proposition~\ref{prop:Case1.2andCase2.2}, the cost of choosing every component $A_i$ for $i\in L_3\cup L_5$, the sets $\widetilde{X}$ and $I\cap \widetilde{X}$ is at most
\[
\frac{\tilde{x}}{2}+\sum_{i\in L_3\cup L_5}\left(\frac{\tilde{g}_i}{2}-\Omega\left(\varepsilon (h_i-b_i)\right)\right)\leq \frac{\tilde{x}}{2}+\sum_{i\in L_3\cup L_5}\left(\frac{\tilde{g}_i}{2}-\frac{\tilde{g}_i}{d}\right),
\] where the last inequality follows from~\eqref{ineq:lower h-b} and~\eqref{eq:defiEpsilon}.
Notice that through this process we have also determined $\widetilde{G}_i\cap I$ for every $A_i\in \mathcal{A}_{\text{com}}$.

For each edge $e$ contained in $V(M_{2d-1})\setminus \left(\left(\cup_{i}\widetilde{G}_i\right)\cup  \widetilde{X}\right)$,
exactly one of the endpoints of $e$ belongs to $I$. Choosing one endpoint for each of these edges to be added to $I$ determines $I\cap V(M_{2d-1})$, which allows us to uniquely determine the rest of the maximal independent set $I$. The number of these edges is $\binom{2d-2}{d-1}-\tilde{x}/2-\sum_{i}\tilde{g}_i/2$.

Therefore, the cost of specifying $I$ given the integers $a_i$ for every $i\in L_0\cup L_1$ and the integer lists $(q_i, h_i, r_i,\tilde{g}_i)$ for every $i\in L_2\cup L_3\cup L_4\cup L_5$ is at most
\[
\log\left(\frac{1}{m!}\left(\frac{(d-1)^2}{2}\binom{2d-2}{d}\right)^m\right)+ \binom{2d-2}{d-1}-(2d-2)m-\sum_{i\in L_1}\frac{da_i}{19}-\sum_{i\notin L_0\cup L_1}\Omega\left(\frac{\tilde{g}_i}{d}\right)
\]
\[
=\log\left(\frac{1}{m!}\left(\frac{(d-1)^2}{2^{2d-1}}\binom{2d-2}{d}\right)^m\right)+ \binom{2d-2}{d-1}-\frac{d}{19} a_s-\Omega\left(\frac{\tilde{g}_{\ell}}{d}\right).
\]
We split the process into cases according to the size of $a_s$.

\noindent \textbf{Case: } $3<a_s\leq d^6$. 

The cost of choosing $a_s$ is at most $6\log d$, and the cost of choosing the composition of $a_s=\sum_{i}a_i$ is at most $a_s$. The cost of choosing the integer lists $(q_{\ell}, h_{\ell},r_{\ell}, \tilde{g}_\ell)$ and a corresponding composition for each of them is at most $2^{O\left(\tilde{g}_{\ell}/d^3\right)}$. The cost of choosing $I$ such that the $3$-linked components of $\mathcal{A}$ have the specified sizes is at most
\[
\log\left(\frac{1}{m!}\left(\frac{(d-1)^2}{2^{2d-1}}\binom{2d-2}{d}\right)^m\right)+\binom{2d-2}{d-1}+a_s+6\log d-\frac{d}{19} a_s+O\left(\frac{\tilde{g}_\ell}{d^3}\right)-\Omega\left(\frac{\tilde{g}_{\ell}}{d}\right)
\]
\[
\leq \log\left(\frac{1}{m!}\left(\frac{(d-1)^2}{2^{2d-1}}\binom{2d-2}{d}\right)^m\right)+\binom{2d-2}{d-1}-\frac{d}{20} a_s-\Omega\left(\frac{\tilde{g}_{\ell}}{d}\right).
\]
Therefore, the number of ways to choose $I$ covered by this case is upper bounded by 
\[
\frac{1}{m!}\left(\frac{(d-1)^2}{2^{2d-1}}\binom{2d-2}{d}\right)^m2^{\binom{2d-2}{d-1}-\frac{d}{20}}.
\]

\noindent \textbf{Case: } $d^6\leq a_s$.

The cost of choosing $a_s$ is trivially upper bounded by $\log \left(2\binom{2d-2}{d}\right)<2d$ and the cost of choosing the composition of $a_s$ is at most $a_s$. The cost of choosing the integer lists $(q_{\ell}, h_{\ell},r_{\ell}, \tilde{g}_\ell)$ and a corresponding composition for each of them is at most $2^{O\left(\tilde{g}_{\ell}/d^3\right)}$. The cost of choosing $I$ such that the $3$-linked components of $\mathcal{A}$ have the specified sizes is at most
\[
\log\left(\frac{1}{m!}\left(\frac{(d-1)^2}{2^{2d-1}}\binom{2d-2}{d-1}\right)^m\right)+\binom{2d-2}{d-1}+a_s+3d-\frac{d}{19} a_s-\Omega\left(\frac{\tilde{g}_{\ell}}{d}\right)
\]
\[
\leq \log\left(\frac{1}{m!}\left(\frac{(d-1)^2}{2^{2d-1}}\binom{2d-2}{d-1}\right)^m\right)+\binom{2d-2}{d-1}-\frac{d}{20} a_s-\Omega\left(\frac{\tilde{g}_{\ell}}{d}\right).
\]
Therefore, the number of ways to choose $I$ covered by this case is upper bounded by 
\[
\frac{1}{m!}\left(\frac{(d-1)^2}{2^{2d-1}}\binom{2d-2}{d-1}\right)^m2^{\binom{2d-2}{d-1}-d^5}.
\]

\noindent \textbf{Case: } $a_s=0$. 

We must have $a_{\ell}\neq 0$, hence $\tilde{g}_{\ell}\geq d^6$. Therefore, the number of ways to choose $I$ covered by this case is upper bounded by 
\[
\frac{1}{m!}\left(\frac{(d-1)^2}{2^{2d-1}}\binom{2d-2}{d-1}\right)^m2^{\binom{2d-2}{d-1}-d^4}.
\]

Summing over the three cases above, we obtain
\[
|\mathcal{U}_{m,>0}^{2d-1}|\leq \left(\frac{1}{2^{d/20}}+\frac{1}{2^{d^4}}+\frac{1}{2^{d^5}}\right)\left(\frac{1}{m!}\left(\frac{(d-1)^2}{2^{2d-1}}\binom{2d-2}{d}\right)^m\right)2^{\binom{2d-2}{d-1}}.
\]
We conclude that
\[
\sum_{m=0}^{\frac{25\log^5 d}{(2d-2)\sqrt{d}}\binom{2d-2}{d-1}}|\mathcal{U}_{m,>0}^{2d-1}|\leq \exp \left(\frac{(d-1)^2}{2^{2d-1}}\binom{2d-2}{d}\right) 2^{\binom{2d-2}{d-1}-\frac{d}{40}}.\qedhere
\]
\end{proof}

\begin{proof}[Proof of Theorem~\ref{thm:TypicalStructure}]
    Let $\mathcal{I}^*\subseteq \mathcal{I}$ denote the set of maximal independent sets $I\in \mathcal{I}$ such that there is an integer $k\in [2d-1]$ for which every $3$-linked component of $I\setminus V(M_k)$ has size $1$ or $2$. By Theorems~\ref{thm:FirstPhase} and~\ref{thm:SecondPhase}
    \[
    |\mathcal{I}\setminus \mathcal{I}^*|\leq |\mathcal{J}_2\setminus \mathcal{I}^*|+|\mathcal{I}\setminus \mathcal{J}_2|=|\mathcal{J}_2\setminus \mathcal{I}^*|+|\mathcal{I}\setminus \mathcal{J}_1|+|\mathcal{J}_2\setminus \mathcal{J}_1|=|\mathcal{J}_2\setminus \mathcal{I}^*|+o(|\mathcal{I}|).
    \]
    Since \[
    \mathcal{J}_2\setminus \mathcal{I}^*\subseteq \bigcup_{k\in[2d-1]}\bigcup_{ m=0}^{\frac{25\log^5 d}{(2d-2)\sqrt{d}}\binom{2d-2}{d-1}} \mathcal{U}^{k}_{m,>0},\] 
    it follows from Propositions~\ref{prop:upperUi} and~\ref{prop:LowerBoundMIS} that
    \[
    \left|\bigcup_{k\in[2d-1]}\bigcup_{ m=0}^{ \frac{25 \log^5 d}{(2d-2)\sqrt{d}}\binom{2d-2}{d-1}} \mathcal{U}^{k}_{m,>0}\right|\leq (2d-1)\exp\left(\frac{(d-1)^2}{2^{2d-1}}\right)2^{\binom{2d-2}{d-1}-\frac{d}{40}}=o(|\mathcal{I}|).\qedhere
    \]
\end{proof}

\begin{prop}\label{prop:upperU0}
\begin{equation}\label{eq:upperU0}
    \sum_{m=0}^{\frac{25\log^5d}{(2d-2)\sqrt{d}}\binom{2d-2}{d-1}}| \mathcal{U}_{m,0}^{k}|\leq (1+o(1))\exp\left(\frac{(d-1)^2}{2^{2d-1}}\binom{2d-2}{d-1}\right)\cdot 2^{\binom{2d-2}{d-1}}.
\end{equation}
\end{prop}
\begin{proof}
We repeat the process in the proof of Proposition~\ref{prop:upperUi}. The number of ways to choose the $3$-linked components $A_i\in\mathcal{A}_{\text{com}}$ with $|A_i|=2$ is upper bounded by $\left(\frac{1}{m!}\left(\frac{(d-1)^2}{2}\binom{2d-2}{d-1}\right)^m\right)$. For each edge $e$ contained in $V(M_{k})\setminus \left(\cup_{i}\widetilde{G}_i\right)$,
exactly one of the endpoints of $e$ belongs to $I$. Choosing one endpoint for each of these edges to be added to $I$ determines $I\cap V(M_{k})$, which allows us to uniquely determine the rest of the maximal independent set $I$. The number of these edges is $\binom{2d-2}{d-1}-2m(d-1)$. Therefore,
\[
|\mathcal{U}_{m,0}^{k}|\leq \frac{1}{m!}\left(\binom{2d-2}{d}\frac{d(d-1)}{2}\right)^m 2^{\binom{2d-2}{d-1}-2m(d-1)}=\frac{1}{m!} \left(\binom{2d-2}{d-1}\frac{(d-1)^2}{2^{2d-1}}\right)^m 2^{\binom{2d-2}{d-1}}.
\]
We conclude that
\[
\sum_{m=0}|\mathcal{U}_{m,0}^{k}|\leq \exp\left(\frac{(d-1)^2}{2^{2d-1}}\binom{2d-1}{d-1}\right) \cdot 2^{\binom{2d-2}{d-1}}.\qedhere
\]    
\end{proof}

\begin{proof}[Proof of Theorem~\ref{thm:main}]
    By Theorem~\ref{thm:stability},
    \[
    \mathcal{J}_2\subseteq \bigcup_{k\in [2d-1]}\bigcup_{m=0}^{\frac{25\log^5 d}{(2d-2)\sqrt{d}}\binom{2d-2}{d-1}}\left(\mathcal{U}_{m,0}^{k}\cup \mathcal{U}_{m,>0}^k\right).
    \]
    Using a union bound, we obtain from Propositions~\ref{prop:upperUi} and~\ref{prop:upperU0} that
    \[
    |\mathcal{J}_2|\leq (1+o(1))(2d-1)\exp\left(\frac{(d-1)^2}{2^{2d-1}}\binom{2d-1}{d-1}\right)\cdot 2^{\binom{2d-2}{d-1}},
    \]
    which, together with Theorems~\ref{thm:FirstPhase} and~\ref{thm:SecondPhase}, implies the desired result. 
\end{proof}


\section{Concluding remarks}
We think that the typical structure of $\text{MIS}(B(n,k))$ can be obtained using an equivalent construction to the one in Section~\ref{sec:lowerbound} as long as $n$ and $k$ satisfy $|n/2-k|\leq C\sqrt{n}$ for some fixed constant $C$. On the other hand, for small values of $k$ we have the following observation, where we did not try to optimize the upper bound on $k$.
\begin{prop}
For every $k\leq \log (n/\log^3 n)$,  $\textup{mis}(B(n,k))=(1+o(1))2^{\binom{n}{k-1}}$.
\end{prop}
\begin{proof}
Observe that $\textup{mis}(B(n,k))\leq 2^{\binom{n}{k-1}}$ because for every $A\subseteq \mathcal{L}_{k-1}^{n}$ there is at most one independent set $I\in \mathcal{I}(B(n,k))$ such that $I \cap\mathcal{L}_{k-1}^{n}=A$. From the other side, we show that almost all $A$ has this property. Let $X\subseteq \mathcal{L}_{k-1}^{n}$ be a random set  such that each vertex $v\in \mathcal{L}_{k-1}^{n}$ belongs to $X$ with probability $1/2$ independently. Let $I_{X}\coloneqq X\cup (\mathcal{L}_{k}^{n}\setminus N(X))$. Call a vertex $v\in \mathcal{L}_{k-1}^{n}$ \emph{bad} if $N(v)\subseteq N(X)$ and $v\notin X$. Then,
\[
\mathbb{P}(v \text{ is bad})=\mathbb{P}(N(v)\subseteq N(X) \text{ and }v\notin X)=2^{-1}\cdot (1-2^{-(k-1)})^{n-k+1}\leq e^{-(n-k+1)/2^{k-1}}< e^{-(\log^3 n)/2}.\] 
The set $I_{X}$ is a maximal independent set if and only if there are no bad vertices in $\mathcal{L}_{k-1}^{n}$. Therefore, 
\[
\mathbb{P}(I_{X}\notin \text{MIS}(B(n,k)))\leq \sum_{v\in \mathcal{L}_{k-1}^{n}}\mathbb{P}(v \text{ is bad})\leq \binom{n}{k-1}e^{-(\log^3 n)/2}.
\]
By $\binom{n}{k-1}\leq \binom{n}{k}\leq \left(\frac{en}{k}\right)^{k}$, we conclude that
\[
\mathbb{P}(I_{X}\notin \text{MIS}(B(n,k)))< e^{k\ln(n/k)+k-(\log^3 n)/2}< e^{\log^2 n+\log n-(\log^3 n)/2}=o(1).\qedhere
\]
\end{proof}
We do not propose a conjecture regarding the typical structure for other ranges of $k$. It would be interesting to analyze the transition in the typical structure of maximal independent sets in $B(n,k)$ for the full range of values of $k$. 

Similar tools as the ones present in this article could potentially be applied to obtain precise asymptotics for the number of balanced independent sets in $B(2d-1,d)$. 
It is also possible that our modified version of Sapozhenko's container method can be applied to obtain more precise estimates for the number of maximal independent sets in the hypercube or to obtain upper bounds for the number of maximal independent sets in other families of graphs.

\section*{Acknowledgements}

The authors thank Robert Krueger, Lina Li and Adam Wagner for helpful discussions during the early stages of the project. 

\bibliographystyle{abbrv}
\bibliography{references} 

\begin{thebibliography}{10}

\bibitem{balogh2021independent}
J.~Balogh, R.~I. Garcia, and L.~Li.
\newblock Independent sets in the middle two layers of {B}oolean lattice.
\newblock {\em Journal of Combinatorial Theory, Series A}, 178:105341, 2021.

\bibitem{balogh2015independent}
J.~Balogh, R.~Morris, and W.~Samotij.
\newblock Independent sets in hypergraphs.
\newblock {\em Journal of the American Mathematical Society}, 28(3):669--709, 2015.

\bibitem{BTW}
J.~Balogh, A.~Treglown, and A.~{\relax Zs}. Wagner.
\newblock Applications of graph containers in the {B}oolean lattice.
\newblock {\em Random Structures \& Algorithms}, 49(4):845--872, 2016.

\bibitem{IsoBey}
C.~Bey.
\newblock An upper bound on the sum of squares of degrees in a hypergraph.
\newblock {\em Discrete Mathematics}, 269(1):259--263, 2003.

\bibitem{2025kang}
M.~Collares, J.~Erde, A.~Geisler, and M.~Kang.
\newblock Counting independent sets in expanding bipartite regular graphs.
\newblock {\em arxiv:2503.22255}.

\bibitem{kk-das}
S.~Das.
\newblock Shifting shadows: the {K}ruskal-{K}atona {T}heorem.
\newblock {\em http://discretemath.imp.fu-berlin.de/DMII-2015-16/kruskal.pdf}.

\bibitem{duffus2011maximal}
D.~Duffus, P.~Frankl, and V.~R{\"o}dl.
\newblock Maximal independent sets in bipartite graphs obtained from {B}oolean lattices.
\newblock {\em European Journal of Combinatorics}, 32(1):1--9, 2011.

\bibitem{galvin2011threshold}
D.~Galvin.
\newblock A threshold phenomenon for random independent sets in the discrete hypercube.
\newblock {\em Combinatorics, Probability and Computing}, 20(1):27--51, 2011.

\bibitem{HujterTuza}
M.~Hujter and {\relax Zs}.~Tuza.
\newblock The number of maximal independent sets in triangle-free graphs.
\newblock {\em SIAM J. Discrete Math.}, 6(2):284--288, 1993.

\bibitem{ilinca2013counting}
L.~Ilinca and J.~Kahn.
\newblock Counting maximal antichains and independent sets.
\newblock {\em Order}, 30(2):427--435, 2013.

\bibitem{Jenssen2023Homorphisms}
M.~Jenssen and P.~Keevash.
\newblock Homomorphisms from the torus.
\newblock {\em Adv. Math.}, 430:Paper No. 109212, 89, 2023.

\bibitem{jenssen2024refined}
M.~Jenssen, A.~Malekshahian, and J.~Park.
\newblock A refined graph container lemma and applications to the hard-core model on bipartite expanders.
\newblock {\em arXiv:2411.03393}.

\bibitem{jenssen2019independent}
M.~Jenssen and W.~Perkins.
\newblock Independent sets in the hypercube revisited.
\newblock {\em Journal of the London Mathematical Society}, 102(2):645--669, 2020.

\bibitem{kahn2020stability}
J.~Kahn and J.~Park.
\newblock Stability for maximal independent sets.
\newblock {\em The Electronic Journal of Combinatorics}, pages P1--59, 2020.

\bibitem{Kahn2022}
J.~Kahn and J.~Park.
\newblock The number of maximal independent sets in the {H}amming cube.
\newblock {\em Combinatorica}, 42(6):853--880, 2022.

\bibitem{katona2009theorem}
{\relax Gy}.~Katona.
\newblock A theorem of finite sets.
\newblock In {\em Classic Papers in Combinatorics}, pages 381--401. Springer, 2009.

\bibitem{kleitmanantichain}
D.~Kleitman.
\newblock On {D}edekind's problem: the number of monotone {B}oolean functions.
\newblock {\em Proceedings of the American Mathematical Society}, 21(3):677--682, 1969.

\bibitem{korshunov1983number}
A.~Korshunov and A.~Sapozhenko.
\newblock The number of binary codes with distance 2.
\newblock {\em Problemy Kibernet}, 40:111--130, 1983.

\bibitem{krueger2024lipschitz}
R.~A. Krueger, L.~Li, and J.~Park.
\newblock Lipschitz functions on weak expanders.
\newblock {\em arXiv:2408.14702}.

\bibitem{kruskal1963number}
J.~B. Kruskal.
\newblock The number of simplices in a complex.
\newblock {\em Mathematical Optimization Techniques}, 10:251--278, 1963.

\bibitem{li2025number}
L.~Li, G.~McKinley, and J.~Park.
\newblock The number of colorings of the middle layers of the hamming cube.
\newblock {\em Combinatorica}, 45(1):1--47, 2025.

\bibitem{lovasz1975ratio}
L.~Lov{\'a}sz.
\newblock On the ratio of optimal integral and fractional covers.
\newblock {\em Discrete Mathematics}, 13(4):383--390, 1975.

\bibitem{lovasz2007combinatorial}
L.~Lov{\'a}sz.
\newblock {\em Combinatorial problems and exercises}, volume 361.
\newblock American Mathematical Soc., 2007.

\bibitem{potukuchi2021abelian}
A.~Potukuchi and L.~Yepremyan.
\newblock Enumerating independent sets in {A}belian {C}ayley graphs.
\newblock {\em arXiv:2109.06152}.

\bibitem{sapozhenko1987}
A.~A. Sapozhenko.
\newblock On the number of connected subsets with given cardinality of the boundary in bipartite graphs.
\newblock {\em Metody Diskret. Analiz}, 45(45):42--70, 1987.

\bibitem{saxton2015hypergraph}
D.~Saxton and A.~Thomason.
\newblock Hypergraph containers.
\newblock {\em Inventiones Mathematicae}, 201(3):925--992, 2015.

\bibitem{stein1974two}
S.~K. Stein.
\newblock Two combinatorial covering theorems.
\newblock {\em Journal of Combinatorial Theory, Series A}, 16(3):391--397, 1974.

\end{thebibliography}

\appendix

\section{Proofs of the $3$-linked container lemmas}

\begin{proof}[Proof of Lemma~\ref{lem:phiaprox2}]
    Let $p=40\frac{\log d}{\varphi d}$. Let $X$ be a random subset of $H$, where each vertex in $H$ belongs to $X$ independently with probability $p$. Let $\Omega_{X}\coloneqq E(X, V(M_{2d-1})^c\setminus Q)$. Notice that $\mathbb{E}(|X|)=ph$ and $\mathbb{E}(|\Omega_X|)=pr$. Set \[
    T_{X}\coloneqq\{v\in V(M_{2d-1}): |N_{Q}(v)\cap N_{Q}(X)|\geq 2\},
    \]
    then we have $T_{X}\subseteq N(Q)$. Let $M\subseteq E(Q,N(Q))$ be the induced matching covering $N(Q)\setminus H$ guaranteed by the definition of $\mathcal{G}_2(q,h,r,\tilde{g})$. Then, for each $v\in V(M)\cap V(M_{2d-1})$ we have $|N_Q(v)|=1$, which implies that $T_{X}\cap V(M)=\emptyset$. We conclude that $T_{X}\subseteq N(Q)\setminus V(M)=H$. 
    
    Since $|N_Q(y)|\geq \varphi$ for every $y\in H^{\varphi}$, by Lemma~\ref{lem:isoperimetry} and the vertex-transitivity of $B(2d-1,d)$, we have $|N(N_{Q}(y))|\geq \varphi d/9$, which implies
    \[
    |N(N_{Q}(y))\cap H|\geq \varphi d/9-|N_{Q}(y)|\geq \varphi d/9-d\geq \varphi d/10.
    \]
    Furthermore, if $N_{Q}(y)\cap N_{Q}(X)=\{w\}$ for $y\in H^{\varphi}$, then $|N_{Q}(y)\setminus\{w\}|\geq \varphi-1$ and $N_{Q}(y)\setminus\{w\}$ has no neighbors in $X$. Similarly, we obtain \[
    |N(N_{Q}(y)\setminus\{w\})\cap H|\geq (\varphi-1)d/9-|N_{Q}(y)\setminus\{w\}|\geq \varphi d/10.\]
    It follows that for every $y\in H^{\varphi}$, we have
    \begin{align*}
    \mathbb{P}(y\notin T_X)&\leq \mathbb{P}(|N_{Q}(y)\cap N_{Q}(X)|=0)+\mathbb{P}(|N_{Q}(y)\cap N_{Q}(X)|=1)\\
    &\leq (1-p)^{\varphi d/10}+d(1-p)^{\varphi d/10}\leq  (d+1)\exp (-40\log d /10)\leq 2d^{-3},
    \end{align*}
    which implies
    \[
    \mathbb{E}(|H^{\varphi}\setminus T_{X}|)\leq  \frac{2h}{d^{3}}.
    \]

  \noindent By Markov's inequality there are vertex sets $T_0\coloneqq X\subseteq H$, $T_1\coloneqq T_X\subseteq H$ and $T_2\coloneqq  H^{\varphi}\setminus (T_0\cup T_1)$ satisfying
    \begin{equation}\label{sec6:uppert0}
    |T_0|\leq \frac{200h\log d}{\varphi d},
    \end{equation}
    \begin{equation}\label{sec6:upperOmegat0'}
    e(T_0, V(M_{2d-1})^c\setminus Q)\leq \frac{200r\log d}{\varphi d},
    \end{equation}and
    \begin{equation}\label{sec6:uppert0'}
    |T_2|\leq |H^{\varphi}\setminus T_1|\leq \frac{10h}{d^{3}}.
    \end{equation}

Let $T^{*}\coloneqq T_0\cup T_1\cup T_2\supseteq H^{\varphi}$. Up to now, if we are provided $T_0$, $\Omega_{T_0}$ and $T_2$, then we can recover $T_1$ and build $T^{*}$ such that $H^{\varphi}\subseteq T^{*}\subseteq H$. This is almost our $\varphi$-approximation, the set $F'$. We are missing the property $A\subseteq N(F')$.

Let $T_3\subseteq H\setminus T^{*}$ be a minimal cover of $Q\setminus N(T^{*})$ and set 
\[
F'\coloneqq T^{*}\cup T_{3}=T_0\cup T_1\cup T_2\cup T_3,
\] 
which is our $\varphi$-approximation. Notice that $T_{i}\subseteq H$ for every $i$, thus $F'\subseteq H$.

Now we proceed to get an upper bound on the number of possible choices of $F'$. 
Let 
\[
    T=T_0\cup T_2\cup T_3\subseteq N(Q).
\]
\noindent We have $d(u,T)\leq 3$ for every $u\in Q$ by the definitions of $T_3$ and $T_1$, and $d(u,Q)=1$ for every $u\in T$. Using Lemma~\ref{lem:linked}, it follows that $T$ is $9$-linked.

For every $T_0$ and $\Omega_{T_0}=E(T_0,V(M_{2d-1})^c\setminus Q)$, we can determine $N_{Q}(T_0)$ and $T_1$. Hence, if we are provided with $T$ and $\Omega_{T_0}$, then we can determine $F'$. 

To get an upper bound on $|T_3|$, notice that $T_3\subseteq H\setminus T^{*} \subseteq H\setminus H^{\varphi}$. We have 
\[
|H\setminus T^{*}|\leq \frac{e(H\setminus T^{*},V(M_{2d-1})^c\setminus Q )}{d-\varphi}\leq \frac{e(H,V(M_{2d-1})^c\setminus Q)}{d-\varphi}=\frac{r}{d-\varphi}.
\] 
For $i\in\{1,2\}$ we also have $d_{H_i\setminus L}(u)\geq d-1$ for every $u\in Q_i\setminus N(T^{*})$ and trivially $d_{Q_i\setminus N(T^{*})}(u)\leq d$ for all $u\in H_i\setminus T^{*}$. Applying Theorem~\ref{thm:cover} twice,  setting $P=Q_i\setminus N(T^*)$, $Q=H_i\setminus T^{*}$, $a=(d-1)$, $b=d$ for each $i\in\{1,2\}$. Summing what we obtain, we get 
    \begin{equation}\label{sec6:uppert1}
    |T_3|\leq  \frac{|H_1\setminus T^{*}|}{(d-1)}(1+\ln d)+\frac{|H_2\setminus T^{*}|}{(d-1)}(1+\ln d)=\frac{|H\setminus T^{*}|}{(d-1)}(1+\ln d) \leq  \frac{2r\log d}{d(d-\varphi)}.
    \end{equation}

\noindent From~\eqref{sec6:uppert0},~\eqref{sec6:uppert0'} and~\eqref{sec6:uppert1}, it follows that
\begin{equation}\label{sec6:upperT}
|T|=|T_0|+|T_2|+|T_3|\leq \frac{200 h\log d}{\varphi d}+\frac{10h}{d^{3}}+\frac{2r\log d}{d(d-\varphi)}.
\end{equation}

\noindent With a fixed $T_0$, since $\Omega_{T_0}\subseteq E(T_0, V(M_{2d-1})^{c})$, from~\eqref{sec6:uppert0} and~\eqref{sec6:upperOmegat0'} it follows that the number of choices of $\Omega_{T_0}$ is at most 
$$\binom{d|T_0|}{\leq \frac{200r\log d}{\varphi d} }\leq\binom{200h\log d/\varphi}{\leq \frac{200r\log d}{\varphi d}}= 2^{O\left(\frac{r\log^2 d}{\varphi d}\right)},$$ 
where we used Proposition~\ref{prop:Entropy} and Claim \ref{claim:quotientrh2} to obtain the last equality.

To bound the cost of specifying $T$, we split the proof into two cases.
    
\noindent $\bullet$ If $r\leq \frac{h(d-\varphi)}{\varphi }$: From~\eqref{sec6:upperT} we get $|T|\leq\frac{203h\log d}{\varphi d}$. Recall that $T$ is $9$-linked. By Corollary~\ref{cor:klinkednumber}, the number of choices of $T$ is $\binom{2d-1}{d}\cdot 2^{O\left(\frac{h\log^2 d}{\varphi d}\right)}=2^{O\left(\frac{h\log^2 d}{\varphi d}\right)}$. Since $T=T_0\cup T_2\cup T_3$, we have at most $3^{|T|}<2^{\frac{406h\log d}{\varphi d}}$ choices of $T_0,T_2$ and $T_3$. Given $T_0$, $T_2$, $T_3$ and $\Omega_{T_0}$, we can reconstruct $F'$. Thus, the cost of specifying $F'$ is
$$O\left(\frac{h\log^2 d}{\varphi d}\right)+O\left(\frac{r\log^2 d}{d \varphi }\right).$$

\noindent$\bullet$ If $r>\frac{h(d-\varphi)}{\varphi }$: From~\eqref{sec6:upperT} we get $|T|\leq \frac{203r\log d}{d(d-\varphi)}$. Recall that $T$ is $9$-linked. By Corollary~\ref{cor:klinkednumber}, the number of choices of $T$ is $\binom{2d-1}{d}\cdot 2^{O\left(\frac{r\log^2 d}{d(d-\varphi)}\right)}=2^{O\left(\frac{r\log^2 d}{d(d-\varphi)}\right)}$. Since $T=T_0\cup T_2\cup T_3$, we have at most $3^{|T|}<2^{\frac{406r\log d}{d(d-\varphi)}}$ choices of $T_0,T_2$ and $T_3$.  As in the previous case, the cost of choosing $\Omega_{T_0}$ is $O\left(\frac{r\log^2 d}{\varphi d}\right)$. Thus, in this case, the cost of specifying $F'$ is
\[O\left(\frac{r\log^2 d}{d(d-\varphi)}\right)+O\left(\frac{r\log^2 d}{d\varphi} \right).\qedhere\]
    
\end{proof}

\begin{proof}[Proof of Lemma~\ref{lem:psiapprox2}]

We split the proof into two steps. In the first step we make sure that the condition (i) in Definition~\ref{defpsi2} is satisfied, while in the second step we make sure that condition (ii) is satisfied.

\textbf{Step 1:} Fix an arbitrary linear order $\prec_2$ on $\mathcal{L}_{d}\cup \mathcal{L}_{d-1}$. Following the ordering $\prec_2$, for each $u\in Q$, if $d_{H\setminus F'}(u)\geq\psi-1$, then update $F'$ to $F'\cup N(u)$. We stop when there is no $u\in Q$ with $d_{H\setminus F'}(u)\geq\psi-1$, i.e., for every $u\in Q$ we have $d_{H\setminus F'}(u)<\psi-1$. At the end we set $F^{*}=F'$.

Notice that each $v\in H\setminus F'\subseteq H\setminus H^{\varphi}$ satisfies $d_{V(M_{2d-1})^c\setminus Q}(v)> d-1-\varphi$, thus 
\[(d-1-\varphi)|H\setminus F'|\leq e(H\setminus F', V(M_{2d-1})^c\setminus Q)\leq r,\]
which implies $|H\setminus F'|\leq r/(d-1-\varphi)$. In each step where we add vertices to $F'$, we remove at least $\psi-1$ vertices from $H\setminus F'$, hence the number of times we could add vertices to $F'$ is at most $ r/((\psi-1)(d-1-\varphi)) $. Since we are choosing vertices from $Q\subseteq N(F')$ and $|N(F')|\leq dh$, the number of possible outcomes of the algorithm described above is at most
    $$\binom{dh}{\leq r/((\psi-1)(d-1-\varphi)}= 2^{O\left(\frac{r\log d}{\psi(d-\varphi))}\right)},$$
where the equality comes from Proposition~\ref{prop:Entropy} together with Claim \ref{claim:quotientrh2}.

Set $\hat{F}^{*}\coloneqq V(M(I))\cap F^{*}$. Notice that in each step where we add vertices to $F'$, we add at most one vertex from $V(M_{2d-1})\setminus H$, which belongs to $V(M(I))$. So, for each $u\in Q$, either $N(u)\cap V(M)=\emptyset$ or we have $d$ choices for $N(u)\cap V(M)$. Therefore, the number of ways to specify $\hat{F}^{*}$ is at most
\[
(d+1)^{r/((\psi-1)(d-\varphi))}=2^{O\left(\frac{r\log d}{\psi(d-\varphi)}\right)}.
\]
After specifying $\hat{F}^{*}$, we set $F''=F^{*}\setminus \hat{F}^{*}\subseteq H$. Since $\hat{F}^{*}\cap H\subseteq V(M)\cap H=\emptyset$ and $d_{H\setminus F^{*}}(u)< \psi-1$ for every $u\in Q$, we have $d_{F''}(u)+d_{\hat{F}^{*}}(u)>d-\psi+1$. Set $S''\coloneqq\{u\in \mathcal{L}_d: d_{F''}(u)\geq d-\psi\}$. Since $u\in Q$ implies $d_{F''}(u)\geq d-\psi$, we have $Q\subseteq S''$.

Now we obtain the desired pair $(F'',S'')$ in Step 1 such that $d_{F''}(u)\geq d-\psi$ for every $u\in S''$ and $Q\subseteq S''$, $F''\subseteq H$.

\textbf{Step 2:} Again, we use the linear ordering $\prec_2$ on $\mathcal{L}_{d}\cup \mathcal{L}_{d-1}$. Following $\prec_2$ for each vertex $w\in V(M_{2d-1})\setminus H$, if $d_{S''}(w)>\psi$, then we update $S''$ with $S''\setminus N(w)$. Note that such $w$ has exactly one neighbor in $Q$ when $w\in V(M)$, and no neighbors in $Q$ otherwise. So, when $w\in V(M)$, in order to satisfy $Q\subseteq S''$, we need to add the vertex in $N(w)\cap V(M)$ back to $S''$. Therefore, for this vertex $w\in V(M_{2d-1})\setminus H$ we have $d+1$ options: add back either one of its $d$ neighbors or none of them. In the end,  set $S\coloneqq S''$ and $F\coloneqq F''\cup \{u\in V(M_{2d-1}): d_{S}(u)>\psi\}$. Notice that for each $u\in S$ we have $d_{F}(u)\geq d_{F''}(u)\geq d-\psi$.

Since $v\in S''$ implies $d_{F''}(v)\geq d-\psi$, we have 
\[
|S''\setminus Q|\cdot (d-\psi)\leq e(S''\setminus Q,  F'')\leq e(H, V(M_{2d-1})^c\setminus Q)=r.
\]
Each time when we remove vertices from $S''$, we remove at least $\psi-1$ vertices from $S''\setminus Q$, that is, the guaranteed $\psi$ vertices from $N(w)$ minus the potential vertex from $V(M)$. Therefore, the number of times we could delete vertices from $S''\setminus Q$ is at most $r/((d-\psi)\psi)$.

Note that each vertex $w$ is chosen from $N(S'')\subseteq N^2(H)\subseteq N^{3}(Q)\subseteq N^4(F)$. Since $|F''|\leq h$, the number of ways to run Step 2 is at most
\[\binom{d^4 h}{\leq r/((\psi-1)(d-\psi))}(d+1)^{r/((\psi-1)(d-\psi))}= 2^{O\left(\frac{r\log d}{\psi(d-\psi)}\right)}.\qedhere\]

\end{proof}

\end{document}